\documentclass[12pt]{amsart}
\usepackage{amssymb,eucal,graphicx}
\usepackage{psfrag}
\usepackage{ifthen}
\usepackage{latexsym}
\usepackage[T1]{fontenc}

\usepackage[all,poly,import,color,ps,dvips]{xy}
\def\colorxy(#1){/xycolor{#1 setrgbcolor}def}

\SelectTips{cm}{}

\numberwithin{equation}{section}

\newtheorem{theorem}[equation]{Theorem}
\newtheorem{theorema}{Theorem}
\newtheorem{remarko}{Remark}

\newtheorem{corollary}[equation]{Corollary}
\newtheorem*{corollario}{Corollary}
\newtheorem{claim}[equation]{Claim}
\newtheorem{lemma}[equation]{Lemma}
\newtheorem{proposition}[equation]{Proposition}

\theoremstyle{definition}
\newtheorem{definition}[equation]{Definition}
\newtheorem{notation}[equation]{Notation}
\newtheorem{example}[equation]{Example}

\newtheorem{remark}[equation]{Remark}

\newtheorem{algo}[equation]{Normal crossing reduction algorithm}
\newtheorem{algog}[equation]{Gorenstein reduction algorithm}

\DeclareMathOperator{\mult}{mult}
\DeclareMathOperator{\emdim}{emdim}
\DeclareMathOperator{\coker}{coker}
\DeclareMathOperator{\codim}{codim}
\DeclareMathOperator{\depth}{depth}
\DeclareMathOperator{\pd}{pd}
\DeclareMathOperator{\Sing}{Sing} \DeclareMathOperator{\im}{Im}

\newcommand{\red}{\text{\rm red}}

\newcommand\Ii{{\mathcal I}}
\newcommand\Oc{{\mathcal O}}

\newcommand\m{\mathfrak m}

\newcommand\X{\mathcal X}
\newcommand\Ha{\mathcal H}
\newcommand\D{\Delta}
\newcommand\eps{\varepsilon}

\newcommand\ZZ{\mathbb{Z}}
\newcommand\CC{\mathbb{C}}
\newcommand\NN{\mathbb{N}}

\newcommand\N{\mathcal N}

\newcommand\Pp{\mathbb P}
\newcommand\PR{\mathbb{P}^r}
\newcommand\Pt{\mathbb{P}^3}

\newcommand\Ff{\mathbb{F}}
\newcommand\A{\mathbb{A}}

\newcommand{\cO}{\mathcal{O}}
\renewcommand{\P}{\mathbb{P}}
\renewcommand{\geq}{\geqslant}
\renewcommand{\leq}{\leqslant}
\renewcommand{\ge}{\geqslant}
\renewcommand{\le}{\leqslant}

\newcommand{\st}{\scriptstyle}
\newcommand{\dt}{\displaystyle}

\begin{document}

\title{On degenerations of surfaces}
\author{A. Calabri, C. Ciliberto, F. Flamini, R. Miranda}

\email{calabri@dmsa.unipd.it}
\curraddr{Dipartimento Metodi e Modelli Matematici per le Scienze Applicate, Universit\`a degli Studi di Padova\\
via Trieste 63, 35121 Padova \\Italy}

\email{cilibert@mat.uniroma2.it}
\curraddr{Dipartimento di Matematica, Universit\`a degli Studi di Roma
Tor Vergata\\ Via della Ricerca Scientifica, 1 - 00133 Roma \\Italy}

\email{flamini@mat.uniroma2.it}
\curraddr{Dipartimento di Matematica, Universit\`a degli Studi di Roma
Tor Vergata\\ Via della Ricerca Scientifica, 1 - 00133 Roma \\Italy}

\email{miranda@math.colostate.edu}
\curraddr{Department of Mathematics, 101 Weber Building,
Colorado State University
\\ Fort Collins, CO 80523-1874 \\USA}

\thanks{
The first three authors are members of M.I.U.R-G.N.S.A.G.A. at I.N.d.A.M. ``F. Severi''}

\begin{abstract}
This paper surveys and gives a uniform exposition of results contained in \cite{CCFMto}, \cite{CCFMLincei}, \cite{CCFMk2} and \cite{CCFMpg}. The subject is degenerations of surfaces, especially to unions of planes. 
More specifically, we deduce some properties of the smooth surface which is the general fibre of
the degeneration from combinatorial features of the central fibre.
In particular we show that there are strong constraints on the invariants of a
smooth surface
which degenerates to configurations of planes.

Finally we consider several examples of embedded degenerations of smooth surfaces
to unions of planes.

Our interest in these problems has been raised by a series of interesting
articles by Guido
Zappa in 1950's.
\end{abstract}

\maketitle
\tableofcontents

\section{Introduction}\label{S:1}

This paper surveys and gives a uniform exposition of results contained in \cite{CCFMto}, \cite{CCFMLincei}, \cite{CCFMk2} and \cite{CCFMpg}, where we study several properties of flat degenerations
of surfaces whose general fibre is a smooth projective algebraic surface
and whose central fibre is a reduced, connected surface $X \subset \PR$, $r \geq
3$; a very interesting case is, in particular, when $X$ is assumed to be a \emph{union of planes}.

As a first application of this approach, we shall see that there are
strong constraints on the invariants of a
smooth projective surface which degenerates to configurations of planes with
global normal crossings or other mild singularities (cf.\ \S~\ref{S:last}).

Our results include formulas on the basic
invariants of smoothable surfaces (see e.g.\ Theorems \ref{thm:4.pgbis}, \ref{thm:k2Gmain} and \ref{thm:pg}).
These formulas are useful in studying a wide range of open problems, as it
happens in the curve case, where
one considers \emph{stick curves}, i.e.\ unions of lines with only nodes as singularities. Indeed, as stick curves
are used to study moduli spaces of smooth curves and are
strictly related to fundamental problems as the \emph{Zeuthen problem}
(cf.\ \cite{Hart} and \cite{Sev}), degenerations of surfaces to unions of planes
naturally arise in several important
instances, like toric geometry (cf.\ e.g.\ \cite{BCK}, \cite{GonciuL} and \cite{Nami})
and the study of the behaviour of components of
moduli spaces of smooth surfaces and
their compactifications. For example, see \cite{Ran}, where the abelian
surface case is considered, or several papers related to the $K3$ surface case
(see, e.g., \cite{CLM}, \cite{CMT} and \cite{FM}).

Using the techniques developed here (cf.\ also the original papers \cite{CCFMk2} and \cite{CCFMpg}),
we are able to extend some results of topological
nature of Clemens-Schmid (see Theorem \ref{thm:pg} and cf.\ e.g.\ \cite{Morr})
and to prove a Miyaoka-Yau type inequality
(see Theorem \ref{thm:1ZappaBis} and Proposition \ref{prop:9chi-k2}).

We expect that degenerations of surfaces to unions of planes
will find many applications. For example, applications to secant varieties and the interpolation problem are contained in \cite{CDM}.  

It is an open problem to understand when a family of surfaces may degenerate
to a union of planes, and in some sense this is one of the most interesting
questions in the subject. The techniques we develop here in some cases allow
us to conclude that this is not possible.  When it is possible, we obtain
restrictions on the invariants which may lead to further theorems on classification. 
The case of scrolls has been treated in \cite{CCFMLincei} (cf. also Theorem 
\ref{thm:introLincei} below). Other possible applications are related, for example, to 
the problem of bounding the irregularity of surfaces in $\Pp^4$.

Further applications include the possibility of performing braid monodromy computations
(see \cite{CS}, \cite{Moi}, \cite{MoiTei}, \cite{Teicher}).
We hope that future work will include an analysis of
higher-dimensional analogues.

Our interest in degenerations to union of planes has been stimulated by a series
of papers by Guido
Zappa appeared in the 1940--50's regarding in particular: (1)~degenerations of scrolls
to unions of planes and (2)~the computation of bounds for the topological
invariants of an arbitrary smooth projective surface which
degenerates to a union of planes (see \cite{Za1, Za1b,
Za1c, Za2, Za2b, Za3, Za3b}).

In this paper we shall consider a reduced, connected, projective surface $X$
which is
a union of planes --- or
more generally a union of smooth surfaces --- whose
singularities are:
\begin{itemize}
\item in codimension one, double curves which are smooth and irreducible,
along which two surfaces meet
transversally;
\item multiple points, which are locally analytically isomorphic
to the vertex of a cone over a stick curve
with arithmetic genus either zero or one and which is projectively
normal in the projective space it spans.
\end{itemize}
These multiple points will be called \emph{Zappatic singularities}
and $X$ will be called a \emph{Zappatic surface}.
If moreover $X \subset \PR$, for some positive $r$, and if all its irreducible
components
are planes, then $X$ is called a \emph{planar Zappatic surface}.

We will mainly concentrate on the so called \emph{good Zappatic surfaces},
i.e.\ Zappatic surfaces having only Zappatic singularities whose associated
stick curve has one of the following
dual graphs (cf.\ Examples \ref{ex:tngraphs} and \ref{ex:zngraphs}, Definition
\ref{def:goodzapp},
Figures \ref{fig:zapsings} and \ref{fig:R3E3R4}):
\begin{itemize}
    \item[$R_n$:] a chain of length $n$, with $n\ge3$;
    \item[$S_n$:] a fork with $n-1$ teeth, with $n\ge4$;
    \item[$E_n$:] a cycle of order $n$, with $n\ge3$.
\end{itemize}
Let us call $R_n$-, $S_n$-, $E_n$-{\em point} the corresponding
multiple point of the Zappatic surface $X$.

We first
study some combinatorial properties of a Zappatic surface $X$ (cf.\ \S~\ref{S:3}).
We then focus on the case in which $X$ is the central fibre of a (an embedded, respectively) flat
degeneration $\X\to\D$, where $\D$ is the complex unit
disk (and where $\X \subseteq \D \times \PR$, $r \geq 3$, is a closed subscheme
of relative dimension two, respectively). In this case, we deduce some
properties of the general fibre $\X_t$, $t \neq 0$, of the degeneration from the
aforementioned properties of the central fibre $\X_0 = X$
(see \S 's \ref{S:zapdeg},
\ref{S:5}, \ref{S:pgzappdeg}, \ref{S:BI} and \ref{S:last}).

A first instance of this approach can be found in \cite{CCFMto},
where we gave some partial results on the computation of $h^0(X, \omega_X)$, when $X$ is
a Zappatic surface with global normal crossings, i.e.\ with only $E_3$-points, and
where $\omega_X$ is its dualizing sheaf (see Theorem 4.15 in \cite{CCFMto}).
In the particular case in which $X$ is smoothable, namely if
$X$ is the central fibre of a flat degeneration, we recalled that
the formula for $h^0(X,\omega_X)$ can be also deduced from the well-known
Clemens-Schmid exact sequence (cf.\ also e.g.\ \cite{Morr}).

In this paper we address three main problems.

We define the \emph{$\omega$-genus} of a projective variety $Y$ to be
\begin{equation}\label{eq:pgX}
p_{\omega}(Y) := h^0 (Y, \omega_Y),
\end{equation}
where $\omega_Y$ is the dualizing sheaf of $Y$. It is just the \emph{arithmetic} genus, if $Y$ is a reduced curve, and the \emph{geometric} genus, if $Y$ is a smooth surface.

We first extend the computation of the $\omega$-genus to the more general case
in which a good Zappatic surface $X$ --- considered as a reduced, connected surface on its own ---
has $R_n$-, $S_n$- and  $E_n$-points, for $n \geq 3$, as Zappatic singularities (cf.\ Theorem \ref{thm:4.pgbis}). When $X$ is the central fibre of a degeneration $\X\to\Delta$, we then
relate the $\omega$-genus of the central fibre to the {\em geometric genus} of the
general one; more precisely, we show that the $\omega$-genus of the fibres of a
flat degeneration of surfaces with Zappatic central fibre as above
is \emph{constant} (cf.\ Theorem \ref{thm:pg} and \cite{CCFMpg}).

As a second main result, we compute the $K^2$ of a smooth surface
which degenerates to a good Zappatic surface,
i.e.\ we compute $K^2_{\X_t}$, where $\X_t$ is the general fibre
of a degeneration $\X \to \Delta$ such that the central
fibre $\X_0$ is a good Zappatic surface (see\ Theorem \ref{thm:k2Gmain} and \cite{CCFMk2}).

We will then prove a basic inequality,
called the {\em Multiple Point Formula} (cf.\ Theorem \ref{thm:BI} and \cite{CCFMk2}),
which can be viewed as a generalization, for good Zappatic singularities,
of the well-known Triple Point Formula
(see Lemma \ref{lem:tpf} and cf.\ \cite{Frie}).

These results follow from a detailed analysis of local properties
of the total space $\X$ of the degeneration at a good Zappatic
singularity of the central fibre $X$ (cf.\ \S's ~\ref{S:MS} and ~\ref{S:resolution}).

Furthermore, we apply the computation of $K^2$ and the Multiple Point Formula
to prove several results concerning degenerations of surfaces.
Precisely, if $\chi$ and $g$ denote, respectively, the
Euler-Poincar\'e characteristic and the sectional genus of the
general fibre $\X_t$, for $t \in \D \setminus \{ 0 \}$, then (cf.\ Definition \ref{def:zappdeg}):

\begin{theorema}[cf.\ Theorem \ref{thm:1ZappaBis}] \label{thm:1}
Let $\X \to \D$ be a good, planar Zappatic degeneration, where the
central fibre $\X_0=X$ has at most $R_3$-, $E_3$-, $E_4$- and
$E_5$-points. Then
\begin{equation}\label{eq:1ZappaBis1}
K^2\leq 8\chi+1-g.
\end{equation}
Moreover, the equality holds in \eqref{eq:1ZappaBis1} if and only
if $\X_t$ is either the Veronese surface in $\Pp^5$ degenerating
to four planes with associated graph $S_4$ (i.e.\ with three
$R_3$-points, see Figure \ref{fig:zappequal1}.a), or an elliptic
scroll of degree $n \geq 5$ in $\Pp^{n-1}$ degenerating to $n$
planes with associated graph a cycle $E_n$
(see Figure \ref{fig:zappequal1}.b).

Furthermore, if $\X_t$ is a surface of general type, then$$K^2 <
8\chi -g.$$
\end{theorema}

\begin{figure}[ht]
\[
\begin{array}{cc}
\qquad
\begin{xy}
0; <35pt,0pt>: 
{\xypolygon3{~>{}}}, "1"*=0{\bullet};"0"*{\bullet}**@{-} ,
"2"*=0{\bullet};"0"**@{-} , "3"*=0{\bullet};"0"**@{-} ,
\end{xy}
\qquad & \qquad
\begin{xy}
0; <35pt,0pt>: 
{\xypolygon9{\bullet}},
\end{xy}
\qquad \\[35pt]
\text{(a)}  & \text{(b)}
\end{array}
\]
\caption{}\label{fig:zappequal1}
\end{figure}

In particular, we have:

\begin{corollario}[cf.\ Corollaries \ref{cor:1Zappabis} and \ref{cor:2Zappa}]
Let $\X$ be a good, planar Zappatic degeneration.
\begin{itemize}
\item[(a)] Assume that $\X_t$, $t \in \D \setminus \{0 \}$, is
a scroll of sectional genus $g\geq 2$. Then $\X_0 = X$ has worse
singularities than $R_3$-, $E_3$-, $E_4$- and $E_5$-points.

\item[(b)] If $\X_t$ is a minimal surface of general type and
$\X_0 = X$ has at most $R_3$-, $E_3$-, $E_4$- and $E_5$-points,
then$$g\leq 6\chi+5.$$
\end{itemize}
\end{corollario}
These improve the main results of Zappa in \cite{Za3}.

Let us describe in more detail the contents of the paper.
Section \ref{S:2} contains some basic
results on reducible curves and their dual graphs.

In Section \ref{S:3}, we give the definition of Zappatic singularities
and of (planar, good) Zappatic surfaces.
We associate to a good Zappatic surface $X$ a graph $G_X$
which encodes the configuration of the irreducible components of $X$
as well as of its Zappatic singularities  (see Definition \ref{def:dualgraph}).

Then we compute
from the combinatorial invariants of the associated graph $G_X$ some of the
invariants of $X$,
e.g.\ the Euler-Poincaré characteristic $\chi(\Oc_X)$, and
--- when $X \subset \PR$, $r \geq 3$ --- the degree
$d=\deg(X)$, the sectional genus $g$, and so on.
These computations will be frequently used in later sections, e.g.\
\S~\ref{S:last}.

In Section \ref{S:4pg} we address the
problem of computing
the $\omega$-genus of a good Zappatic surface $X$.
Precisely, we compute the cohomology of its structure sheaf,
since $p_\omega(X)=h^ 2(X,\cO_X)$, and we prove the following:

\begin{theorema} \label{thm:intro0} (cf.\ Theorem \ref{thm:4.pgbis})
Let $X=\bigcup_{i=1}^v X_i$ be a good Zappatic surface and let $G_{X}$ be its
associated graph (cf.\ Definition \ref{def:dualgraph}). Consider the natural map
$$\Phi_X:\bigoplus_{i=1}^v H^1(X_i,\Oc_{X_i}) \to \bigoplus_{1 \leq i< j \leq v}
H^1(C_{ij},\Oc_{C_{ij}}),$$where $C_{ij}=X_i\cap X_j$ if $X_i$ and $X_j$ meet along a curve,
or $C_{ij}=\emptyset$ otherwise
(cf.\ Definition \ref{rem:indices}). Then:
\begin{equation}\label{eq:intro}
p_{\omega}(X) = h^ 2(G_X, \CC)+\sum_{i=1}^v p_g(X_i)+
\dim({\rm coker}(\Phi_X)). \end{equation}
\end{theorema}

In particular, we have:

\begin{corollario}
Let $X$ be a good
planar Zappatic surface. Then,
\[
p_{\omega}(X)= b_2(G_X).
\]
\end{corollario}

\begin{remarko}\label{rem:topinv}
\normalfont{
It is well-known that, for smooth surfaces $S$, the  geometric genus $p_g(S)$ is a topological invariant
of $S$. From Formula \eqref{eq:intro}, it follows that also the $\omega$-genus $p_{\omega}(X)$ is a topological invariant
of any good Zappatic surface $X$.
}
\end{remarko}

In order to prove the above results, we exploit the
natural injective resolution of the sheaf $\cO_X$ in terms
of the structure sheaves of the irreducible components of $X$ and of its
singular locus.  An alternative, and in some sense dual, approach is via
the interpretation of the global sections of $\omega_X$ as collections of
meromorphic 2-forms on the irreducible components of $X$, having poles
along the double curves of $X$ with suitable matching conditions. This
interpretation makes it possible, in principle, to compute $h^
0(X,\omega_X)$ by computing the number of such independent collections of
forms. This is the viewpoint  taken in \cite {CCFMto}, where we discussed
only the normal crossings case. However, the approach taken here leads
more quickly and neatly to our result.

In Section \ref{S:zapdeg} we give the definition of
Zappatic degenerations of surfaces and we recall some properties of smooth surfaces which degenerate to
Zappatic ones.

In Section \ref{S:MS} we recall the notions of
{\em minimal singularity} and {\em quasi-minimal singularity},
which are needed to study the singularities of the total space $\X$ of a degeneration of surfaces
at a good Zappatic singularity of its central fibre $\X_0 = X$
(cf.\ also \cite{Kollar} and \cite{KolSB}).

Section \ref{S:resolution} is devoted to studying (partial and total)
resolutions of the singularities that
the total space $\X$ of a degeneration of surfaces
at a good Zappatic singularity of its central fibre.

The local analysis of minimal and quasi-minimal
singularities of $\X$ is fundamental in  \S\ \ref{S:5}, where we
compute $K_{\X_t}^2$,
for $t \in \D \setminus \{0\}$, when $\X_t$ is the general fibre of a degeneration
such that the central fibre is a good Zappatic surface.
More precisely, we prove the following main result (see Theorem
\ref{thm:k2Gmain} and \cite{CCFMk2}):

\begin{theorema}\label{thm:intro1}
Let $\X\to\D$ be a degeneration of surfaces whose central fibre
is a good Zappatic surface $X=\X_0=\bigcup_{i=1}^v X_i$.
Let $C_{ij} := X_i \cap X_j$ be a smooth
(possibly reducible) curve of the double locus of $X$, considered
as a curve on $X_i$,  and let $g_{ij}$ be its geometric genus, $1
\leq i \neq j \leq v$.
Let $v$ and $e$ be the number of vertices and edges of the graph
$G_X$ associated to $X$. Let $f_n$, $r_n$, $s_n$ be the number of
$E_n$-, $R_n$-, $S_n$-points of $X$, respectively.
If $K^2 := K^2_{\X_t}$, for $t\ne0$, then:
\begin{equation}
K^2  =  \sum_{i=1}^v \left( K_{X_i}^2 + \sum_{j \ne i} (4 g_{ij} -
C_{ij}^2) \right)-8e+\sum_{n \ge 3} 2n f_n+r_3+ k
\end{equation}
where $k$ depends only on the presence of $R_n$- and $S_n$-points,
for $n\ge 4$, and precisely:
\begin{equation}\label{eq:introc}
\sum_{n \geq 4} (n-2)(r_n +s_n) \le k
 \le \sum_{n \ge 4} \left( (2n-5) r_n + \binom{n-1}{2} s_n \right).
\end{equation}
\end{theorema}

In the case that the central fibre is also planar, we have the
following:

\begin{corollario}[cf.\ Corollary \ref{cor:k2Gmain}]
Let $\X\to\D$ be an embedded degeneration of surfaces whose central fibre
is a good, planar Zappatic surface
$X=\X_0=\bigcup_{i=1}^v \Pi_i$. Then:
\begin{equation}
K^2 = 9v-10e+\sum_{n \ge 3} 2n f_n+r_3+ k
\end{equation}
where $k$ is as in \eqref{eq:introc} and depends only on the
presence of $R_n$- and $S_n$-points, for $n\ge 4$.
\end{corollario}

The inequalities in the theorem and the corollary above reflect
deep geometric properties of the degeneration. For example, if $\X
\to \D$ is a degeneration with central fibre
$X$ a Zappatic surface which is the union of four planes
having only a $R_4$-point, Theorem \ref{thm:intro1} states
that $8 \leq K^2 \leq 9$. The two values of $K^2$ correspond to
the fact that $X$, which is the cone over a stick curve $C_{R_4}$ (cf.\ Example
\ref{ex:tngraphs}), can be smoothed either to the Veronese surface,
which has $K^2=9$, or to a rational normal quartic scroll in
$\Pp^5$, which has $K^2=8$ (cf.\ Remark \ref{rem:R4expli}). This
in turn corresponds to  different local structures of the total
space of the degeneration at the $R_4$-point. Moreover, the local
deformation space of a $R_4$-point is reducible.

Section \ref{S:BI} is devoted to the {\em Multiple Point Formula}
\eqref{eq:B1} below (cf.\ Definition \ref{def:zappdeg}, see Theorem \ref{thm:BI} and \cite{CCFMk2}):

\begin{theorema}\label{thm:intro2}
Let $X$ be a good Zappatic surface which is the central fibre of a
good Zappatic degeneration $\X \to \D$. Let $\gamma=X_1\cap X_2$
be the intersection of two irreducible components $X_1$, $X_2$ of
$X$. Denote by $f_n(\gamma)$ [$r_n(\gamma)$ and $s_n(\gamma)$,
respectively] the number of $E_n$-points [$R_n$-points and
$S_n$-points, respectively] of $X$ along $\gamma$. Denote by
$d_\gamma$ the number of double points of the total space $\X$
along $\gamma$, off the Zappatic singularities of $X$. Then:
\begin{equation}\label{eq:B1}
\deg(\N_{\gamma|X_1}) + \deg(\N_{\gamma|X_2}) + f_3(\gamma) -
r_3(\gamma) -\sum_{n\ge4} (r_n(\gamma)+ s_n(\gamma) +f_n(\gamma))
\ge d_\gamma \ge 0.
\end{equation}
In particular, if $X$ is also planar, then:
\begin{equation}\label{eq:B2}
2 + f_3(\gamma) - r_3(\gamma) -\sum_{n\ge4} (r_n(\gamma) +
s_n(\gamma) + f_n(\gamma)) \ge d_\gamma \ge 0.
\end{equation}
Furthermore, if $d_\X$ denotes the total number of double points
of $\X$, off the Zappatic singularities of $X$, then:
\begin{equation}\label{eq:B3}
2e+3f_3 - 2 r_3 -\sum_{n\ge4} nf_n -\sum_{n\ge 4} (n-1) (s_n +
r_n) \ge d_\X \ge 0.
\end{equation}
\end{theorema}

In \S\ \ref{S:last} we apply Theorem \ref{thm:intro1} and \ref{thm:intro2} above
to prove several generalizations of statements given by Zappa. For example we show that
worse singularities than normal crossings are needed in order to degenerate as
many surfaces as possible to unions of planes.

In Section \ref{S:pgzappdeg}, we apply the result in \S~\ref{S:4pg}
on the $\omega$-genus (i.e.\ Theorem \ref{thm:intro0} above) to the case of
$X$ a smoothable, good Zappatic surface, namely
$X=\X_0$ is the central fibre of a flat degeneration $\X\to\D$ of surfaces, where $\D$ is the spectrum
of a DVR (or equivalently the complex unit disk)
and each fibre $\X_{t}=\pi^{-1}(t)$, $0 \neq t \in \D$, is smooth.

Precisely, with the same hypotheses of Theorem \ref{thm:intro0}, we prove:

\begin{theorema}\label{thm:intro3} (cf.\ Theorem \ref{thm:pg})
Let $\X\to\D$ be a flat degeneration of surfaces parametrized by a disk,
such that the central fibre $\X_0=X$ is good Zappatic and each fibre $\X_t$,
$t\ne0$, is smooth. Then, for any $t\neq0$, one has:
\begin{equation}\label{eq:intro3}
p_g(\X_t) = p_{\omega}(X).
\end{equation}In particular, the $\omega$-genus of the fibres of $\X \to \D$ is
\emph{constant}.
\end{theorema}

\begin{remarko}\label{rem:smussab}
\normalfont{Recall that, when $X$ has only $E_3$-points as Zappatic singularities and it is smoothable,
with smooth total space $\X$, (i.e.\ $X$ is the central fibre of a
{\em semistable degeneration}), it is well--known that \eqref{eq:intro3} holds.
This has been proved via the {\em Clemens-Schmid exact sequence} approach, which
relates the mixed Hodge theory of the central fibre $X$ to that of the general one
$\X_t$, $t \in \D \setminus \{ 0 \}$, by means of the monodromy of the total space $\X$
(cf.\ e.g.\ \cite{Morr} for details).

We remark that Theorems \ref{thm:intro0} and \ref{thm:intro3} above not only show
that \eqref{eq:intro3} more generally
holds for a smoothable good Zappatic surface, i.e.\ with $R_n$-, $S_n$- and $E_n$-points,
$n\ge3$, as Zappatic singularities and with (possibly) singular total space $\X$,
but mainly they extend the Clemens-Schmid approach
since the computation of $p_{\omega}(X)$ is independent of the fact that $X$ is
the central fibre of a degeneration.

To prove Theorem \ref{thm:intro3}, we use the construction performed in \S\ \ref{S:resolution}
of a normal crossing reduction $\bar\pi:\bar\X\to\D$ of $\pi$, i.e.\ $\bar\X\to\X$ is a resolution of singularities of $\X$ and the support of its central fibre $\bar\X_0$ has global normal crossings
(cf.\ Remark \ref{def:degen2}). Then we apply the results in Chapter II of \cite{Kempf}
in order to get a semistable reduction $\tilde\pi:\tilde\X\to\D$ of $\pi$.
This enables us to deduce the topological properties of the fibres of $\tilde{\X}$ from those of $X$,
with the assistance of the Clemens-Schmid exact sequence
(cf.\ e.g.\ \cite{Morr}).
}
\end{remarko}

In the last section we  exhibit several examples of  degenerations of smooth surfaces
to good Zappatic ones, some of them contained also in \cite{CCFMLincei},
in particular with the central fibre having only $R_3$-, $E_n$-, $3\le n\le 6$,
points.

We conclude the paper with Appendix \ref{S:pCMpG}, where we collect
several definitions and well-known results concerning  connections between
commutative homological Algebra and projective Geometry.

\medskip
{\it Acknowledgments.} The authors would like
to thank L. Badescu, A. Beauville and J. Koll\'ar,
for fundamental discussions and references.

\section{Reducible curves and associated graphs}\label{S:2}

Let $C$ be a projective curve and let $C_i$, $i=1,\ldots,n,$ be its
irreducible components. We will assume that:

\begin{itemize}
\item $C$ is connected and reduced;
\item $C$ has at most nodes as
singularities;
\item the curves $C_i, i=1,\ldots,n,$ are smooth.\end{itemize}

If two components $C_i, \; C_j, \; i<j,$ intersect at $m_{ij}$ points, we
will denote by $P_{ij}^h, \; h=1,\ldots,m_{ij}$, the corresponding nodes of $C$.

We can associate to this situation a simple (i.e.\ with no loops),
weighted connected graph $G_C$, with vertex $v_i$ weighted by the genus $g_i$ of $C_i$:
\begin{itemize}

\item whose
vertices  $v_1,\ldots,v_n,$ correspond to the components $C_1,\ldots, \;C_n$;

\item whose edges $\eta^h_{ij}$, $i<j, \; h=1,\ldots,m_{ij}$, joining the
vertices $v_i$ and $v_j$, correspond to the nodes $P_{ij}^h$ of
$C$.\end{itemize}

We will assume the graph to be {\em lexicographically oriented}, i.e.\ each edge is assumed to be oriented from the vertex with lower index to
the one with higher.

We will use the following notation:

\begin{itemize}
\item $v$ is the number of vertices of $G_C$, i.e.\ $v=n$;
\item $e$ is the number of edges of $G_C$;
\item $\chi(G_C)=v-e$ is the Euler-Poincar\'e characteristic of $G_C$;
\item $h^1(G_C)=1-\chi(G_C)$ is the first Betti number of $G_C$.
\end
{itemize}

Notice that conversely, given any simple, connected, weighted (oriented)
graph $G$, there is some curve $C$ such that $G=G_C$.

One has the following basic result:

\begin{theorem}\label{thm:pacurves} (cf. \cite[Theorem 2.1]{CCFMto}) 
In the above situation
\begin{equation}\label{eq:chicurves}
\chi(\Oc_C)=\chi(G_C) - \sum_{i=1}^v g_i = v - e - \sum_{i=1}^v g_i.
\end{equation}
\end {theorem}

\begin{proof}  Let $\nu: \tilde{C} \to C$ be the normalization morphism; this
defines
the exact sequence of sheaves on $C$:
\begin{equation}\label{eq:4gen}
0 \to \Oc_C \to \nu_* (\Oc_{\tilde{C}}) \to \underline{\tau} \to 0,
\end{equation}
where $\underline{\tau} $ is a skyscraper sheaf supported on $\Sing(C)$.
Since the singularities of $C$ are only nodes, one easily determines
$H^0(C, \underline{\tau} ) \cong \CC^e$. Therefore, by the exact sequence
\eqref{eq:4gen}, one gets
$$\chi(\Oc_C) = \chi(\nu_* (\Oc_{\tilde{C}})) -e.$$By the Leray isomorphism and
by the
fact that $\nu$ is finite, one has
$\chi(\nu_* (\Oc_{\tilde{C}})) = \chi(\Oc_{\tilde{C}})$. Since $\tilde{C}$ is
a disjoint union of the $v=n$ irreducible components of $C$, one has
$\chi(\Oc_{\tilde{C}}) = v - \sum_{i=1}^v g_i$, which proves
\eqref{eq:chicurves}.
(Cf.\ also \cite{Bar} for another proof.)
\end{proof}

We remark that Formula \eqref{eq:chicurves} is equivalent to (cf.\ Proposition \ref{prop:ghypsect}):
\begin{equation} \label{eq:genuscurves}
p_a(C)= h^1(G_C)+ \sum_{i=1}^v g_i.
\end{equation}

Notice that $C$ is Gorenstein, i.e.\ the dualizing sheaf
$\omega_C$ is invertible. We define the $\omega$-{\it genus} of $C$
to be
\begin{equation}\label{eq:ggcurves}
p_{\omega}(C):= h^0(C,\omega_C).
\end{equation}Observe that, when $C$ is smooth, the $\omega$-genus coincides
with the geometric genus of $C$.

In general, by the Riemann-Roch theorem, one has
\begin{equation}\label{eq:genuscurves2}
p_{\omega}(C)=p_a(C)=h^1(G_C)+\sum_{i=1}^v g_i = e - v + 1 +\sum_{i=1}^v g_i.
\end{equation}

If we have a flat family ${\mathcal C}\to \D$ over a disc $\D$ with general
fibre ${\mathcal C}_t$ smooth and irreducible of genus $g$ and special fibre
${\mathcal C}_0=C$,
then we can combinatorially compute $g$ via the formula:
$$
g=p_a(C)=h^1(G_C)+\sum_{i=1}^vg_i.
$$

Often we will consider $C$ as above embedded in a projective space $\PR$.
In this situation each curve $C_i$ will have a certain degree $d_i$, so that
the graph $G_C$ can be considered as {\em double weighted}, by attributing to
each vertex the pair of weights $(g_i,d_i)$. Moreover one can attribute to the
graph a further marking number, i.e.\ $r$ the embedding dimension of $C$.

The total degree of $C$ is
$$
d=\sum_{i=1}^v d_i
$$
which is also invariant by flat degeneration.

More often we will consider the case in which each curve $C_i$ is a line. The
corresponding curve $C$ is called a {\it stick curve}. In this case the
double weighting is
$(0,1)$ for each vertex, and it will be omitted if no confusion arises.

\begin{figure}[ht]
\[
\begin{xy}
0; <35pt,0pt>: 
(0,1)*=0{\bullet};(0,0)*=0{\bullet}**@{-};
(1,0)*=0{\bullet}**@{-}="a";(1,1)*=0{\bullet}**@{-};(0,1)**@{-};"a"**@{-} ,
\end{xy}
\]
\caption{Dual graph of an ``impossible'' stick curve.}\label{fig:nostick}
\end{figure}
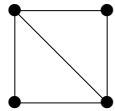

It should be stressed  that it is not true that for any simple,
connected, double weighted graph $G$  there is a curve $C$ in a projective
space such that $G_C=G$. For example there is no stick curve corresponding
to the graph of Figure \ref{fig:nostick}.

We now give two examples of stick curves which will be frequently used in this
paper.

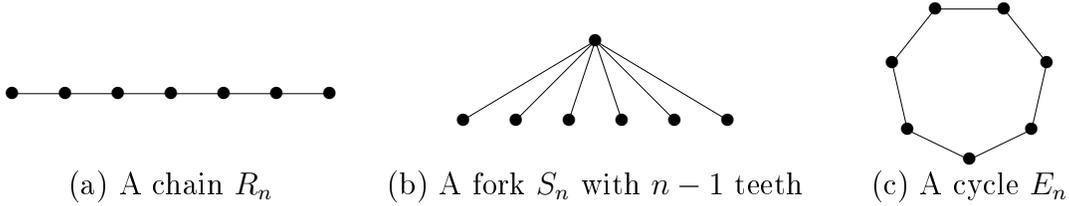
\begin{figure}[ht]
\[
\begin{array}{ccc}
\quad
\raisebox{10pt}{$ %
\begin{xy}
0; <20pt,0pt>: 
(0,0)*=0{\bullet};(1,0)*=0{\bullet}**@{-};(2,0)*=0{\bullet}**@{-};
(3,0)*=0{\bullet}**@{-};
(4,0)*=0{\bullet}**@{-};(5,0)*=0{\bullet}**@{-};(6,0)*=0{\bullet}**@{-}
,
\end{xy}
$}
\quad & \qquad
\begin{xy}
0; <20pt,0pt>: 
(0,1.5)*=0{\bullet}="u" ,
(-2.5,0)*=0{\bullet};"u"**@{-} ,
(-1.5,0)*=0{\bullet};"u"**@{-} ,
(-0.5,0)*=0{\bullet};"u"**@{-} ,
(2.5,0)*=0{\bullet};"u"**@{-} ,
(1.5,0)*=0{\bullet};"u"**@{-} ,
(0.5,0)*=0{\bullet};"u"**@{-} ,
\end{xy}
\qquad & \qquad
\raisebox{45pt}{$ %
\begin{xy}
0; <30pt,0pt>: <0pt,-30pt>::
{\xypolygon7{\bullet}}
\end{xy}
$}
\qquad
\\[4mm]
\text{(a) A chain $R_n$} & \text{(b) A fork $S_n$ with $n-1$ teeth} & \text{(c)
A cycle $E_n$}
\end{array}
\]
\caption{Examples of dual graphs.}\label{fig:zapsings}
\end{figure}

\begin{example}\label{ex:tngraphs}
Let $T_n$ be any connected tree with $n \geq 3$ vertices.
This corresponds to a non-degenerate stick curve of degree $n$ in
$\P^n$, which we denote by $C_{T_n}$.
Indeed one can check that, taking a general point $p_i$
on each component of $C_{T_n}$, the line bundle $\Oc_{C_{T_n}} (p_1+\cdots+p_n)$
is very ample.
Of course $C_{T_n}$ has arithmetic genus $0$ and is a
flat limit of rational normal curves in $\P^n$.

We will often consider two particular kinds of trees $T_n$:
a chain $R_n$ of length $n$ and the fork $S_n$ with $n-1$ teeth,
i.e.\ a tree consisting of $n-1$ vertices
joining a further vertex (see Figures \ref{fig:zapsings}.(a) and (b)).
The curve $C_{R_n}$ is the union of $n$ lines $l_1,
l_2, \ldots, l_n$ spanning $\P^n$, such that $l_i\cap l_j=\emptyset$ if and only
if $1<|i-j|$.
The curve $C_{S_n}$ is the union of $n$
lines $l_1, l_2, \ldots, l_n$ spanning $\P^n$, such that
$l_1,\ldots,l_{n-1}$ all intersect $l_n$ at distinct points (see Figure
\ref{fig:zapsings2}).
\end{example}


\begin{example}\label{ex:zngraphs}
Let $Z_n$ be any  simple, connected graph  with $n\geq 3$ vertices
and $h^1(Z_n, \CC)=1$. This corresponds to an arithmetically normal stick
curve of degree $n$ in $\P^{n-1}$, which we denote by $C_{Z_n}$
(as in Example \ref{ex:tngraphs}).
The curve $C_{Z_n}$ has arithmetic genus $1$ and it is a flat
limit of elliptic normal curves in $\P^{n-1}$.

We will often consider the particular case of a cycle $E_n$ of order $n$
(see Figure \ref{fig:zapsings}.(c)).
The curve $C_{E_n}$ is the union of $n$ lines $l_1, l_2, \ldots, l_n$
spanning $\P^{n-1}$, such that $l_i\cap l_j=\emptyset$ if and only if $1<|i-
j|<n-1$
(see Figure \ref{fig:zapsings2}).

We remark that $C_{E_n}$ is projectively Gorenstein
(i.e.\ it is projectively Cohen-Macaulay and sub-canonical, cf.\ Proposition \ref{prop:4} in
Appendix \ref{S:pCMpG}),
because $\omega_{C_{E_n}}$ is trivial, since there is an
everywhere non-zero global section of $\omega_{C_{E_n}}$, given by
the meromorphic 1-form on each component with residues 1 and $-1$
at the nodes (in a suitable order).

All the other $C_{Z_n}$'s, instead, are not Gorenstein
because $\omega_{C_{Z_n}}$, although of degree zero, is not
trivial. Indeed a graph $Z_n$, different from $E_n$, certainly has
a vertex with valence 1. This corresponds to a line $l$ such that
$\omega_{C_{Z_n}}\otimes\Oc_l$ is not trivial.
\end{example}

\vspace{-.7cm}
\begin{figure}[ht]
\[
\begin{array}{ccc}
\begin{xy}
0; <15pt,0pt>: 
(1.6,1.6)*=0{\bullet},
(2.8,0.4)*=0{\bullet},
(4.0,1.6)*=0{\bullet},
(5.2,0.4)*=0{\bullet},
(6.4,1.6)*=0{\bullet},
(7.6,0.4)*=0{\bullet},
(0,0);(2,2)**@{-} ,
(1.2,2);(3.2,0)**@{-} ,
(2.4,0);(4.4,2)**@{-} ,
(3.6,2);(5.6,0)**@{-} ,
(4.8,0);(6.8,2)**@{-} ,
(6.0,2);(8.0,0)**@{-} ,
(7.2,0);(9.2,2)**@{-} ,
\end{xy}
&
\begin{xy}
0; <15pt,0pt>: 
(0,1.5)*=0{\bullet},
(1.5,1.5)*=0{\bullet},
(3,1.5)*=0{\bullet},
(4.5,1.5)*=0{\bullet},
(6,1.5)*=0{\bullet},
(7.5,1.5)*=0{\bullet},
(-1,1.5);(8.5,1.5)**@{-} ,
(0,0);(0,2)**@{-} ,
(1.5,0);(1.5,2)**@{-} ,
(3,0);(3,2)**@{-} ,
(4.5,0);(4.5,2)**@{-} ,
(6.0,0);(6.0,2)**@{-} ,
(7.5,0);(7.5,2)**@{-} ,
\end{xy}
&
\raisebox{20pt}{$ %
\begin{xy}
0; <35pt,0pt>: 
0,
{\xypolygon7"B"{~:{(-1,0):}~*{\bullet}~>{}}} ,
{\xypolygon7"D"{~={285}~:{(1.2,0):}~>{}}} ,
{\xypolygon7"E"{~={255}~:{(1.2,0):}~>{}}} ,
"E7";"D1"**@{-} ,
"E6";"D7"**@{-} ,
"E5";"D6"**@{-} ,
"E4";"D5"**@{-} ,
"E3";"D4"**@{-} ,
"E2";"D3"**@{-} ,
"E1";"D2"**@{-} ,
\end{xy}
$}
\\[6mm]
\text{$C_{R_n}$: a chain of $n$ lines,}
 & \text{$C_{S_n}$: a comb with $n-1$ teeth,}
 & \text{$C_{E_n}$: a cycle of $n$ lines.}
\end{array}
\]
\caption{Examples of stick curves.}\label{fig:zapsings2}
\end{figure}

\section{Zappatic surfaces and associated graphs}\label{S:3}

We will now give a parallel development, for surfaces,
to the case of curves recalled in the previous
section. Before doing this, we need to introduce the singularities we will
allow (cf. \cite[\S\;3]{CCFMto}).

\begin{definition}[Zappatic singularity] \label{def:zappaticsing}
Let $X$ be a surface and let $x\in X$ be a point. We will say that $x$ is a
{\it Zappatic singularity} for $X$ if $(X,x)$ is
locally analytically isomorphic to a pair $(Y,y)$ where $Y$ is the cone over
either a curve $C_{T_n}$ or a curve $C_{Z_n}, n\geq 3$, and $y$ is the vertex of
the cone.
Accordingly we will say that $x$ is either a $T_n$- or a $Z_n$-{\em point} for
$X$.
\end{definition}

Observe that either $T_n$- or $Z_n$-points are not classified by $n$, unless
$n=3$.

We will consider the following situation.

\begin{definition}[Zappatic surface] \label{def:Zappsurf}
Let $X$ be a projective surface with its irreducible components
$X_1,\ldots,X_v$. We will assume that $X$ has the following properties:

\begin{itemize}
\item $X$ is reduced and connected in codimension one;

\item $X_1,\ldots, \; X_v$ are smooth;

\item the singularities in codimension one of $X$ are at most
double curves which are smooth and irreducible along which two surfaces
meet transversally;

\item the further singularities of $X$ are Zappatic singularities.
\end {itemize}

A surface like $X$ will be called a {\it Zappatic surface}.
If moreover $X$ is embedded in a projective space $\P ^r$ and
all of its irreducible components are planes,
we will say that $X$ is a {\it planar Zappatic surface}. In this case,
the irreducible components
of $X$ will sometimes be denoted by $\Pi_i$ instead of $X_i$, $ 1 \leq i \leq
v$.
\end{definition}

\begin{notation}\label{not:Cij}
Let $X$ be a Zappatic surface.
Let us denote by:
\begin{itemize}
\item $X_i$: an irreducible component of $X$, $ 1 \leq i \leq v$;
\item $C_{ij}:=X_i \cap X_j$, $1 \leq i \neq j \leq v$, if $X_i$ and $X_j$ meet
along a curve,
otherwise set $C_{ij}=\emptyset$. We assume that each
$C_{ij}$ is smooth but not necessarily irreducible;
\item $g_{ij}:$ the geometric
genus of $C_{ij}$, $1 \leq i \neq j \leq v$; i.e.\ $g_{ij}$ is the sum of the
geometric genera of the irreducible (equiv., connected) components of $C_{ij}$;
\item $C:=\Sing(X)=\cup_{i<j} C_{ij}$: the union of all the double curves of $X$;
\item $\Sigma_{ijk}: = X_i \cap X_j \cap X_k$, $1 \leq i \neq j \neq k \leq v$,
if  $X_i \cap X_j \cap X_k\ne\emptyset$, otherwise $\Sigma_{ijk} = \emptyset$;
\item $m_{ijk}:$ the cardinality of the set $\Sigma_{ijk}$;
\item $P_{ijk}^h:$ the Zappatic singular point belonging to $\Sigma_{ijk}$, for
$h=1,\ldots,m_{ijk}$.
\end{itemize}
Furthermore, if $X\subset\PR$, for some $r$, we denote by
\begin{itemize}
\item $d=\deg(X):$ the degree of $X$;
\item $d_i=\deg(X_i):$ the degree of $X_i$, $i \leq i \leq v$;
\item $c_{ij}=\deg(C_{ij})$: the degree of $C_{ij}$, $1 \leq i \neq j \leq v$;
\item $D:$ a general hyperplane section of $X$;
\item $g:$ the arithmetic genus of $D$;
\item $D_i:$ the (smooth) irreducible component of $D$ lying in $X_i$,
which is a general hyperplane section of $X_i$, $1 \leq i  \leq v$;
\item $g_i:$ the genus of $D_i$, $1 \leq i \leq v$.
\end{itemize}

Notice that if $X$ is a planar Zappatic surface,
then each $C_{ij}$, when not empty, is a line and each
non-empty set $\Sigma_{ijk}$ is
a singleton.
\end{notation}

\begin{remark}\label{rem:pg}
Observe that a Zappatic surface $X$ is  Cohen-Macaulay.
More precisely, $X$ has global normal crossings except at points $T_n$, $n \geq
3$,
and $Z_m$, $m \geq 4$.
Thus the dualizing sheaf $\omega_X$ is well-defined.
If $X$ has only $E_n$-points as Zappatic singularities,
then $X$ is Gorenstein,
hence $\omega_X$ is an invertible sheaf.
\end{remark}

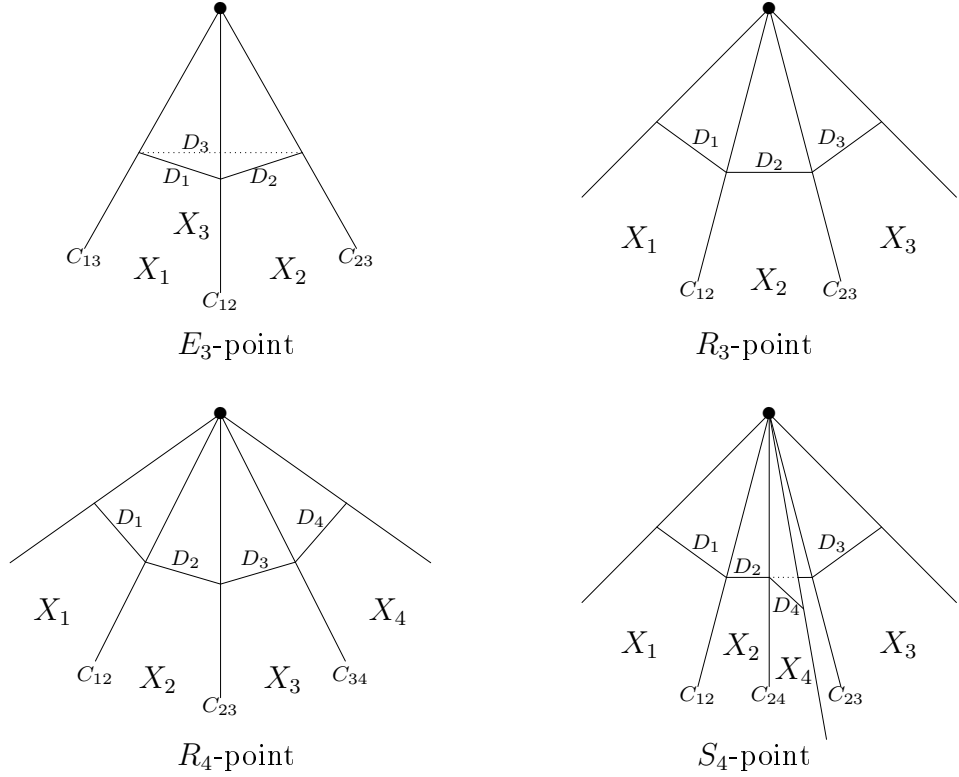
\begin{figure}[ht]
\vspace{-24mm}
\[
\begin{array}{cc}
\begin{xy}
<50pt,0pt>:
0; {\xypolygon10{~={18}~:{(1.75,0):}~>{}}},
(0,.4)="a"*=0{\bullet};"7"**@{-}?(.6)*{}="b" ,
"a";"8"**@{-}?(.6)*{}="c" ,
"a";"9"**@{-}?(.6)*{}="d" ,
"b";"c"**@{-}?+<-1pt,-1pt>*!U{\st D_1} ,
"c";"d"**@{-}?+<1pt,-1pt>*!U{\st D_2} ,
"b";"d"**@{.}?(.35)*!D{\st D_3} ,
"7";"8"**@{}?*{X_1} ,
"8";"9"**@{}?*{X_2} ,
"7";"9"**@{}?*+!DR{X_3} ,
"7"*!U{\st C_{13}} ,
"8"*!U{\st C_{12}} ,
"9"*!U{\st C_{23}} ,
\end{xy}
\quad & \quad
\begin{xy}
<50pt,0pt>:
0; {\xypolygon10{~:{(1.75,0):}~>{}}},
(0,.4)="a"*=0{\bullet};"7"**@{-}?(.6)*{}="b" ,
"a";"8"**@{-}?(.6)*{}="c" ,
"a";"9"**@{-}?(.6)*{}="d" ,
"a";"10"**@{-}?(.6)*{}="e" ,
"b";"c"**@{-}?*!DL{\st D_1} ,
"c";"d"**@{-}?+<0pt,4pt>*{\st D_2} ,
"d";"e"**@{-}?*!DR{\st D_3} ,
"7";"8"**@{}?*{X_1} ,
"8";"9"**@{}?*{X_2} ,
"9";"10"**@{}?*{X_3} ,
"8"*!U{\st C_{12}} ,
"9"*!U{\st C_{23}} ,
\end{xy}
\\[33mm]
E_3 \text{-point} & R_3 \text{-point}
\\[-16mm]
\begin{xy}
<50pt,0pt>:
0; {\xypolygon11{~={8.1818}~:{(1.75,0):}~>{}}},
(0,.4)="a"*=0{\bullet};"7"**@{-}?(.6)*{}="b" ,
"a";"8"**@{-}?(.6)*{}="c" ,
"a";"9"**@{-}?(.6)*{}="d" ,
"a";"10"**@{-}?(.6)*{}="e" ,
"a";"11"**@{-}?(.6)*{}="f" ,
"b";"c"**@{-}?+<4pt,2pt>*!D{\st D_1} ,
"c";"d"**@{-}?+<1pt,2pt>*!D{\st D_2} ,
"d";"e"**@{-}?+<-1pt,2pt>*!D{\st D_3} ,
"e";"f"**@{-}?+<-4pt,2pt>*!D{\st D_4} ,
"7";"8"**@{}?*{X_1} ,
"8";"9"**@{}?*{X_2} ,
"9";"10"**@{}?*{X_3} ,
"10";"11"**@{}?*{X_4} ,
"8"+<0pt,-2pt>*!U{\st C_{12}} ,
"9"*!U{\st C_{23}} ,
"10"+<2pt,-2pt>*!U{\st C_{34}} ,
\end{xy}
\quad & \quad
\begin{xy}
<50pt,0pt>:
0; {\xypolygon10{~:{(1.75,0):}~>{}}},
(0,.4)="a"*=0{\bullet};"7"**@{-}?(.6)="b" ,
"a";"8"**@{-}?(.6)="c" ,
"a";"9"**@{-}?(.6)="d" ,
"8";"9"**@{}?="m" ,
"a";"10"**@{-}?(.6)="e" ,
"a";"m"**@{-} ?!{"c";"d"}="n" ,
"b";"c"**@{-}?*!DL{\st D_1} ,
"c";"n"**@{-}?+<0pt,4pt>*{\st D_2} ,
"d";"e"**@{-}?*!DR{\st D_3} ,
"7";"8"**@{}?*{X_1} ,
"8";"m"**@{}?+<3pt,10pt>*!D{X_2} ,
"9";"10"**@{}?*{X_3} ,
"8"*!U{\st C_{12}} ,
"9"+<2pt,0pt>*!U{\st C_{23}} ,
"m";"9"**@{}?(.8)-(0,0.4)="p" ,
"a";"p"**@{-} ?(.6)="q" ,
"m";"9"**@{}?(.35)+<0pt,1pt>*!D{\dt X_4} ,
"m"*!U{\st C_{24}} ,
"n";"q"**@{-}?-<0pt,5pt>*{\st D_4} ,
"n";"d"**@{} ?!{"a";"p"}="r" ,
"n";"r"**@{.} ,
"r";"d"**@{-} ,
\end{xy}
\\[4mm]
R_4 \text{-point} & S_4 \text{-point}
\end{array}
\]
\caption{Examples of good Zappatic singularities.}
\label{fig:R3E3R4}
\end{figure}

\begin{definition}[Good Zappatic surface] \label{def:goodzapp}
The \emph{good Zappatic singularities} are the
\begin{itemize}
\item $R_n$-points, for $n\ge3$,
\item $S_n$-points, for $n\ge4$,
\item $E_n$-points, for $n\ge3$,
\end{itemize}
which are the Zappatic singularities whose associated stick
curves are respectively $C_{R_n}$, $C_{S_n}$, $C_{E_n}$
(see Examples \ref{ex:tngraphs} and \ref{ex:zngraphs},
Figures \ref{fig:zapsings}, \ref{fig:zapsings2} and \ref{fig:R3E3R4}).

A \emph{good Zappatic surface} is a Zappatic surface
with only good Zappatic singularities.
\end{definition}

To a good Zappatic surface $X$ we can associate an oriented
complex $G_X$,
which we will also call the {\it associated graph} to $X$.

\begin{figure}[ht]
\[
\begin{array}{cccc}
\,\,\,
\begin{xy}
0; <30pt,0pt>: 
(0,1)="a"*=0{\bullet};(1,0)*=0{\bullet}**@{-}="b";(2,1)*=0{\bullet}**@{-}="c";
"a"**@{--} ,
"a"*+!DR{\st v_1} ,
"b"*+!U{\st v_2} ,
"c"*+!DL{\st v_3} ,
\end{xy}
\,\,\, & \,\,\,
\begin{xy}
\xyimport(2,1)(0,0){\rotatebox{45}{\begin{xy} 0; <30pt,0pt>: 
(0,1.4142).(1.4142,0)*[!\colorxy(0.875 0.875 0.875)]\frm{*} , \end{xy}}}
0; <30pt,0pt>: 
(0,1)="a".(2,2)*[white]\frm{*} ,
"a"*=0{\bullet};(1,0)="b"*=0{\bullet}**@{-};
(2,1)*=0{\bullet}**@{-}="c";"a"**@{-} ,
"a"*+!DR{\st v_1} ,
"b"*+!U{\st v_2} ,
"c"*+!DL{\st v_3} ,
\end{xy}
\,\,\, & \,\,\,
\raisebox{-7.5pt}{$
\begin{xy}
0; <45pt,0pt>: 
(0,1)="a"*=0{\bullet};(0,0)*=0{\bullet}**@{-}="d";
(1,0)*=0{\bullet}**@{-}="b";(1,1)*=0{\bullet}**@{-}="c";
"a"**@{--} ,
"a"*+!DR{\st v_1} ,
"b"*+!UL{\st v_3} ,
"c"*+!DL{\st v_4} ,
"d"*+!UR{\st v_2} ,
\end{xy}
$}
\,\,\, & \,\,\,
\begin{xy}
0; <30pt,0pt>: 
(0,1)="a"*=0{\bullet};(1,0)*=0{\bullet}**@{-}="b";(2,1)*=0{\bullet}**@{-}="c" ,
"b";(1,1.4142)*=0{\bullet}**@{-}="d" ,
"a"*+!DR{\st v_1} ,
"b"*+!U{\st v_2} ,
"c"*+!DL{\st v_4} ,
"d"*+!D{\st v_3} ,
"b"*\cir<7pt>{ul^dl}*\cir<9pt>{ul^dl} ,
\end{xy}
\,\,\,
\\[2mm]
R_3\text{-point} & E_3\text{-point}
& R_4\text{-point} & S_4\text{-point}
\end{array}
\]
\caption{Associated graphs of $R_3$-, $E_3$-, $R_4$- and $S_4$-points (cf.\
Figure \ref{fig:R3E3R4}).}
\label{fig:graphR3E3R4}
\end{figure}
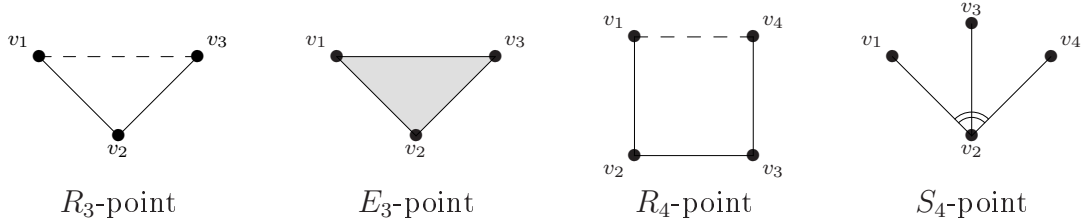

\begin{definition}[The associated graph to $X$] \label{def:dualgraph}
Let $X$ be a good Zappatic surface with Notation \ref{not:Cij}.
The graph $G_X$ associated to $X$ is defined as follows (cf.\ Figure
\ref{fig:graphR3E3R4}):

\begin{itemize}

\item each surface $X_i$ corresponds to a vertex $v_i$;

\item each irreducible component of the
double curve $C_{ij} = C_{ij}^1 \cup \cdots \cup C_{ij}^{h_{ij}}$
corresponds to an edge $e^t_{ij}$, $1 \leq t \leq h_{ij}$, joining $v_i$ and $v_j$.
The edge $e^t_{ij}, i<j$, is oriented from the vertex $v_i$ to the one $v_j$.
The union of all the edges $e_{ij}^t$ joining $v_i$ and $v_j$ is denoted by
$\tilde{e}_{ij}$, which corresponds to the (possibly reducible) double curve $C_{ij}$;

\item each $E_n$-point $P$ of $X$ is a face of the graph whose $n$
edges correspond to the double curves concurring at $P$. This is
called a {\em $n$-face} of the graph;

\item for each $R_n$-point $P$, with $n\ge 3$, if $P \in
X_{i_1}\cap X_{i_2}\cap\cdots\cap X_{i_n}$, where $X_{i_j}$ meets
$X_{i_k}$ along a curve $C_{i_ji_k}$ only if $1=|j-k|$, we add in
the  graph a \emph{dashed edge} joining the vertices corresponding
to $X_{i_1}$ and $X_{i_n}$. The dashed edge $e_{i_1,i_n}$,
together with the other $n-1$ edges $e_{i_j,i_{j+1}}$,
$j=1,\ldots,n-1$, bound an \emph{open $n$-face} of the graph;

\item for each $S_n$-point $P$, with $n\ge 4$, if $P \in
X_{i_1}\cap X_{i_2}\cap \cdots\cap X_{i_n}$, where $X_{i_1},
\ldots, X_{i_{n-1}}$ all meet $X_{i_n}$ along curves $C_{i_ji_n}$,
$j=1,\ldots,n-1$, concurring at $P$, we mark this in the graph by
an \emph{$n$-angle} spanned by the edges corresponding to the curves
$C_{i_ji_n}$, $j=1,\ldots,n-1$.
\end{itemize}
In the sequel, when we speak of {\it faces} of $G_X$ we always mean closed faces.
Of course each vertex $v_i$ is weighted with the relevant
invariants of the corresponding surface $X_i$.
We will usually omit these weights if $X$ is planar,
i.e.\ if all the $X_i$'s are planes.

By abusing notation, we will sometimes denote by $G_X$ also
the natural CW-complex associated to the graph $G_X$.
\end{definition}

Since each $R_n$-, $S_n$-, $E_n$-point is an element of some set of points $\Sigma_{ijk}$
(cf.\ Notation \ref{not:Cij}),
we remark that there can be different
faces (as well as open faces and angles) of $G_X$ which
are incident on the same set of vertices and edges.
However this cannot occur if $X$ is planar.


\begin{figure}[ht]
\[
\begin{xy}
0; <30pt,0pt>: 
(0,1)="a"*=0{\bullet};(1,0)*=0{\bullet}**@{-}="b";(2,1)*=0{\bullet}**@{-}="c",
"a"*+!DR{\st v_1} ,
"b"*+!U{\st v_2} ,
"c"*+!DL{\st v_3} ,
\end{xy}
\]
\caption{Associated graph to a $R_3$-point in a good, planar Zappatic surface.}
\label{fig:graphR3plan}
\end{figure}
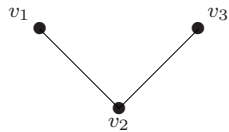

Consider three vertices $v_i, v_j, v_k$ of $G_X$ in such a way
that $v_i$ is joined with $v_j$ and $v_k$. Assume for simplicity that the double
curves $C_{ij}$, $1 \leq i < j \leq v$,
are irreducible. Then, any point in
$C_{ij}\cap C_{ik}$ is either a $R_n$-, or an $S_n$-, or an
$E_n$-point, and the curves $C_{ij}$ and $C_{ik}$ intersect
transversally, by definition of Zappatic singularities.
Hence we can compute the intersection number $C_{ij}\cdot C_{ik}$
by adding the number of closed and open faces and of angles
involving the edges $e_{ij}, e_{ik}$. In particular, if $X$ is
planar, for every pair of adjacent edges only one of
the following possibilities occur: either they belong to an open
face, or to a closed one, or to an angle.
Therefore for good, planar Zappatic surfaces
we can avoid marking open $3$-faces without losing any information
(see Figures \ref{fig:graphR3E3R4} and \ref{fig:graphR3plan} ).

\vspace{-2mm}

As for stick curves, if $G$ is a given graph as above, there
does not necessarily exist a good planar
Zappatic surface $X$ such that its associated graph is $G=G_X$.

\begin{figure}[ht]
\[
\begin{xy}
0; <45pt,0pt>: 
(0,1)="a", (1,0)="b" ,
"a"."b"*[!\colorxy(0.875 0.875 0.875)]\frm{**} ,
"a"*=0{\bullet};(0,0)*=0{\bullet}**@{-}="d";
"b"*=0{\bullet}**@{-};(1,1)*=0{\bullet}**@{-}="c";
"a"**@{-};"b"**@{-} ,
"a"*+!DR{\st v_1} ,
"b"*+!UL{\st v_3} ,
"c"*+!DL{\st v_4} ,
"d"*+!UR{\st v_2} ,
\end{xy}
\]
\caption{Graph associated to an impossible planar Zappatic surface.}
\label{fig:impo}
\end{figure}
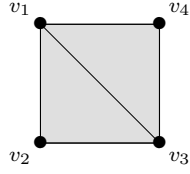

\begin{example}\label{ex:impo}
Consider the graph $G$ of Figure \ref{fig:impo}.
If $G$ were the associated graph to a good planar Zappatic surface
$X$, then $X$ should be a global normal crossing union
of $4$ planes with $5$ double lines and two $E_3$ points,
$P_{123}$ and $P_{134}$, both lying on the double line $C_{13}$.
Since the lines $C_{23}$ and $C_{34}$ (resp.\ $C_{14}$ and
$C_{12}$) both lie on the plane $X_3$ (resp.\ $X_1$), they should
intersect. This means that the planes $X_2,X_4$ also should
intersect along a line, therefore the edge $e_{24}$ should appear in the graph.
\end{example}


Analogously to Example \ref{ex:impo}, one can easily see that,
if the $1$-skeleton of $G$ is $E_3$ or $E_4$, then in order to have a planar
Zappatic surface $X$ such that $G_X=G$, the $2$-skeleton of $G$ has to consist
of the face bounded by the $1$-skeleton.

Let us see two more examples of planar Zappatic surfaces.

\begin{example}\label{ex:zapEn}
In $\Pp^4$, with homogeneous coordinates $x_0,\ldots,x_4$, consider
the good planar Zappatic surface $X$ which is union of the five planes
$$
X_0= \{x_4=x_0=0\}, \; X_i= \{ x_i=x_{i-1}=0\}, \; i=1,2,3,4.
$$
The associated graph is a cycle $E_5$ with no closed faces and the
Zappatic singularities are five $R_3$-points, which, according to the previous remark,
we do not mark with open $3$-faces.
\end{example}

\begin{example}\label{ex:zapEn2}
In $\Pp^5$, with affine coordinates $x_1,\ldots,x_5$,
the planar Zappatic surface $X$, which is the
union of the five planes
\begin{align*}
X_1&=\{ x_4=x_5=x_1=0 \}, \qquad\qquad X_2 =\{ x_5=x_1=x_2=0\}, \\
X_i&=\{ x_{i-2}=x_{i-1}=x_i=0 \},  \quad i  = 3,4,5,
\end{align*}
has an $E_5$-point at $(0,0,0,0,0)$.  The associated graph is a again a cycle
$E_5$ but with a closed $5$-face.
\end{example}

It would be interesting to characterize all the graphs which can be associated
to a good Zappatic surface.

Let us see some examples of a good, non-planar, Zappatic surface.

\begin{example}
Consider $X \subset \Pt$ the union of two general quadrics $X_1$
and $X_2$ and a general plane $X_3$. Then, $C_{12} = C_{21}$ is a smooth
elliptic
quartic in $\Pt$ whereas $C_{13} = C_{31}$ and $C_{23} = C_{32}$ are smooth
conics;
moreover,
$$
X_{1} \cap X_{2} \cap X_{3} =
\Sigma_{123} = \Sigma_{213} = \cdots = \Sigma_{321}
$$
consists of four distinct points. Hence, $G_X$ has three vertices,
three edges (in a cycle) and four triangles (i.e.\ $3$-faces)
which are incident on the same set of vertices (equiv.\ edges).
\end{example}

We can also consider an example of a good Zappatic surface
with reducible double curves.

\begin{example}\label{ex:controes}
Consider $D_1$ and $D_2$ two general
plane curves of degree $m$ and $n$, respectively. Therefore, they are smooth, irreducible and
they transversally intersect each other in $mn$ points. Consider the surfaces:
$$X_1 = D_1 \times \Pp^1 \; \; {\rm and} \;\; X_2 = D_2 \times \Pp^1.$$The union of these two
surfaces, together with the plane $\Pp^2 = X_3$ containing the two curves, determines a good
Zappatic surface $X$ with only $E_3$-points as Zappatic singularities.

More precisely, by using Notation \ref{not:Cij}, we have:
\begin{itemize}
\item $C_{13} = X_1 \cap X_3 = D_1$, $C_{23} = X_2 \cap X_3 = D_2$,
$C_{12} = X_1 \cap X_2 = \sum_{k=1}^{mn}F_k$, where each $F_k$ is a fibre isomorphic to
$\Pp^1$;
\item $\Sigma_{123} = X_1 \cap X_2 \cap X_3$ consists of the
$mn$ points of the intersection of $D_1$ and $D_2$ in $X_3$.
\end{itemize}Observe that $C_{12}$ is smooth but not irreducible.
Therefore, the associated graph to $X$, i.e.\ $G_X$, consists of
3 vertices, $mn+2$ edges and $mn$ triangles incident on them.
\end{example}

In order to combinatorially compute some of the invariants of a good Zappatic
surface,
we need some notation.

\begin{notation}\label{def:vefg}
Let $X$ be a good Zappatic surface (with invariants as in
Notation \ref{not:Cij}) and let $G= G_X$ be its associated graph. We denote by
\begin{itemize}

\item $V :$ the (indexed) set of vertices of $G$;

\item $v :$ the cardinality of $V$, i.e.\ the number
of irreducible components of $X$;

\item $E :$ the set of edges of $G$; this is indexed by the ordered triples
$(i,j,t) \in V \times V \times \NN$, where $i <j$ and $ 1 \leq t \leq h_{ij}$, such that
the corresponding surfaces $X_i$, $X_j$ meet along the curve
$C_{ij} = C_{ji} = C_{ij}^1 \cup \cdots \cup C_{ij}^{h_{ij}}$;

\item $e :$ the cardinality of $E$, i.e.\ the number of irreducible components of double curves in $X$;

\item $\tilde{E}:$ the set of double curves $C_{ij}$ of $X$;
this is indexed by the ordered pairs
$(i,j) \in V \times V $, where $i <j$, such that
the corresponding surfaces $X_i$, $X_j$ meet along the curve
$C_{ij} = C_{ji}$;

\item $\tilde{e} :$ the cardinality of $\tilde{E}$, i.e.\ the pairs of vertices
of $G_X$ which are joined by at least one edge;

\item $f_n :$ the number of $n$-faces of
$G$, i.e.\ the number of $E_n$-points of $X$, for $n \geq 3$;

\item $f := \sum_{n \geq 3} f_n$, the number of faces of $G$,
i.e.\ the total number of $E_n$-points of $X$, for all $n \geq3$;

\item $r_n :$ the number of open $n$-faces of $G$, i.e.\ the number of
$R_n$-points of $X$, for $n\ge 3$;

\item $r$:$= \sum_{n \geq 3} r_n$,  the total number of $R_n$-points of $X$, for
all $n \geq3$;

\item $s_n :$ the number of $n$-angles of $G$, i.e.\ the number of
$S_n$-points of $X$, for $n\ge 4$;

\item $s$: $= \sum_{n \geq 4} s_n$: the total number of $S_n$-points of $X$, for
all $n \geq4$;

\item $\rho_n$: $= s_n + r_n$, for $n \geq 4$, and $\rho_3 = r_3$;

\item $\rho$: $= s + r = \sum_{n \geq 3} \rho_n$;

\item $\tau$: $=\rho +f$, the total number of good Zappatic singularities;

\item $w_i$: the valence of the $i^{\rm th}$ vertex $v_i$ of $G$, i.e.\ the number
of irreducible double curves lying on $X_i$;

\item $\chi(G) :=v-e+f$, i.e.\ the Euler-Poincar\'e characteristic of $G$;

\item $G^{(1)} :$ the \emph{$1$-skeleton} of $G$, i.e.\ the graph
obtained from $G$ by forgetting all the faces, dashed edges and angles;

\item $\chi(G^{(1)}) =v-e$, i.e.\ the Euler-Poincar\'e characteristic of
$G^{(1)}$.
\end{itemize}
\end{notation}

\begin{remark}\label{rem:G^(1)}
Observe that, when $X$ is a good, planar Zappatic surface, $E = \tilde{E}$ and
the $1$-skeleton $G^{(1)}_X$ of $G_X$
coincides with the dual graph $G_D$ of the general hyperplane section $D$ of
$X$.
\end{remark}

Now we can compute some of the invariants of good Zappatic surfaces.

\begin{proposition}\label{prop:ghypsect} (cf. \cite[Proposition 3.12]{CCFMto})  Let $X = \bigcup_{i=1}^v X_i \subset \PR$ be a good Zappatic
surface and let $G = G_X$ be its associated graph. Let $C$ be the
double locus of $X$, i.e.\ the union of the double curves of $X$,
$C_{ij} = C_{ji} = X_i \cap X_j$ and let $c_{ij} = \deg(C_{ij})$.
Let $D_i$ be a general hyperplane section of $X_i$, and denote by $g_i$ its
genus. Then
the arithmetic genus of a general hyperplane section $D$ of $X$ is:
\begin{equation}
g = \sum_{i=1}^v g_i + \sum_{1\le i <j \le v} \; c_{ij} -v +1.
\label{eq:g}
\end{equation}
In particular, when $X$ is a good, planar Zappatic surface, then
\begin{equation}\label{eq:gplanar}
g = e- v + 1 = 1 - \chi(G^{(1)}).
\end{equation}
\end{proposition}
\begin{proof} Denote by $d_i$ the degree of $X_i$, $1 \leq i \leq v$. Then,
$D$ is the union of the $v$ irreducible components $D_i$, $1 \leq
i \leq v$, such that $\deg(D_i) = d_i$ and $d := \deg(D) =
\sum_{i=1}^v d_i$. Consider its associated graph $G_D$, defined as
in \S \ \ref{S:2}.

Take $G$, whose indexed set of edges is denoted by $E$, and
consider an edge $e^t_{ij} \in E$ joining its vertices $v_i$ and $v_j$, $i
<j$, which correspond to the irreducible components $X_i$ and
$X_j$, respectively. The edge $e^t_{ij}$ in $G$ corresponds to an irreducible
component $C_{ij}^t$ of the
double curve $C_{ij}$, $1 \leq t \leq h_{ij}$; its degree is denoted by
$c_{ij}^t$, so that $c_{ij} = \sum_{t=1}^{h_{ij}} c^t_{ij}$.

Thus, we have exactly $c_{ij}$ oriented edges in
the graph $G_D$ joining its vertices $v_i$ and $v_j$, which now
correspond to the irreducible components of $D$, $D_i$ and $D_j$,
respectively. These $c_{ij}$ oriented edges correspond to the
$c_{ij}$ nodes of the reducible curve $D_i \cup D_j$, which is
part of the hyperplane section $D$.

Now, recall that the Hilbert polynomial of $D$ is, with our
notation, $P_{D}(t) = d t + 1 - g$. On the other hand, $P_{D}(t)$
equals the number of independent conditions imposed on
hypersurfaces $\Ha$ of degree $t \gg 0$ to contain $D$.

From what observed above on $G_D$, it follows that the number of singular
points of $D$ is $\sum_{\tilde{e}_{ij}\in \tilde{E}} c_{ij}$. These points impose independent
conditions on hypersurfaces $\Ha$ of degree $t \gg 0$.

Since $t\gg0$ by assumption, we get that the map
$$
H^0(\Oc_{\PR} (t)) \to H^0(\Oc_{D_i}(t))
$$
is surjective and that the line bundle $\Oc_{D_i}(t)$ is
non-special on $D_i$, for each $1 \leq i \leq v$. Thus, in order
for $\Ha$ to contain $D_i$ we have to impose $d_i t - g_i + 1 -
\sum_{j \; {\rm s.t.} \;\tilde{e}_{ij} \in \tilde{E} } c_{ij}$ conditions. Therefore the total
number of conditions for $\Ha$ to contain $D$ is:
\begin{align*}
\sum_{\tilde{e}_{ij}\in \tilde{E}} c_{ij} + \sum_{i=1}^v \biggl( d_i t -g_i +1 -
\sum_{j, \tilde{e}_{ij}\in \tilde{E}} c_{ij} \biggr) & = \sum_{\tilde{e}_{ij}\in \tilde{E}} c_{ij}
+ dt - \sum_{i=1}^v g_i + v
 - \sum_{i=1}^v \sum_{j, \tilde{e}_{ij}\in \tilde{E}} c_{ij}= \\
&= d t + v - \sum_{i=1}^v g_i - \sum_{\tilde{e}_{ij} \in \tilde{E}} c_{ij},
\end{align*}
since $\sum_{i=1}^v \sum_{j , \tilde{e}_{ij}\in \tilde{E}} c_{ij} = 2
\sum_{\tilde{e}_{ij}\in \tilde{E}} c_{ij}$.
This proves \eqref{eq:g} (cf.\ Formula \eqref{eq:genuscurves2}).

The second part of the statement directly follows from the above
computations and from the fact that, in the good planar Zappatic
case $g_i=0$ and $c_{ij} = 1$, for each $i < j$, i.e.\ $G_D$ coincides with
$G^{(1)}$ (cf.\ Remark \ref{rem:G^(1)}).
\end{proof}

By recalling Notation \ref{def:vefg}, one also has:

\begin{proposition}\label{prop:4.chi} (cf. \cite[Proposition 3.15]{CCFMto}) Let $X = \bigcup_{i=1}^v X_i $ be a good Zappatic surface and
$G_X$ be its associated graph, whose number of faces is $f$.
Let $C$ be the double locus of $X$, which is the
union of the curves $C_{ij} = X_i \cap X_j$. Then:
\begin{equation} \label{eq:chi}
\chi(\Oc_X) = \sum_{i=1}^v \chi(\Oc_{X_i}) - \sum_{1 \leq i < j \leq v}
\chi(\Oc_{C_{ij}}) + f .
\end{equation}

In particular, when $X$ is a good, planar Zappatic surface, then
\begin{equation} \label{eq:chiplan}
\chi(\Oc_X) = \chi(G_X) = v - e + f .
\end{equation}
\end{proposition}
\begin{proof}

We can consider the sheaf morphism:
\begin{equation}\label{eq:aiuto2}
\bigoplus_{i=1}^v \Oc_{X_i}
\xrightarrow{\;\lambda\;} \bigoplus_{1 \le i<j \le v} \Oc_{C_{ij}} ,
\end{equation}
defined in the following way: if
$$
\pi_{ij}: \bigoplus_{1 \le i < j \le v} \Oc_{C_{ij}} \to
\Oc_{C_{ij}}
$$
denotes the projection on the $(ij)^{\rm th}$-summand, then
$$
(\pi_{ij} \circ \lambda) (h_1, \ldots , h_v) := h_i - h_j.
$$
Notice that the definition of $\lambda$ is consistent with the
lexicographic order of the indices and with the lexicographic
orientation of the edges of the graph $G_X$.

Observe that, if $\tilde{X}$ denotes the minimal desingularization of $X$,
then $\tilde{X}$ is isomorphic to the disjoint union of the
smooth, irreducible components $X_i$, $1 \leq i \leq v$, of $X$.
Therefore, by the very definition of $\Oc_X$, we see that
$$
\ker(\lambda) \cong \Oc_X.
$$

We claim that the morphism $\lambda$ is not surjective and that
its cokernel is a sky-scraper sheaf supported at the $E_n$-points
of $X$, for $n \geq 3$. To show this, we focus on any irreducible component of
$C
= \bigcup_{1 \leq i < j \leq v} C_{ij}$, the double locus of $X$. For simplicity,
we shall assume that each curve $C_{ij}$ is irreducible; one can easily extend the same
computations to the general case.

Fix any index pair $(i,j)$, with $i <j$, and consider the
generator
\begin{equation}\label{eq:gener}
(0, \ldots, 0 , 1 , 0 , \ldots, 0) \in \bigoplus_{1 \leq l < m
\leq v} \Oc_{C_{lm}},
\end{equation}
where $1 \in \Oc_{C_{ij}}$, the $(ij)^{\rm th}$-summand. The
obstructions to lift up this element to an element of
$\bigoplus_{1 \leq t \leq v} \Oc_{X_t}$ are given by the presence
of good Zappatic singularities of $X$ along $C_{ij}$.

For what concerns the irreducible components
of $X$ which are not involved in the intersection determining a
good Zappatic singularity on $C_{ij}$, the element in
\eqref{eq:gener} trivially lifts-up to $0$ on each of them.
Thus, in the sequel, we shall focus only on the irreducible components
involved in the Zappatic singularity, which will be denoted by $X_i, \; X_j, \;
X_{l_t}$, for $ 1 \leq t \leq n-2$.

We have to consider different cases, according to the good
Zappatic singularity type lying on the curve $C_{ij} = X_i \cap
X_j$.

\begin{itemize}
\item Suppose that $C_{ij}$ passes through a $R_n$-point $P$ of
$X$, for some $n$; we have two different possibilities. Indeed:

\noindent (a) let $X_i$ be an ``external'' surface for $P$ ---
i.e.\ $X_i$ corresponds to a vertex of the associated graph to $P$
which has valence $1$. Therefore, we have:
\[
\begin{picture}(150,15)
\thinlines \put(0,10){\line(1,0){75}}
\put(150,10){\line(-1,0){45}} \put(0,10){\circle*{3}}
\put(30,10){\circle*{3}} \put(60,10){\circle*{3}}
\put(120,10){\circle*{3}} \put(150,10){\circle*{3}}
\put(84,7.7){$\cdots$} \put(-5,2){$\scriptstyle X_i$}
\put(25,2){$\scriptstyle X_j$} \put(55,2){$\scriptstyle X_{l_1}$}
\put(110,2){$\scriptstyle X_{l_{n-3}}$} \put(140,2){$\scriptstyle
X_{l_{n-2}}$}
\end{picture}
\]
%
In this situation, the element in \eqref{eq:gener} lifts up to
$$
(1, 0, \ldots, 0) \in \Oc_{X_i} \oplus \Oc_{X_j} \oplus
\bigoplus_{1\leq t \leq n-2} \Oc_{X_{l_t}}.
$$

\noindent (b) let $X_i$ be an ``internal'' surface for $P$ ---
i.e.\ $X_i$ corresponds to a vertex of the associated graph to $P$
which has valence $2$. Thus, we have a picture like:
\[
\begin{picture}(210,15)
\thinlines \put(0,10){\line(1,0){135}}
\put(210,10){\line(-1,0){45}} \put(0,10){\circle*{3}}
\put(30,10){\circle*{3}} \put(60,10){\circle*{3}}
\put(90,10){\circle*{3}} \put(120,10){\circle*{3}}
\put(180,10){\circle*{3}} \put(210,10){\circle*{3}}
\put(144,7.7){$\cdots$} \put(-5,2){$\scriptstyle {X_{l_1}}$}
\put(25,2){$\scriptstyle {X_{l_2}}$} \put(55,2){$\scriptstyle
{X_{l_3}}$} \put(85,2){$\scriptstyle {X_i}$}
\put(115,2){$\scriptstyle {X_j}$} \put(170,2){$\scriptstyle
{X_{l_{n-3}}}$} \put(200,2){$\scriptstyle {X_{l_{n-2}}}$}
\end{picture}
\]

In this case, the element in \eqref{eq:gener} lifts up to the
$n$-tuple having components:
\begin{itemize}
\item[] $1 \in \Oc_{X_i}$,
\item[] $0 \in \Oc_{X_j}$,
\item[] $1 \in \Oc_{X_{l_t}}$, for those $X_{l_t}$'s corresponding to
vertices in the graph associated to $P$ which are on the left of $X_i$ and,
\item[] $0 \in \Oc_{X_{l_k}}$ for those
$X_{l_k}$'s corresponding to vertices in the graph
associated to $P$ which are on the right of $X_j$.
\end{itemize}
\item Suppose that $C_{ij}$ passes through a $S_n$-point $P$ of
$X$, for any $n$; as before, we have two different possibilities.
Indeed:

\noindent (a) let $X_i$ correspond to the vertex of valence $n-1$
in the associated graph to $P$, i.e.\ \[
\begin{picture}(240,50)
\thinlines \put(120,40){\line(0,-1){30}}
\put(120,40){\line(-1,-1){30}} \put(120,40){\line(-3,-1){90}}
\put(120,40){\line(-4,-1){120}} \put(120,40){\line(1,-1){30}}
\put(120,40){\line(2,-1){60}} \put(120,40){\line(4,-1){120}}
\put(120,40){\circle*{3}} \put(0,10){\circle*{3}}
\put(30,10){\circle*{3}} \put(90,10){\circle*{3}}
\put(120,10){\circle*{3}} \put(150,10){\circle*{3}}
\put(180,10){\circle*{3}} \put(240,10){\circle*{3}}
\put(204,7.7){$\cdots$} \put(54,7.7){$\cdots$}
\put(116,43){$\scriptstyle {X_i}$} \put(-3,2){$\scriptstyle
{X_{l_1}}$} \put(25,2){$\scriptstyle {X_{l_2}}$}
\put(85,2){$\scriptstyle {X_{l_k}}$} \put(115,2){$\scriptstyle
{X_j}$} \put(140,2){$\scriptstyle {X_{l_{k+1}}}$}
\put(170,2){$\scriptstyle {X_{l_{k+2}}}$}
\put(230,2){$\scriptstyle {X_{l_{n-2}}}$}
\end{picture}
\]

In this situation, the element in \eqref{eq:gener} lifts up to the
$n$-tuple having components:
\begin{itemize}
\item[] $1 \in \Oc_{X_i}$, \item[] $0 \in \Oc_{X_j}$, \item[] $1
\in \Oc_{X_{l_t}}$, for all $1 \leq t \leq n-2$.
\end{itemize}

\noindent (b) let $X_i$ correspond to a vertex of valence $1$
in the associated graph to $P$. Since $C_{ij} \neq \emptyset$ by
assumption, then $X_j$ has to be the vertex of valence $n-1$, i.e.\ we have the following picture:
\[
\begin{picture}(240,50)
\thinlines \put(120,40){\line(0,-1){30}}
\put(120,40){\line(-1,-1){30}} \put(120,40){\line(-3,-1){90}}
\put(120,40){\line(-4,-1){120}} \put(120,40){\line(1,-1){30}}
\put(120,40){\line(2,-1){60}} \put(120,40){\line(4,-1){120}}
\put(120,40){\circle*{3}} \put(0,10){\circle*{3}}
\put(30,10){\circle*{3}} \put(90,10){\circle*{3}}
\put(120,10){\circle*{3}} \put(150,10){\circle*{3}}
\put(180,10){\circle*{3}} \put(240,10){\circle*{3}}
\put(204,7.7){$\cdots$} \put(54,7.7){$\cdots$}
\put(116,43){$\scriptstyle {X_j}$} \put(-3,2){$\scriptstyle
{X_{l_1}}$} \put(25,2){$\scriptstyle {X_{l_2}}$}
\put(85,2){$\scriptstyle {X_{l_k}}$} \put(115,2){$\scriptstyle
{X_i}$} \put(140,2){$\scriptstyle {X_{l_{k+1}}}$}
\put(170,2){$\scriptstyle {X_{l_{k+2}}}$}
\put(230,2){$\scriptstyle {X_{l_{n-2}}}$}
\end{picture}
\]
Thus, the element in \eqref{eq:gener} lifts up to the $n$-tuple
having components
\begin{itemize}
\item[] $1 \in \Oc_{X_i}$, \item[] $0 \in \Oc_{X_j}$, \item[] $0
\in \Oc_{X_{l_t}}$, for all $1 \leq t \leq n-2$.
\end{itemize}

\item Suppose that $C_{ij}$ passes through an $E_n$-point $P$ for
$X$. Then, each vertex of the associated graph to $P$ has valence
$2$. Since such a graph is a cycle, it is clear that no lifting of
\eqref{eq:gener} can be done.
\end{itemize}

To sum up, we see that $\coker (\lambda)$ is supported at the
$E_n$-points of $X$. Furthermore, if we consider
\[
\bigoplus_{1 \le i < j \le v} \Oc_{C_{ij}}
\xrightarrow{\;ev_P\;} \Oc_P = \underline{\CC}_P, \qquad \bigoplus
f_{ij} \mapsto \sum f_{ij} (P)
\]
it is clear that, if $P$ is an $E_n$-point then
$$
ev_P \left(\bigoplus_{1 \le i < j \le v} \Oc_{C_{ij}} /
\im(\lambda)\right) \cong \underline{\CC}_P.
$$
This means that
$$\coker (\lambda) \cong \underline{\CC}^f.$$

By the exact sequences
\[
0 \to \Oc_X  \to  \bigoplus_{1 \le i \le v} \Oc_{X_i}  \to
\im(\lambda)  \to 0, \qquad 0 \to \im (\lambda) \to  \bigoplus_{1
\le i < j \le v} \Oc_{C_{ij}}  \to \underline{\CC}^f  \to 0,
\]
we get \eqref{eq:chi}.
\end{proof}

Not all of the invariants of $X$ can be directly computed by the
graph $G_X$. For example, if $\omega_X$ denotes the dualizing
sheaf of $X$, the computation of the $\omega$-genus
$h^0 (X, \omega_X)$, which plays a fundamental role in degeneration theory,
is actually much more involved, even if $X$ has mild Zappatic singularities,
as we shall see in the next section (cf.\ also \cite{CCFMto}).

To conclude this section, we observe that in the particular case of good,
planar Zappatic surface one can determine a simple relation
among the numbers of Zappatic singularities, as the next lemma shows.

\begin{lemma}\label{lem:r3} (cf. \cite[Lemma 3.16]{CCFMk2})
Let $G$ be the associated graph to a
good, planar Zappatic surface
$X=\bigcup_{i=1}^v X_i$. Then, with Notation \eqref{def:vefg}, we have
\begin{equation} \label{eq:r3}
\sum_{i=1}^v \frac{w_i(w_i-1)}{2}=
\sum_{n \ge 3}\left(nf_n+(n-2)r_n \right)
+ \sum_{n \ge 4} \binom{n-1}{2} s_n.
\end{equation}
\end{lemma}

\begin{proof}
The associated graph to three planes which form a $R_3$-point consists of two
adjacent edges (cf.\ Figure \ref{fig:graphR3plan}).
The total number of two adjacent edges in $G$ is the left hand side
member of \eqref{eq:r3} by definition of valence $w_i$.
On the other hand, an $n$-face (resp.\ an open $n$-face, resp.\ an $n$-angle)
clearly contains exactly $n$
(resp.\ $n-2$, resp.\ $\binom{n-1}{2}$) pairs
of adjacent edges.
\end{proof}

\section{The $\omega$-genus of a Zappatic surface}\label{S:4pg}

The aim of this section is to compute the $\omega$-genus of a good Zappatic
surface $X$, as defined in Formula \eqref{eq:pgX}.
What we will actually do will be to compute the cohomology of the
structure sheaf $\cO_X$, which is sufficient, since
$p_\omega(X)=h^ 2(X,\cO_X)$.

We first define the map $\Phi_X$
which appears in the statement of Theorem \ref{thm:intro0}.

\begin{definition}\label{rem:indices}
Let $X=\bigcup_{i=1}^v X_i$ be a good Zappatic surface. Let $r_{ij}: H^1(X_i,
\Oc_{X_i}) \to H^1 (C_{ij},\Oc_{C_{ij}})$ be the restriction map to
$C_{ij}$ as a divisor in $X_i$. We define the natural map:
\begin{equation}\label{eq:4.fi}
\Phi_X: \bigoplus_{i=1}^v H^1(X_i, \Oc_{X_i}) \to
\!\!\!\!\bigoplus_{1\leq i< j\leq v} \!\!\!\! H^1(C_{ij},
\Oc_{C_{ij}}), \;\;\Phi_X(a_i)=-\sum_{j=1}^{i-1}
r_{ij}(a_i)+\sum_{j=i+1}^{v} r_{ij}(a_i)
\end{equation}
if $a_i\in H^1(X_i, \Oc_{X_i})$ and extend $\Phi_X$ linearly. When
$X$ is clear from the context, we will write simply $\Phi$ instead
of $\Phi_X$.
\end{definition}

The main result of this section is the following (cf. \cite[Theorem 3.1]{CCFMpg} together with the beginning of $\S\;2$ and Definition 2.4 in \cite{CCFMpg}):

\begin{theorem}\label{thm:4.pgbis}  Let $X=\bigcup_{i=1}^v X_i$ be a good Zappatic surface. Then:
\begin{equation}\label{eq:Fpgbis}
p_{\omega}(X)=h^ 2(X,\cO_X)= h^2(G_X, \CC)+\sum_{i=1}^v p_g(X_i)+
\dim({\rm coker}(\Phi)),
\end{equation}
and
\begin{equation}\label {eq:Fq}
h^ 1(X,\omega_X)=h^ 1(X,\cO_X)=h^1(G_X,\CC) + \dim(\ker(\Phi))
\end{equation}
where $G_{X}$ is the associated graph to $X$
and $\Phi=\Phi_X$ is the map of Definition \ref{rem:indices}.
\end{theorem}

\begin{proof} Let $p_1,...,p_f$ be the $E_n$--points of $X$, $n\geq 3$.
As in the proof of Proposition \ref{prop:4.chi} (cf.\ also Proposition 3.15 in \cite {CCFMto}),
one has the exact sequence:
\begin{equation}\label {injres}
0 \to \cO_X \to
\bigoplus_{i=1}^v \cO_{X_i}\buildrel d_G^0\over\longrightarrow
\bigoplus_{1\leq i<j\leq v} \cO_{C_{ij}}
\buildrel d_G^1\over\longrightarrow \bigoplus_{h=1}^f \cO_{p_h}
\to 0
\end{equation}
where the maps are as follows:

\begin{itemize}
\item $\cO_X\to \bigoplus_{i=1}^v \cO_{X_i}$
is the direct sum of the natural restriction maps.

\item recall that $d_G^0:\bigoplus_{i=1}^v \cO_{X_i}\to \bigoplus
_{1\leq i<j\leq v} \cO_{C_{ij}}$, can be described by considering
the composition of its restriction to each summand $\cO_{X_i}$
with the projection to any summand $\cO_{C_{hk}}$, with $h<k$.
This map sends $g\in \cO_{X_i}$ to:

\begin{enumerate}
\item $0\in
\cO_{C_{hk}}$, if both $h,k$ are different from $i$;

\item $g_{\vert C_{ik}} \in \cO_{C_{ik}}$ if $k>i$;

\item $-g_{\vert C_{ki}} \in \cO_{C_{ki}}$ if $k<i$;
\end{enumerate}

\item the map
$d_G^1:\bigoplus _{1\leq i<j\leq v} \cO_{C_{ij}} \to \bigoplus_{h=1}^ f\cO _{p_h}$
again can be described by considering the composition
of its restriction to each summand $\cO_{C_{ij}}$, with $i<j$,
with the projection to any summand $\cO_{p_h}$.
Suppose $p_h$ is an $E_n$--point
corresponding to a face $F_h$ of $G_X$
such that $\partial F_h= \sum_{1\leq i<j\leq v} e_{ij}C_{ij}$,
where either $e_{ij}=0$ or $e_{ij}=\pm 1$.
Then this map sends $g\in \cO_{C_{ij}}$ to $e_{ij}g(p_h)$.

\end{itemize}

We note that the induced maps on global sections in each case
are the corresponding cochain map for the graph $G_X$;
this motivates the notation for these maps used in (\ref{injres}).

Let $\Lambda$ be the kernel of the sheaf map
$d_G^1$,
so that we have two short exact sequences
\begin{equation}\label{sesLambda1}
0\to \cO_X \to \bigoplus_{i=1}^v \cO_{X_i}\to \Lambda \to 0
\end{equation}
and
\begin{equation}\label{sesLambda2}
0 \to \Lambda \to \bigoplus_{1\leq i<j\leq v} \cO_{C_{ij}}
\buildrel d_G^1\over\longrightarrow
\bigoplus_{h=1}^f \cO_{p_h}\to 0.
\end{equation}
The latter gives the long exact sequence:
\[
0\to H^ 0(\Lambda)\to
\bigoplus_{1\leq i<j\leq v}H^0(\cO_{C_ij})
\buildrel d_G^1\over\longrightarrow \bigoplus_{h=1}^f H^0(\cO_{p_h})\to
\]
\[
\to H^1(\Lambda)
\buildrel\beta\over\longrightarrow
\bigoplus_{1\leq i<j\leq v}H^1(\cO_{C_{ij}})
\to 0
\]
and since the cokernel of the map $d_G^1$ is $H^2(G_X,\CC)$,
we derive the short exact sequence
\begin{equation}\label{lambda2}
0\to H^ 2(G_X,\CC)\to H^1(\Lambda)
\buildrel\beta\over\longrightarrow
\bigoplus_{1\leq i<j\leq v}H^1(\cO_{C_{ij}})
\to 0.
\end{equation}

 From the short exact sequence (\ref{sesLambda1})
we have the long exact sequence:
\[
0 \to H^0(\cO_X) \to \bigoplus_{i=1}^vH^0(\cO_{X_i}) \to H^0(\Lambda) \to
\]
\[
\to H^ 1(\cO_X)\to \bigoplus_{i=1}^v H^1(\cO_{X_i})
\buildrel\alpha\over\longrightarrow H^1(\Lambda)\to
H^2(\cO_X)\to \bigoplus_{i=1}^v H^2(\cO_{X_i})\to 0.
\]

Now $H^1(G_X,\CC)$ is the kernel of $d^1_G$
(which is $H^0(\Lambda)$)
modulo the image of $d^0_G$,
which is the image of the map $\bigoplus_{i=1}^vH^0(\cO_{X_i}) \to H^0(\Lambda)$
in the first line above.
Hence we recognize $H^1(G_X,\CC)$ as the cokernel of this map,
and therefore the second line of the above sequence becomes
\[
0\to H^1(G_X,\CC)\to H^ 1(\cO_X)\to \bigoplus_{i=1}^v H^1(\cO_{X_i})
\buildrel\alpha\over\longrightarrow H^1(\Lambda)\to
H^2(\cO_X)\to \bigoplus_{i=1}^v H^2(\cO_{X_i}) \to 0.
\]
Now the composition of the map $\beta$ with the map $\alpha$
is exactly the map $\Phi$: $\Phi=\beta\circ\alpha$.
We claim that $\alpha$ and $\Phi$ have the same kernel,
which by (\ref{lambda2}) is equivalent to having
${\rm Im}(\alpha) \cap H^2(G_X,\CC) (=\ker(\beta)) = \{0\}$.

If we are able to show this, then the leftmost part of the
above sequence would split off as
\[
0\to H^1(G_X,\CC)\to H^ 1(\cO_X)\to \ker(\alpha) = \ker(\Phi) \to 0
\]
which would prove the $H^1$ statement of the theorem.
In addition, if this is true, then the natural surjection
from the cokernel of $\alpha$ to the cokernel of $\Phi$
would have $\ker(\beta) = H^2(G_X,\CC)$ as its kernel, and we would have
$\dim(\coker(\alpha)) = \dim H^2(G_X,\CC) + \dim(\coker(\Phi))$.
Since the rightmost part of the long exact sequence above
splits off as
\[
0 \to \coker(\alpha) \to H^2(\cO_X)\to \bigoplus_{i=1}^v H^2(\cO_{X_i}) \to 0
\]
we see that the $H^2$ statement of the theorem follows also.

To prove that
${\rm Im}(\alpha)\cap H^2(G_X,\CC)=\{0\}$,
notice that the sheaf map $d^0_G$
(which has $\Lambda$ as its image)
factors through obvious maps:
\[
\bigoplus_{i=1}^v \cO_{X_i} \to \bigoplus_{i=1}^v \cO_{C_i}
\to \bigoplus_{1\leq i<j\leq v} \cO_{C_{ij}}
\]
and therefore the map $\alpha$ on the $H^1$ level factors as:
\[
\bigoplus _{i=1}^v H^1(\cO_{X_i})
\to \bigoplus_{i=1}^v H^1(\cO_{C_i}) \to H^1(\Lambda).
\]

Moreover one has the short exact sequence:
\[
0\to \cO_C \to \bigoplus_{i=1}^v \cO_{C_i} \to \Lambda\to 0
\]
where $C$ is the singular locus of $X$, and thus we have an exact sequence:
\begin{equation} \label {lambda3}
H^ 1(C,\cO_C)\to \bigoplus_{i=1}^v H^1(\cO_{C_i}) \to H^1(\Lambda)\to 0.
\end{equation}

We remark that $H^1(C_i,\cO_{C_i})$ [resp.\ $H^ 1(C,\cO_C)$]
is the tangent space at the origin to ${\rm Pic}^0(C_i)$ [resp.\
to ${\rm Pic}^0(C)$] which is a $(\CC^*)^{\delta_i}$--extension
[resp.\ a $(\CC^*)^{\delta}$--extension] of $\bigoplus_{j=1}^v{\rm
Pic}^0(C_{ij})$ [resp.\ of $\bigoplus_{1\leq i<j\leq v} {\rm
Pic}^0(C_{ij})$], where $\delta_i$ [resp.\ $\delta$] depends on
the singular points of $C_i$ [resp.\ of $C$].

There are natural restriction maps:
\[
a: \bigoplus_{i=1}^v {\rm Pic}^0(X_i) \to \bigoplus_{i=1}^v {\rm Pic}^0(C_i)
\]
and
\[
b: {\rm Pic}^0(C) \to \bigoplus_{i=1}^v {\rm Pic}^0(C_i)
\]
which are maps of $\CC^*$--extensions of abelian varieties;
their differentials at the origin are
\[
\bigoplus_{i=1}^v H^ 1(\cO_{X_i}) \to \bigoplus_{i=1}^v H^1(\cO_{C_i})
\]
and
\[
H^1(C,\cO_{C})\to \bigoplus_{i=1}^vH^1(\cO_{C_i})
\]
respectively;
the latter is the leftmost map of the sequence (\ref{lambda3}).

The map $b$
appears in the following exact diagram:
\[
\xymatrix{%
0 \ar[r] & {(\CC^*)^\delta} \ar[r] \ar[d] & {\rm Pic}^0(C) \ar[r] \ar[d]^{b}
& \bigoplus_{i<j}{\rm Pic}^0(C_{ij}) \ar[r] \ar[d] & 0 \\
0 \ar[r] & \bigoplus_i{(\CC^*)^{\delta_i}} \ar[r] & \bigoplus_i{\rm Pic}^0(C_i) \ar[r]
& \bigoplus_{i,j}{\rm Pic}^0(C_{ij}) \ar[r] & 0
}
\]
The vertical map on the right is an injection;
indeed, it is the direct sum of diagonal maps
${\rm Pic}^0(C_{ij}) \to {\rm Pic}^0(C_{ij})\bigoplus {\rm Pic}^0(C_{ji})$.
Therefore, if we denote by $V$ the cokernel of the central map $b$,
we have a short exact sequence of cokernels
\[
0  \to  {(\CC^*)^\gamma}  \to  V  \to \bigoplus_{i<j}{\rm
Pic}^0(C_{ij})  \to  0
\]
for some $\gamma$; in particular,
V is again a $\CC^*$-extension of abelian varieties.
We now recognize by (\ref{lambda3}) that
$H^1(X,\Lambda)$ is the tangent space at the origin to $V$;
moreover the sequence (\ref{lambda2}) is the map on tangent spaces
for the above sequence of groups.
In particular the map $\beta$ is the tangent space map for
the projection $V \to \bigoplus_{i<j}{\rm Pic}^0(C_{ij})$.

Composing $a$ with the projection of
$\bigoplus_{i=1}^v {\rm Pic}^0(C_i)$ to $V$ gives a map
\[
c: \bigoplus_{i=1}^v {\rm Pic}^0(X_i)\to V
\]
whose differential at the origin is the previously encountered map
\[
\alpha: \bigoplus_{i=1}^vH^1(\cO_{X_i}) \to H^1(X,\Lambda).
\]

Now $\bigoplus_{i=1}^v {\rm Pic}^0(X_i)$ is compact,
and therefore the image of $c$ in $V$
has finite intersection with the kernel of the projection
$V\to \bigoplus_{1\leq i<j\leq v}{\rm Pic}^0 (C_{ij})$.
At the tangent space level, this means that the
image of $\alpha$ has trivial intersection
with the kernel of the map $\beta$,
which we have identified as $H^2(G_X,\CC)$,
which was to be proved.
\end{proof}

\begin{remark} Note that, in particular, Formulas \eqref{eq:Fpgbis} and
\eqref{eq:Fq} agree with, and imply, Formula \eqref{eq:chi}
that
we proved in Proposition \ref{prop:4.chi}.
\end{remark}

In case $X$ is a planar Zappatic surface, Theorem \ref{thm:4.pgbis}
implies the following:

\begin{corollary}
Let $X$ be a good, planar Zappatic surface. Then,
\begin{align}
 p_\omega(X)  &=  b_2(G_X), \\
 q(X)  &=  b_1(G_X).
\end{align}
\end{corollary}

\section{Degenerations to Zappatic surfaces}\label{S:zapdeg}

In this section we will focus on flat degenerations of smooth
surfaces to Zappatic ones.

\begin{definition}\label{def:degen}
Let $\D$ be the spectrum of a DVR (equiv.\ the complex unit disk).
A {\em degeneration} of relative dimension $n$ is a proper
and flat  morphism
\[
\xymatrix@R=5mm{\X \ar[d]^\pi \\
\Delta}
\]
such that $\X_t = \pi^{-1}(t)$ is a smooth, irreducible,
$n$-dimensional, projective variety, for $t \neq 0$.

If $Y$ is a smooth, projective variety, the degeneration
\[
\xymatrix@C=0mm@R=5mm{\X \ar[d]_\pi & \subseteq & \Delta \times Y \ar[dll]^{\text{pr}_1} \\
\Delta}
\]
is said to be an {\em embedded degeneration} in $Y$ of relative dimension $n$.
When it is  clear from the context, we will omit
the term {\em embedded}.

We will say that $\X\to\D$ is a \emph{normal crossing degeneration}
if the total space $\X$ is smooth and the support $X_\red$ of the central
fibre $X=\X_0$ is a divisor in $\X$ with global normal crossings,
i.e.\ $X_\red$ is a good Zappatic surface
with only $E_3$-points as Zappatic singularities.

A normal crossing degeneration is called \emph{semistable} (see, e.g., \cite{Morr})
if the central fibre is reduced.
\end{definition}

\begin{remark}\label{def:degen2}
Given a degeneration $\pi:\X\to\D$, Hironaka's Theorem on the resolution
of singularities implies that there exists a birational morphism
$\bar\X\to\X$ such that $\bar\X\to\D$ is a normal crossing degeneration,
which we will call a \emph{normal crossing reduction} of $\pi$.

Given a degeneration $\pi:\X\to\D$, the Semistable
Reduction Theorem (see Theorem on p.\ 53--54 in \cite{Kempf}) states that there exists a
base change $\beta : \D \to \D$, defined by $\beta(t) = t^m$, for
some $m$, a semistable degeneration $\tilde\pi: \tilde\X \to \D$
and a diagram
\begin{equation} \label{eq:s}
\xymatrix@C-1mm@R-1mm{\tilde\X \ar[r]^\psi \ar[dr]_{\tilde\pi} & \X_\beta
\ar[d] \ar[r] & \X \ar[d]^\pi \\ & \Delta \ar[r]^\beta & \Delta }
\end{equation}
such that the square is Cartesian and $\psi:\tilde\X\to\X_\beta$
is a birational morphism
obtained by blowing-up a suitable sheaf of ideals on $\X_\beta$. This is called
a \emph{semistable reduction} of $\pi$.
\end{remark}

From now on, we will be concerned with degenerations of relative
dimension two, namely degenerations of smooth, projective
surfaces.

\begin{definition}\label{def:zappdeg} (cf. \cite[Definition 4.2]{CCFMk2}) Let $\X \to \D$ be a degeneration (equiv.\ an embedded degeneration) of surfaces. Denote by $\X_t$ the general fibre,
which is by definition a smooth, irreducible and projective
surface; let $X=\X_0$ denote the central fibre. We will say that
the degeneration is {\em Zappatic} if $X$ is a Zappatic surface,
the total space $\X$ is smooth except for:
\begin{itemize}
\item ordinary double points at points of the double locus of $X$,
which are not the Zappatic singularities of $X$;
\item further singular points at the Zappatic singularities of $X$ of type
$T_n$, for $n \geq 3$, and $Z_n$, for $n \geq 4$,
\end{itemize}
and there exists a birational morphism $\X'\to\X$,
which is the composition of blow-ups at points of the central fibre,
such that $\X'$ is smooth.

A Zappatic degeneration will be called \emph{good} if the central
fibre is moreover a good Zappatic surface. Similarly, an embedded
degeneration will be called a \emph{planar Zappatic degeneration}
if its central fibre is a planar Zappatic surface.

Notice that we require the total space $\X$ of a good Zappatic degeneration to be smooth
at $E_3$-points of $X$.
\end{definition}

On the other hand, if $\pi: \X \to \D$ is an arbitrary degeneration of surfaces
such that $\pi^{-1}(t) = \X_t$, for $t \neq 0$, is by definition a smooth, irreducible
and projective surface and the central fibre $\X_0$ is a good Zappatic surface,
then the total space $\X$ of $\pi$ may have the following singularities:
\begin{itemize}
\item double curves, which are double curves also for $X$;
\item isolated double points along the double curves of $X$;
\item further singular points at the Zappatic
singularities of $X$, which can be isolated or
may occur on double curves of the total space.
\end{itemize}The singularities of the total space $\X$ of either an arbitrary degeneration
with good Zappatic central fibre or of a good Zappatic degeneration
will be described in details in Sections \ref{S:MS} and \ref{S:resolution}.

\begin{notation}\label{not:zappdeg}
Let $\X \to \D$ be a degeneration of surfaces and let $\X_t$ be
the general fibre, which is a smooth, irreducible
and projective surface. Then, we consider the following intrinsic
invariants of $\X_t$:
\begin{itemize}
\item $\chi : = \chi(\Oc_{\X_t})$; \item $K^2 := K_{\X_t}^2$.
\end{itemize}If the degeneration is assumed to be embedded in $\PR$, for some
$r$, then we also have:
\begin{itemize}
\item $d : = \deg(\X_t)$; \item $g := (K+H)H/2 + 1,$ the sectional
genus of $\X_t$.
\end{itemize}
\end{notation}

We will be mainly interested in computing these invariants in
terms of the central fibre $X$. For some of them, this is quite
simple. For instance, when $\X \to \D$ is an embedded degeneration
in $ \PR$, for some $r$,  and if the central fibre $\X_0 = X =
\bigcup_{i=1}^v X_i$, where the $X_i$'s are smooth, irreducible
surfaces of degree $d_i$, $ 1 \leq i \leq v$, then by the flatness
of the family we have
$$
d = \sum_{i=1}^v d_i.
$$

When $X=\X_0$ is a good Zappatic surface (in particular a
good, planar Zappatic surface), we can easily compute some of
the above invariants by using our results of \S~\ref{S:3}.
Indeed, by Propositions \ref{prop:ghypsect} and \ref{prop:4.chi} by the flatness of
the family, we get:

\begin{proposition}\label{prop:gzapdeg}
Let $\X \to \D$ be a degeneration of surfaces and suppose that the central fibre
$\X_0 = X = \bigcup_{i=1}^v X_i$ is a good Zappatic surface.
Let $G = G_X$ be its associated graph (cf.\ Notation \ref{def:vefg}).
Let $C$ be the
double locus of $X$, i.e.\ the union of the double curves of $X$,
$C_{ij} = C_{ji} = X_i \cap X_j$ and let $c_{ij} = \deg(C_{ij})$.

\noindent (i) If $f$ denotes the number of (closed) faces of $G$,
then
\begin{align}
\chi &= \sum_{i=1}^v \chi(\Oc_{X_i}) - \sum_{1 \leq i < j \leq v}
\chi(\Oc_{C_{ij}}) + f. \label{eq:chizapdeg}
\end{align}

Moreover, if $X=\X_0$ is a good, planar Zappatic surface,
then
\begin{align}
\chi &= \chi(G) = v -e + f, \label{eq:chiplanzapdeg}
\end{align}where $e$ denotes the number of edges of $G$.

\noindent (ii) Assume further that $\X \to \D$ is embedded in
$\PR$. Let $D$ be a general hyperplane section of $X$; let $D_i$
be the $i^{\rm th}$-smooth, irreducible component of $D$, which is a
general hyperplane section of $X_i$, and let $g_i$ be its genus.
Then
\begin{align}
g &= \sum_{i=1}^v g_i + \sum_{1 \leq i < j \leq v} \; c_{ij} -v
+1. \label{eq:gzapdeg}
\end{align}
When $X$ is a good, planar Zappatic surface, if $G^{(1)}$ denotes
the $1$-skeleton of $G$, then:
\begin{align}
g &= 1 - \chi(G^{(1)}) = e- v + 1 .\label{eq:gplanzapdeg}
\end{align}
\end{proposition}

In the particular case that $\X \to \D$ is a semistable
degeneration, i.e.\ if $X$ has only $E_3$-points as Zappatic
singularities and the total space $\X$ is smooth, then $\chi$ can
be computed also in a different way by topological methods via the Clemens-Schmid's exact
sequence (cf.\ e.g.\ \cite{Morr}).

Proposition \ref{prop:gzapdeg} is indeed more general: $X$ is
allowed to have any good Zappatic singularity, namely $R_n$-,
$S_n$- and $E_n$-points, for any $n\ge3$, the total space $\X$ is
possibly singular, even in dimension one, and, moreover, our computations do not depend on
the fact that $X$ is smoothable, i.e.\ that $X$ is the central
fibre of a degeneration.

\section{Minimal and quasi-minimal singularities}\label{S:MS}

In this section we shall describe the singularities that the total
space of a degeneration of surfaces has at the Zappatic singularities
of its central fibre.
We need to recall a few general facts about
reduced Cohen-Macaulay singularities and two fundamental concepts
introduced and studied by Koll\'ar in \cite{Kollar} and \cite{KolSB}.

Recall that $V = V_1 \cup \cdots \cup V_r \subseteq \P^n$, a
reduced, equidimensional and non-degenerate scheme is said to be
{\em connected in codimension one} if it is possible to arrange
its irreducible components $V_1$, \ldots, $V_r$ in such a way that
$$
\codim_{V_j} V_j \cap (V_1 \cup \cdots \cup V_{j-1})=1, \; {\rm
for} \; 2 \leq j \leq r.
$$

\begin{remark}\label{rem:pCM}
Let $X$ be a surface in $\PR$ and $C$ be a hyperplane section of
$X$. If $C$ is a projectively Cohen-Macaulay curve, then $X$ is
connected in codimension one. This immediately follows from the
fact that $X$ is projectively Cohen-Macaulay (cf.\ Appendix \ref{S:pCMpG}).
\end{remark}

Given $Y$ an arbitrary algebraic variety, if $y\in Y$ is a
reduced, Cohen-Macaulay singularity then
\begin{equation}\label{eq:CMSing}
\emdim_y(Y) \leq \mult_y(Y)+\dim_y(Y)-1,
\end{equation}
where $\emdim_y(Y) = \dim (\m_{Y,y}/\m^2_{Y,y})$ is the
\emph{embedding dimension} of $Y$ at the point $y$, where
$\m_{Y,y} \subset \Oc_{Y,y}$ denotes the maximal ideal of $y$ in
$Y$ (see, e.g., \cite{Kollar}).

For any singularity $y \in Y$ of an algebraic variety $Y$, let us set
\begin{equation}\label{eq:deltay}
\delta_y(Y)=\mult_y(Y)+\dim_y(Y)-\emdim_y(Y)-1.
\end{equation}
If $y\in Y$ is reduced and Cohen-Macaulay,
then Formula \eqref{eq:CMSing} states that $\delta_y(Y)\geq 0$.

Let $H$ be any effective Cartier
divisor of $Y$ containing $y$. Of course one has
$$\mult_y(H)\geq \mult_y(Y).$$

\begin{lemma}\label{lem:emdimmult} (cf. \cite[Lemma 5.4]{CCFMk2}) In the above setting, if
$\emdim_y(Y)=\emdim_y(H)$, then $\mult_y(H) > \mult_y(Y)$.
\end{lemma}
\begin{proof} Let $f \in \Oc_{Y,y}$ be a local equation defining $H$ around $y$.
If $f \in \m_{Y,y} \setminus \m^2_{Y,y}$ (non-zero), then $f$
determines a non-trivial linear functional on the Zariski tangent
space $T_{y}(Y) \cong (\m_{Y,y} /\m^2_{Y,y})^{\vee}$. By the
definition of $\emdim_y(H)$ and the fact that $f \in \m_{Y,y}
\setminus \m^2_{Y,y}$, it follows that $\emdim_y(H)=\emdim_y(Y)
-1$. Thus, if $\emdim_y(Y)=\emdim_y(H)$, then $f \in \m^h_{Y,y}$,
for some $h \geq 2$. Therefore, $\mult_y(H) \geq h \mult_y(Y) >
\mult_y(Y)$.
\end{proof}

We let
\begin{equation}\label{eq:nu}
\nu:=\nu_y(H) = {\rm min} \{ n \in \NN \; | \; f \in m^n_{Y,y} \}.
\end{equation}
Notice that:
\begin{equation}\label{eq:mult,emdim}
\mult_y(H) \geq \nu \mult_y(Y),\qquad \emdim_y(H)=
\begin{cases}
\emdim_y(Y) & \text{if }\nu>1, \\
\emdim_y(Y)-1 & \text{if }\nu=1.
\end{cases}
\end{equation}

\begin{lemma}\label{lem:delta} (cf. \cite[Lemma 5.7]{CCFMk2}) One has
$$
\delta_y(H)\geq \delta_y(Y).
$$
Furthermore:

\begin{itemize}
\item[(i)] if the equality holds, then either

\begin{enumerate}
\item $\mult_y(H)=\mult_y(Y)$, $\emdim_y(H)=\emdim_y(Y)-1$ and
$\nu_y(H)=1$, or

\item $\mult_y(H)=\mult_y(Y)+1$, $\emdim_y(H)=\emdim_y(Y)$, in
which case $\nu_y(H)=2$ and $\mult_y(Y)=1$;
\end{enumerate}

\item[(ii)] if $\delta_y(H)= \delta_y(Y)+1$, then either
\begin{enumerate}
\item $\mult_y(H)=\mult_y(Y)+1$, $\emdim_y(H)=\emdim_y(Y)-1$, in
which case $\nu_y(H)=1$, or

\item $\mult_y(H)=\mult_y(Y)+2$ and $\emdim_y(H)=\emdim_y(Y)$, in
which case either

\noindent (a) $2 \leq \nu_y(H)\leq 3$ and $\mult_y(Y)=1$, or

\noindent (b) $\nu_y(H)= \mult_y(Y)= 2$.
\end{enumerate}
\end{itemize}
\end{lemma}

\begin{proof} It is a straightforward consequence of \eqref{eq:deltay}, of Lemma
\ref{lem:emdimmult} and of \eqref{eq:mult,emdim}.
\end{proof}

We will say that $H$ has {\it general behaviour} at $y$ if
\begin{equation}\label{eq:genbeh}
\mult_y(H)=\mult_y(Y).
\end{equation}

We will say that $H$ has {\it good behaviour} at $y$ if
\begin{equation}\label{eq:goodbeh}
\delta_y(H)=\delta_y(Y).
\end{equation}

Notice that if $H$ is a general hyperplane section through $y$,
than $H$ has both general and good behaviour.

We want to discuss in more details the relations between the two
notions. We note the following facts:

\begin{lemma}\label{lem:multemb} (cf. \cite[Lemma 5.10]{CCFMk2}) In the above setting:

\begin{itemize}
\item[(i)] if $H$ has general behaviour at $y$, then it has also
good behaviour at $y$;

\item[(ii)] if $H$ has good behaviour at $y$, then either
\begin{enumerate}
\item $H$ has also general behaviour and
$\emdim_y(Y)=\emdim_y(H)+1$,  or \item $\emdim_y(Y)=\emdim_y(H)$,
in which case $\mult_y(Y)=1$ and $\nu_y(H)=\mult_y(H)=2$.
\end{enumerate}
\end{itemize}
\end{lemma}
\begin{proof} The first assertion is a trivial consequence of Lemma
\ref{lem:emdimmult}.

If $H$ has good behaviour and $\mult_y(Y)=\mult_y(H)$, then it is
clear that $\emdim_y(Y)=\emdim_y(H)+1$. Otherwise, if
$\mult_y(Y)\neq \mult_y(H)$, then $\mult_y(H)=\mult_y(Y)+1$ and
$\emdim_y(Y)=\emdim_y(H)$. By Lemma \ref{lem:delta}, $(i)$, we
have the second assertion.
\end{proof}

As mentioned above, we can now give two fundamental definitions
(cf.\ \cite{Kollar} and \cite{KolSB}):

\begin{definition}\label{def:min} (cf. \cite[Definition 5.11]{CCFMk2})
Let $Y$ be an algebraic variety. A reduced,
Cohen-Macaulay singularity $y\in Y$ is called \emph{minimal} if
the tangent cone of $Y$ at $y$ is geometrically reduced and
$\delta_y(Y)=0$.
\end{definition}

\begin{remark} Notice that if $y$ is a smooth point for $Y$, then
$\delta_y(Y)=0$ and we are in the minimal case.\end{remark}

\begin{definition}\label{def:emin} (cf. \cite[Definition 5.13]{CCFMk2}) 
Let $Y$ be an algebraic variety. A reduced,
Cohen-Macaulay singularity $y\in Y$ is called \emph{quasi-minimal}
if the tangent cone of $Y$ at $y$ is geometrically reduced and
$\delta_y(Y)=1$.
\end{definition}

It is important to notice the following fact:

\begin{proposition}\label{prop:min-quasimin} (cf. \cite[Proposition 5.14]{CCFMk2}) 
Let $Y$ be a projective threefold and $y \in Y$ be a point. Let $H$ 
be an effective Cartier divisor of $Y$ passing through $y$.

\begin{itemize}
\item[(i)] If $H$ has a minimal singularity at $y$, then $Y$ has
also a minimal singularity at $y$. Furthermore $H$ has general
behaviour at $y$, unless $Y$ is smooth at $y$ and $\nu_y(H) =
\mult_y(H)=2$.

\item[(ii)] If $H$ has a quasi-minimal, Gorenstein singularity at $y$
then $Y$ has also a quasi-minimal singularity at $y$, unless
either
\begin{enumerate}
\item $\mult_y(H)=3$ and $1\le \mult_y(Y) \le 2$, or
\item $\emdim_y(Y)=4$, $\mult_y(Y)=2$ and $\emdim_y(H)=\mult_y(H)=4$.
\end{enumerate}

\end{itemize}
\end{proposition}

\begin{proof}
Since $y\in H$ is a minimal (resp.\ quasi-minimal) singularity,
hence reduced and Cohen-Macaulay,
the singularity $y\in Y$ is reduced and Cohen-Macaulay too.

Assume that $y\in H$ is a minimal singularity,
i.e.\ $\delta_y(H)=0$.
By Lemma \ref{lem:delta}, (i), and by
the fact that $\delta_y(Y)\geq 0$, one has $\delta_y(Y)=0$.
In particular, $H$ has good behaviour at $y$.
By Lemma \ref{lem:multemb}, (ii), either $Y$ is smooth at $y$
and $\nu_y(H)=2$, or $H$ has general behaviour at $y$.
In the latter case, the tangent cone of $Y$ at $y$
is geometrically reduced, as is the tangent cone of $H$ at $y$.
Therefore, in both cases $Y$ has a minimal singularity at $y$,
which proves (i).

Assume that  $y\in H$ is a quasi-minimal singularity, namely
$\delta_y(H) = 1$. By Lemma \ref{lem:delta}, then
either $\delta_y(Y)=1$ or $\delta_y(Y)=0$.

If $\delta_y(Y)=1$, then the case (i.2) in Lemma \ref{lem:delta}
cannot occur, otherwise we would have $\delta_y(H)=0$, against the
assumption.
Thus $H$ has general behaviour and, as above,
the tangent cone of $Y$ at $y$
is geometrically reduced, as the tangent cone of $H$ at $y$ is.
Therefore $Y$ has a quasi-minimal singularity at $y$.

If $\delta_y(Y)=0$, we have the possibilities listed in Lemma
\ref{lem:delta}, (ii). If (1) holds, we have $\mult_y(H)= 3$,
i.e.\  we are in case (ii.1) of the statement.
Indeed, $Y$ is Gorenstein at $y$ as $H$ is, and therefore
$\delta_y(Y)=0$ implies that $\mult_y(Y) \leq 2$ by Corollary 3.2 in \cite{Sally},
thus $\mult_y(H) \leq 3$, and in fact $\mult_y(H)= 3$ because $\delta_y(H)=1$.
Also the possibilities listed in Lemma \ref{lem:delta}, (ii.2)
lead to cases listed in the statement.
\end{proof}

\begin{remark}\label{rem:minquasimin} From an analytic viewpoint, case $(1)$ in
Proposition \ref{prop:min-quasimin} $(ii)$, when $Y$ is smooth at $y$, can be thought of as
$Y= \P^3$ and $H$ a cubic surface with a triple point at $y$.

On the other hand, case $(2)$ can be thought of as $Y$ being a
quadric cone in $\P^4$ with vertex at $y$ and as $H$ being cut out
by another quadric cone with vertex at $y$. The resulting
singularity is therefore the cone over a quartic curve $\Gamma$ in
$\P^3$ with arithmetic genus $1$, which is the complete
intersection of two quadrics.
\end{remark}

Now we describe the relation between minimal and quasi-minimal
singularities and Zappatic singularities. First we need the
following straightforward remark:

\begin{lemma}\label{prop:mingminzap}
Any $T_n$-point (resp.\ $Z_n$-point) is a minimal
(resp.\ quasi-minimal) surface singularity.
\end{lemma}

The following direct consequence of Proposition
\ref{prop:min-quasimin} will be important for us:

\begin{proposition}\label{prop:mgmz3f} (cf. \cite[Proposition 5.17]{CCFMk2}) 
Let $X$ be a surface with a Zappatic singularity at a point $x\in
X$ and let $\X$ be a threefold containing $X$ as a Cartier divisor.
\begin{itemize}
\item If $x$ is a $T_n$-point for $X$, then $x$ is a minimal
singularity for $\X$ and $X$ has general behaviour at $x$;

\item If $x$ is an $E_n$-point for $X$, then $\X$ has a
quasi-minimal singularity at $x$ and $X$ has general behaviour at
$x$, unless either:
\begin{itemize}

\item[(i)] $\mult_x(X)=3$ and $1\le \mult_x(\X)\le 2$, or

\item[(ii)] $\emdim_x(\X)=4$, $\mult_x(\X)=2$ and $\emdim_x(X)=\mult_x(X)= 4$.
\end{itemize}
\end{itemize}
\end{proposition}

In the sequel, we will need a description of a surface having as a
hyperplane section a stick curve of type $C_{S_n}$, $C_{R_n}$, and
$C_{E_n}$ (cf.\ Examples \ref{ex:tngraphs} and \ref{ex:zngraphs}).

First of all, we recall well-known results about {\em minimal
degree surfaces} (cf.\ \cite{GH}, page 525).

\begin{theorem}[del Pezzo]\label{thm:DelP1}
Let $X$ be an irreducible, non-degenerate surface of minimal
degree in $\PR$, $r\ge3$. Then $X$ has degree $r-1$ and is one of
the following:
\begin{itemize}
\item[(i)] a rational normal scroll; \item[(ii)] the Veronese
surface, if $r=5$.
\end{itemize}
\end{theorem}

Next we recall the result of Xamb\'o concerning reducible minimal
degree surfaces (see \cite{X}).

\begin{theorem}[Xamb\'o]\label{thm:X}
Let $X$ be a non-degenerate surface which is connected in
codimension one and of minimal degree in $\PR$, $r\ge3$. Then, $X$
has degree $r-1$, any irreducible component of $X$ is a minimal
degree surface in a suitable projective space and any two
components intersect along a line.
\end{theorem}

In what follows, we shall frequently
refer to Appendix \ref{S:pCMpG}.
Let $X \subset \Pp^r$ be an irreducible, non-degenerate,
projectively Cohen-Macaulay surface with canonical singularities,
i.e.\ with Du Val singularities. We recall that $X$ is called a
\emph{del Pezzo surface} if $\Oc_X(-1)\simeq\omega_X$. We note
that a del Pezzo surface is projectively Gorenstein
(cf.\ Definition \ref{def:6bis} in Appendix \ref{S:pCMpG}).

\begin{theorem}\label{thm:DelP2} (del Pezzo, \cite{DelP}. Compare also with \cite[Theorem 5.20]{CCFMk2}) 
Let $X$ be an irreducible, non-degenerate, linearly normal surface
of degree $r$ in $\Pp^r$. Then one of the following occurs:
\begin{itemize}
\item[(i)] one has $3\le r \le 9$ and $X$ is either

\begin{itemize}
\item[a.] the image of the blow-up of $\Pp^2$ at $9-r$ suitable
points, mapped to $\Pp^r$ via the linear system of cubics through
the $9-r$ points, or

\item[b.] the 2-Veronese image in $\Pp^8$ of a quadric in $\Pp^3$.
\end{itemize}In each case, $X$ is a del Pezzo surface.

\item[(ii)] $X$ is a cone over a smooth elliptic normal curve of
degree $r$ in $\Pp^{r- 1}$.
\end{itemize}
\end{theorem}
\begin{proof}
Let $f: Y \to X$ be a minimal desingularization of $X$. Let $H$ be a general
hyperplane section of $X$ and let $C := f^*(H)$. One has
$0 \leq p_a(H) \leq 1$. On the other hand, $C$ is smooth and irreducible (by Bertini's
theorem) of genus $g \leq p_a(H)$. By the
linear normality, one has $g= 1$ and, therefore,
also $p_a(H) = g= 1$. So $H$ is smooth and irreducible, which means that $X$ has
isolated singularities.

Assume that $X$ is not a scroll. If $r \geq 5$, Reider's Theorem states that
$K_Y + C$ is b.p.f. on $Y$. Thus, $(K_Y + C)^2 \geq 0$. On the other hand,
$(K_Y + C) C = 0$. Then, the Index Theorem implies that $K_Y$ is
numerically equivalent to $-C$. Therefore,
$H^1(Y, \Oc_Y (K_Y)) = (0)$ and so $Y$ is rational and $K_Y$ is linearly
equivalent to $-C$.

By the Adjunction Formula, the above relation is trivially true also
if $r = 3$ and $r =4$.

Now, $r = C^2 = K_Y^2 \leq 9$.
If $ r=9$, then $Y = \Pp^2$ and $ C \in |\Oc_{\Pp^2} (3)|$. If $Y$ has no $(-1)$-curve,
the only other possibility is $ r = K^2_Y = 8$ which leads right away to case
$(i)-b$.

Suppose now that $E$ is a $(-1)$-curve in $Y$, then $CE = 1$. We claim that
$|C + E|$ is b.p.f., it contracts $E$ and maps $Y$ to a surface of degree
$r+1$ is $\Pp^{r+1}$. Then, the assertion immediately
follows by the description of the cases $r = 8$ and $9$. To prove the claim,
consider the exact sequence
$$0 \to \Oc_Y(C) \to  \Oc_Y(C+E) \to \Oc_E \to 0$$and remark that
$h^1(\Oc_Y(C)) = h^1 (\Oc_Y(-2C)) = 0$.

Now suppose that $X$ is a scroll which is not a cone.
Note that $Y$ is an elliptic ruled surface. Let $R$ be the pull-back via $f$ of a line $L$.
We claim that all curves in the algebraic system $\{ R \}$ are irreducible. Otherwise,
we would have some $(-1)$-curve on $Y$ contracted by $f$, against the minimality
assumption. Therefore $Y$ is a minimal, elliptic ruled surface. Moreover,
$f : Y \to X$ is finite since $f$ cannot contract any curve
transversal to $R$, otherwise $X$ would be a cone, and
cannot contract any curve $R$.

Now, $Y = \Pp_E (\mathcal{E})$, where $E$ is an elliptic curve and
$\mathcal{E}$ is a rank-two vector bundle on $E$. If $\mathcal{E}$ is
indecomposable, then the $e$-invariant of the scroll is either
$e = 0$ or $e = -1$ (cf.\ Theorem 2.15, page 377 in \cite{H}).
Let $C_0$ be a section of the ruling with $C_0^2 = e$. Then,
$ C \equiv C_0 + \alpha R$, with $\alpha \geq 2 $.

More precisely, we have $\alpha \geq 3$. Otherwise, we would have
$r = C_0^2 + 4 \leq 4$, then the hyperplane section $H$ of $X$ would
be a complete intersection so $X$ would be a complete intersection
hence a cone. Furthermore, when $\alpha =3$
then $e=0$. Indeed, assume
$e=-1$ and $\alpha=3$, so $C_0C =2$. This would imply that
$X$ is a surface of degree $5$ in $\Pp^5$ with a double line (which is
the image of $C_0$). If we project $X$ from this line, we have a curve
of degree $3$ in $\Pp^3$ which contradicts that $Y$ is an elliptic scroll.

Notice that $K_Y \equiv - 2 C_0 - e R$ and therefore
$K_Y - C \equiv - 3 C_0 - (e+ \alpha) R$.
Since $(C - K_Y) C_0 = \alpha + 4 e$, from what observed above, in any case
$(C-K_Y) C_0 >0$. Since
$(C - K_Y)^2 = 6 \alpha + 15 e > 0$, we have that $C-K_Y$ is big and
nef. Therefore, $h^1 (Y, \Oc_Y(C)) = h^1 (Y, \Oc_Y(K_Y - C) ) =0$. Look now at
the sequence
$$0 \to \Oc_Y \to  \Oc_Y(C) \to \Oc_C(C) \to 0.$$Since
$h^1(Y, \Oc_Y) = 1$, then the restriction map
$$H^0(\Oc_Y(C)) \to H^0(\Oc_C(C))$$is not
surjective, against the hypotesis that $X$ is a surface of degree
$r$ in $\Pp^r$.

Finally, assume that $\mathcal{E}$ is decomposable. Then, $\mathcal{E} = {\mathcal L}_1
\oplus {\mathcal L}_2$, where ${\mathcal L}_i$ is a line bundle of degree $d_i$ on
$E$, $ 1 \leq i \leq 2$. Observe that $d_1 + d_2 = r$; furthermore, since $X$ is
not a cone, then $h^0 (E, {\mathcal L}_i) \geq 2$, for each $1 \leq i \leq 2$, hence
$d_i \geq 2$, for $ 1 \leq i \leq 2$. Thus, $h^0 (E, \mathcal{E}) = h^0 (E, {\mathcal L}_1 )
+ h^0 (E, {\mathcal L}_2) = d_1 + d_2 = r$, a contradiction.
\end{proof}

Since cones as in (ii) above are projectively Gorenstein surfaces (see Appendix \ref{S:pCMpG}),
the surfaces listed in Theorem \ref{thm:DelP2} will be called
\emph{minimal Gorenstein surfaces}.

We shall make use of the following easy consequence of the
Riemann-Roch theorem.

\begin{lemma}\label{lem:pa=d} (Compare also with \cite[Lemma 5.21]{CCFMk2}) 
Let $D \subset \PR$ be a reduced (possibly reducible),
non-degenerate and linearly normal curve of degree $r+d$ in $\PR$,
with $0\le d<r$. Then $p_a(D)=d$.
\end{lemma}
\begin{proof}
Let $\Oc_D(H)$ be the hyperplane line bundle on $D$.
By assumption $h^0(D,\Oc_D(H))=r+1$ and $\deg(\Oc_D(H))=r+d$.
Riemann-Roch Theorem then gives:
\begin{equation}\label{eq:pa=d}
p_a(D)=h^1(D, \Oc_D(H))+d=h^0(D,\omega_D\otimes \Oc_D(-H))+d,
\end{equation}where $\omega_D$ denotes the dualizing sheaf.
Suppose that $h^0(\omega_D\otimes \Oc_D(-H))>0$.
Thus, the effectiveness of $\Oc_D(H)$ and $\omega_D\otimes \Oc_D(-H)$ would
imply that:
\[
p_a(D)=h^0(D,\omega_D)\ge h^0(\Oc_D(H))+h^0(\omega_D \otimes \Oc_D(-H))-1
= r + h^0(\omega_D\otimes \Oc_D(-H)),
\]
which contradicts \eqref{eq:pa=d}, since $d<r$ by hypothesis.
\end{proof}

\begin{theorem}\label{thm:Gsurf} (cf. \cite[Theorem 5.22]{CCFMk2}) 
Let $X$ be a non-degenerate, projectively Cohen-Macaulay surface
of degree $r$ in $\PR$, $r\ge3$, which is connected in codimension
one. Then, any irreducible component of $X$ is either

\begin{itemize}
\item[(i)] a minimal Gorenstein surface, and there is at most one
such component, or \item[(ii)] a minimal degree surface.
\end{itemize}

If there is a component of type (i), then the
intersection in codimension one of any two distinct components can
be only a line.

If there is no component of type (i), then the
intersection in codimension one of any two distinct components is
either a line or a (possibly reducible) conic. Moreover, if two
components meet along a conic, all the other intersections are
lines.

Furthermore, $X$ is projectively Gorenstein if and only if either
\begin{itemize}
\item[(a)] $X$ is irreducible of type (i), or
\item[(b)] $X$ consists of only two components of type (ii) meeting along a conic, or
\item[(c)] $X$ consists of $\nu$, $3\le \nu \le r$, components of type (ii)
meeting along lines
and the dual graph $G_D$ of a general hyperplane section $D$ of $X$
is a cycle $E_\nu$.
\end{itemize}
\end{theorem}

\begin{proof}
Let $D$ be a general hyperplane section of $X$. Since $X$ is
projectively Cohen-Macaulay, it is arithmetically
Cohen-Macaulay (cf.\ Proposition \ref{prop:74}). This implies that $D$ is an arithmetically
Cohen-Macaulay (equiv.\ arithmetically normal) (equiv. arithmetically normal) curve (cf.\ Theorem \ref{thm:1ter} in Appendix \ref{S:pCMpG}). By Lemma
\ref{lem:pa=d}, $p_a(D)=1$. Therefore, for each connected subcurve
$D'$ of $D$, one has $0\le p_a(D')\le 1$ and there is at most one
irreducible component $D''$ with $p_a(D'')=1$. In particular two
connected subcurves of $D$ can meet at most in two points. This
implies that two irreducible components of $X$ meet either along a
line or along a conic. The linear normality of $X$ immediately
implies that each irreducible component is linearly normal too. As
a consequence of Theorem \ref{thm:DelP2} and of Lemma
\ref{lem:pa=d}, all this proves the statement about the components
of $X$ and their intersection in codimension one.

It remains to prove the final part of the statement.

If $X$ is irreducible, the assertion is trivial,
so assume $X$ reducible.

Suppose that all the intersections in codimension one
of the distinct components of $X$ are lines.
If either the dual graph $G_D$ of a general hyperplane section $D$ of $X$
is not a cycle or there is an irreducible component of $D$ which is not rational,
then $D$ is not Gorenstein
(see the discussion at the end of Example \ref{ex:zngraphs}),
contradicting the assumption that $X$ is Gorenstein.

Conversely, if $G_D$ is a cycle $E_{\nu}$ and each component of $D$ is rational,
then $D$ is projectively Gorenstein.
In particular, if all the components of $D$ are lines, then $D$
isomorphic to $C_{E_{\nu}}$ (cf.\ again Example \ref{ex:zngraphs}).
Therefore $X$ is projectively Gorenstein too (cf.\ Proposition \ref{lem:PG} in Appendix
\ref{S:pCMpG}).

Suppose that $X$ consists of two irreducible components meeting along a conic.
Then $D$ consists of two rational normal curves meeting at
two points; thus the dualizing sheaf $\omega_D$ is trivial, i.e.\ $D$ is
projectively Gorenstein and Gorenstein,
therefore so is $X$ (cf.\ Proposition \ref{lem:PG} in Appendix
\ref{S:pCMpG}).

Conversely, let us suppose that $X$ is projectively Gorenstein and there are two
irreducible components $X_1$ and $X_2$ meeting along a conic.
If there are other components, then there is a component $X'$ meeting all the
rest along a line.
Thus, the hyperplane section contains a rational curve meeting all the
rest at a point. Therefore the dualizing sheaf of $D$ is not trivial,
hence $D$ is not Gorenstein, thus $X$ is not Gorenstein.
\end{proof}

By using Theorems \ref{thm:DelP1}, \ref{thm:X} and
\ref{thm:DelP2}, we can prove the following result:

\begin{proposition}\label{prop:hypsect} (cf. \cite[Proposition 5.23]{CCFMk2}) 
Let $X$ be a non-degenerate surface in $\Pp^r$, for some $r$, and
let $n \ge 3$ be an integer.

\begin{itemize}

\item[(i)] If $r=n+1$ and if a hyperplane section of $X$ is
$C_{R_n}$, then either:

\begin{itemize}
\item[a.] $X$ is a smooth rational cubic scroll, possible only if
$n=3$, or

\item[b.] $X$ is a Zappatic surface, with $\nu$ irreducible
components of $X$ which are either planes or smooth quadrics,
meeting along lines, and the
Zappatic singularities of $X$ are $h\ge1$ points of type
$R_{m_i}$, $i=1,\ldots,h$, such that
\begin{equation}\label{eq:nu-2}
\sum_{i=1}^h (m_i-2) = \nu-2.
\end{equation}
In particular $X$ has global normal crossings if and only if
$\nu=2$, i.e.\ if and only if either $n=3$ and $X$ consists of a
plane and a quadric meeting along a line, or $n=4$ and $X$
consists of two quadrics meeting along a line.
\end{itemize}

\item[(ii)] If $r=n+1$ and if a hyperplane section of $X$ is
$C_{S_n}$, then either:

\begin{itemize}
\item[a.] $X$ is the union of a smooth rational normal scroll
$X_1=S(1,d-1)$ of degree $d$, $2\le d \le n$, and of $n-d$
disjoint planes each meeting $X_1$ along different lines of the
same ruling, in which case $X$ has global normal crossings; or

\item[b.] $X$ is planar Zappatic surface
with $h\ge1$ points of type $S_{m_i}$, $i=1,\ldots,h$, such that
\begin{equation}\label{eq:n-1su2}
\sum_{i=1}^h \binom{m_i-1}{2} =\binom{n-1}{2}.
\end{equation}
\end{itemize}

\item[(iii)] If $r=n$ and if a hyperplane section of
$X$ is $C_{E_n}$ then either:

\begin{itemize}
\item[a.] $X$ is an irreducible del Pezzo surface of degree $n$ in
$\Pp^n$, possible only if $n\le6$; in particular $X$ is smooth if
$n=6$; or

\item[b.] $X$ has two irreducible components $X_1$ and $X_2$, meeting
along a (possibly reducible) conic; $X_i$, $i=1,2$, is either
a smooth rational cubic scroll, or a quadric, or a plane;
in particular $X$ has global normal crossings
if $X_1\cap X_2$ is a smooth conic
and neither $X_1$ nor $X_2$ is a quadric cone;

\item[c.] $X$ is a Zappatic surface whose irreducible
components $X_1,\ldots,X_\nu$ of $X$ are either planes
or smooth quadrics.
Moreover $X$ has a unique $E_\nu$-point, and no other Zappatic singularity,
the singularities in codimension one being double lines.

\end{itemize}

\end{itemize}

\end{proposition}

\begin{proof}
(i)  According to Remark \ref{rem:pCM} and Theorem \ref{thm:X},
$X$ is connected in codimension one and is a union of minimal
degree surfaces meeting along lines. Since a hyperplane section is
a $C_{R_n}$, then each irreducible component $Y$ of $X$ has to contain some line and
therefore it is a rational normal scroll, or a plane. Furthermore
$Y$ has a hyperplane section which is a
connected subcurve of $C_{R_n}$. It is then clear that $Y$ is
either a plane, or a quadric or a smooth rational normal cubic scroll.

We claim that $Y$ cannot be a quadric cone.
In fact, in this case, the hyperplane sections of $Y$ consisting
of lines pass through the vertex $y\in Y$.
Since $Y\cap\overline{(X\setminus Y)}$ also consists of lines passing through $y$,
we see that no hyperplane section of $X$ is a $C_{R_n}$.

Reasoning similarly, one sees that if a component $Y$ of $X$ is a
smooth rational cubic scroll, then $Y$ is the only component of
$X$, i.e.\ $Y=X$, which proves statement a.

Suppose now that $X$ is reducible, so its components are either planes or smooth quadrics.
The dual graph $G_D$ of a general hyperplane section $D$ of $X$
is a chain of length $\nu$ and any connecting edge
corresponds to a double line of $X$.
Let $x\in X$ be a singular point and let $Y_1,\ldots,Y_m$ be the irreducible
components of $X$ containing $x$.
Let $G'$ be the subgraph of $G_D$ corresponding to $Y_1\cup\cdots\cup Y_m$.
Since $X$ is projectively Cohen-Macaulay, then clearly $G'$ is connected,
hence it is a chain.
This shows that $x$ is a Zappatic singularity of type $R_m$.

Finally we prove Formula \eqref{eq:nu-2}.
Suppose that the Zappatic singularities of $X$ are $h$ points
$x_1,\ldots,x_h$ of type $R_{m_1},\ldots,R_{m_h}$, respectively.
Notice that the hypothesis that a hyperplane section of $X$ is a $C_{R_n}$
implies that two double lines of $X$
lying on the same irreducible component have to meet at a point,
because they are either lines in a plane or fibres of different
rulings on a quadric.

\begin{figure}[ht]
  \[
  \xymatrix@C=30pt{%
  *=0{\bullet} \ar@{-}[rrrrrr] \ar@/^1pc/ @{--}[]+<2pt,1.5pt>;
  [rrr]+<-2pt,1.5pt>^{x_i}
   & *=0{\bullet} & *=0{\bullet} \ar@/_1pc/ @{--}[]+<2pt,-1pt>;
   [rrrr]-<2pt,1pt>_{x_{i+1}}
   & *=0{\bullet} & *=0{\bullet} & *=0{\bullet} & *=0{\bullet}
  }
  \]
  \caption{The points $x_i$ and $x_{i+1}$ share a common edge in the associated
graph $G_X$.}
  \label{fig:commonedge}
\end{figure}
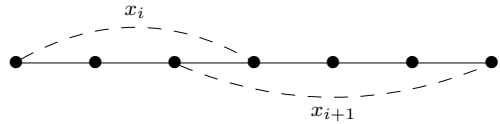
So the graph $G_X$ consists of $h$ open faces corresponding to the
points $x_i$, $1 \leq i \leq h$, and two contiguous open faces
must share a common edge, as shown in Figure \ref{fig:commonedge}.
Thus, both Formula \eqref{eq:nu-2} and the last part of statement
b.\ immediately follow.

\medskip
\noindent
(ii) Arguing as in the proof of (i), one sees that any irreducible component $Y$
of $X$ is either a plane, or a smooth quadric or a smooth
rational normal scroll with a line as a directrix.

If $Y$ is a rational normal scroll $S(1,d-1)$ of degree $d\ge 2$,
the subgraph of $S_n$ corresponding to the hyperplane section of
$Y$ is $S_d$. Then a.\ follows in this case, namely all the other
components of $X$ are planes meeting $Y$ along lines of the
ruling. Note that, since $X$ spans a $\Pp^{n+1}$, these planes are
pairwise skew and therefore $X$ has global normal crossings.

Suppose now that $X$ is a union of planes. Then $X$ consists of a
plane $\Pi$ and of $n-1$ more planes meeting $\Pi$ along distinct lines.
Arguing as in part (i), one sees that the planes different from $\Pi$
pairwise meet only at a point in $\Pi$.
Hence $X$ is smooth off $\Pi$.
On the other hand, it is clear that the singularities $x_i$ in $\Pi$
are Zappatic of type $S_{m_i}$, $i=1,\ldots,h$.
This corresponds to the fact that $m_i-1$ planes different from $\Pi$
pass through the same point $x_i\in\Pi$. Formula \eqref{eq:n-1su2}
follows by suitably counting the number of pairs of double lines in the
configuration.

\medskip
\noindent (iii) If $X$ is irreducible, then a.\ holds by elementary
properties of lines on a del Pezzo surface.

Suppose that $X$ is reducible. Every irreducible component $Y$ of
$X$ has a hyperplane section which is a stick curve strictly
contained in $C_{E_n}$. By an argument we already used in part
(i), then $Y$ is either a plane, or a quadric or a smooth
rational normal cubic scroll.

Suppose that an irreducible component $Y$ meets
$\overline{X\setminus Y}$ along a conic. Since $C_{E_n}$ is
projectively Gorenstein, then also $X$ is projectively Gorenstein (cf.\ Proposition \ref{lem:PG} in Appendix
\ref{S:pCMpG});
so, by Theorem \ref{thm:Gsurf}, $X$ consists of only two
irreducible components and b.\ follows.

Again by Theorem \ref{thm:Gsurf} and reasoning as in part (i),
one sees that all the irreducible components of $X$ are either planes or
smooth quadrics and
the dual graph $G_D$ of a general hyperplane section
$D$ of $X$ is a cycle $E_\nu$ of length $\nu$.

As we saw in part (i), two double lines of $X$ lying on the same
irreducible component $Y$ of $X$ meet at a point of $Y$.
Hence $X$ has some singularity besides the general points on the double lines.
Again, as we saw in part (i), such singularity can be either of type $R_m$
or of type $E_m$, where $R_m$ or $E_n$ are subgraphs of the dual graph $G_D$
of a general hyperplane section $D$ of $X$.
Since $X$ is projectively Gorenstein,
it has only Gorenstein singularities (cf.\ Remark \ref{rem:90} in Appendix \ref{S:pCMpG}),
in particular $R_m$-points are excluded.
Thus, the only singularity compatible with the above graph is an $E_\nu$-point.
\end{proof}

\begin{remark}
At the end of the proof of part $(iii)$, instead of using the
Gorenstein property, one can prove by a direct computation that a
surface $X$ of degree $n$, which is a union of planes and smooth
quadrics and such that the dual graph $G_D$
of a general hyperplane section $D$ of $X$ is a cycle of length $\nu$,
must have an $E_\nu$-point and no other Zappatic singularity
in order to span a $\Pp^n$.
\end{remark}

\begin{corollary}\label{cor:hypsect}
Let $\X\to\Delta$ be a degeneration of surfaces whose central fibre $X$
is Zappatic.
Let $x\in X$ be a $T_n$-point. Let $\X'$ be the
blow-up of $\X$ at $x$. Let $E$ be the exceptional divisor, let
$X'$ be the proper transform of $X$, $\Gamma=C_{T_n}$ be the
intersection curve of $E$ and $X'$. Then $E$ is a minimal degree
surface of degree $n$ in $\Pp^{n+1}=\Pp(T_{\X,x})$, and $\Gamma$
is one of its hyperplane sections.

In particular, if $x$ is either a $R_n$- or a $S_n$-point, then
$E$ is as described in Proposition \ref{prop:hypsect}.
\end{corollary}
\begin{proof}
The first part of the statement directly follows from Lemma
\ref{prop:mingminzap}, Proposition \ref{prop:mgmz3f} and Theorem
\ref{thm:X}.
\end{proof}

We close this section by stating a result which will be useful in
the sequel:

\begin{corollary}\label{cor:Gor}
Let $y$ be a point of a projective threefold $Y$. Let $H$ be an
effective  Cartier divisor on $Y$ passing through $y$. If $H$ has
an $E_n$-point at $y$, then $Y$ is Gorenstein at $y$.
\end{corollary}

\begin{proof}
Recall that $H$ is Gorenstein at $y$ (cf.\ Remark \ref{rem:pg})
and $H$ locally behaves as a hyperplane section of $Y$ at $y$
(cf.\ the proof of Proposition \ref{prop:min-quasimin}), therefore
$Y$ is Gorenstein at $y$ (cf.\ Theorem \ref{thm:1tris} in Appendix \ref{S:pCMpG}).
\end{proof}

\begin{remark}\label{rem:Gor}
Let $\X\to\D$ be a degeneration of surfaces whose central fibre $X$
is good Zappatic.
From Definition \ref{def:Zappsurf} and Corollary \ref{cor:Gor},
it follows that $\X$ is Gorenstein at all the points of X,
except at its $R_n$- and $S_n$-points.
\end{remark}

\section{Resolutions of the total space of a degeneration of surfaces to a Zappatic one}\label{S:resolution}

Given $\pi: \X \to \D$ a degeneration of surfaces with good Zappatic central fibre $X = \X_0$,
the aim of this section is to describe partial and total desingularizations
of the total space $\X$ of the degeneration. These will be fundamental tools in in Sections \ref{S:5} and \ref{S:BI}, where
we shall combinatorially compute
the $K^2$ of the smooth fibres of $\X$ (cf.\ Theorem \ref{thm:k2Gmain})
and prove the {\em Multiple Point Formula}
(cf.\ Theorem \ref{thm:BI}), respectively, as well as in
Section \ref{S:pgzappdeg}, where we shall compute the {\em
geometric genus} of the smooth fibres of $\X$ (cf.\ Theorem \ref{thm:pg}).

As recalled in Remark \ref{rem:pg}, a good Zappatic surface $X$ is
Gorenstein only at the $E_n$-points, for $n \geq 3$. Therefore,
when $X$ is the central fibre of a degeneration $\X \to \D$, from Remark \ref{rem:Gor}, then
also $\X$ is Gorenstein at the $E_n$-points of its central fibre. Thus one can first
consider a partial resolution of the total space $\X$ at the $R_n$- and $S_n$-points of $X$,
for $n \geq 3$, in order to make both the total space and the central fibre Gorenstein. More precisely,
one wants to produce a birational model of $\X$, denoted by
\begin{equation}\label{def:gorred}
\X^G \to \D,
\end{equation}such that:

\begin{itemize}
\item[(i)] $\X^G$ is isomorphic to $\X$ off the central fibre;
\item[(ii)] $\X^G$ is Gorenstein;
\item[(iii)] for each irreducible component $X_i$ of the central fibre $\X_0 = X$
of $\X$ there is some irreducible component of the central fibre
$\X^G_0 = X^G$ of $\X^G$ which dominates $X_i$, $1 \leq i \leq v$.
\end{itemize}$\X^G$ is said to be a {\em Gorenstein reduction} of $\X$. As it will be clear from the steps of Algorithm \ref{gr}, $X^G$ is a good Zappatic surface having only $E_n$-points as
Zappatic singularities, so it is Gorenstein. In particular, both the dualizing sheaves
$\omega_{\X^G}$ and $\omega_{X^G}$ will be invertible.

For our aims, we also need to completely resolve the total space of the degeneration.
In this case, as it will be shown in Algorithm \ref{ssr}, one can
get a {\em normal crossing reduction} of $\X$ (cf.\ Remark \ref{def:degen2}), say $\X^s$, such that:
\begin{itemize}
\item[(i)] $\X^s$ is smooth and it is isomorphic to $\X$ off the central fibre;
\item[(ii)] its central fibre $\X_0^s=X^s$ is a Zappatic surface whose support is with global
normal crossings, i.e.\ $X^s_{red}$ is with only $E_3$-points as Zappatic singularities;
\item[(iii)] for each irreducible component $X_i$ of the central fibre $\X_0 = X$
of $\X$ there is some irreducible component of the central fibre
$X^s$ which dominates $X_i$, $1 \leq i \leq v$.
\end{itemize}

Given a degeneration of surfaces $\pi: \X \to \D$, in order to determine both
a Gorenstein and a normal crossing reduction of $\X$
it is necessary to carefully analyze the process,
basically described in Chapter II of \cite{Kempf}, which produces the
semistable reduction.

As we said in Remark \ref{def:degen2},
Hironaka's result implies the existence of a normal crossing
reduction of $\pi$. The birational transformation involved in resolving the
singularities can be taken to be a sequence of blow-ups
(which one can arrange to be at isolated points and along smooth curves)
interspersed with normalization maps.  For general singularities such a
procedure may introduce components and double curves which can affect the computation of the
invariants of the central fibre, like e.g.\ the $\omega$-genus, etc..  Our next task is to show that,
under the assumption that the central fibre is good Zappatic, we have very
precise control over its invariants.  For this we will need to
more explicitly describe algorithms which produce Gorenstein and normal crossing reductions
of $\pi$, respectively. In order to do this, we will use,
as common in programming languages, the word ``while'' to indicate that
the statement following it is repeated until it becomes false.

\begin{algog}\label{gr}
Let $\X\to\D$ be a degeneration of surfaces with good Zappatic central fibre.
While $\X_0$ has a point $p$ of type either $R_n$ or $S_n$, $n\ge3$,
replace $\X$ by its blow-up at $p$.
\end{algog}

\begin{algo}\label{ssr} (cf. \cite[Algorithm 4.8]{CCFMpg}) 
Let $\X\to\D$ be a degeneration of surfaces with good Zappatic central fibre having
only $E_n$-points, $n \geq 3$, as Zappatic singularities.

\begin{enumerate}
\item[Step 1:] while $\X_0$ has a point $p$ of type $E_n$
and $\X$ has multiplicity $n\ge3$ at $p$,
replace $\X$ by its blow-up at $p$;

\item[Step 2:] while $\X$ has a double curve $\gamma$,
replace $\X$ by its blow-up along $\gamma$;

\item[Step 3:] if $\X$ has a double point $p$,
then replace $\X$ by the normalization of its blow-up at $p$
and go back to Step 2;

\item[Step 4:]
while there is a component of $\X_0$
with a double point $p$, replace $\X$ by its blow-up at $p$;

\item[Step 5:]
while there are two components $X_1$ and $X_2$ of $\X_0$
meeting along a curve with a node $p$,
first blow-up $\X$ at $p$, then
blow-up along the line which is the intersection
of the exceptional divisor with the proper transform of $\X_0$,
and finally replace $\X$ with the resulting threefold.

\end{enumerate}
\end{algo}

The following proposition is devoted to prove that the above algorithms work (cf. the proof of \cite[Theorem 6.1]{CCFMk2} and \cite[Proposition 4.9]{CCFMpg}).

\begin{proposition}\label{prop:ssr}
Let $\pi:\X\to\D$ be a degeneration of surfaces
with good Zappatic central fibre $X=\X_0=\bigcup_{i=1}^v X_i$.

\begin{itemize}
\item[(1)] Run the Gorenstein reduction algorithm \ref{gr};
the algorithm stops after finitely many steps
and its output gives
a Gorenstein reduction $\pi^G : \X^G \to\D$ of $\pi$.

\item[(2)] Consider $\pi : \X \to\D$ such that $X$ has only $E_n$-points, $n \geq 3$, as Zappatic
singularities and run
the normal crossing reduction algorithm \ref{ssr}.
The algorithm stops after finitely many steps
and its output gives
a normal crossing reduction $\bar\pi:\bar\X\to\D$ of $\pi$.
\end{itemize}
\end{proposition}

\begin{proof} (1) As observed in \S\ \ref{S:zapdeg}, since $\pi: \X \to \D$ is an arbitrary
degeneration of smooth surfaces to a good Zappatic one, the total space $\X$ of $\pi$ may have the following singularities:
\begin{itemize}
\item double curves, which are double curves also for $X$;
\item isolated double points along the double curves of $X$;
\item further singular points at the Zappatic
singularities of $X$, which can be isolated or
may occur on double curves of the total space.
\end{itemize}
Our aim is to prove that the Gorenstein reduction algorithm \ref{gr}
produces a total space which is Gorenstein and a central fibre
which has only $E_n$-points as Zappatic singularities.

\bigskip
\noindent
By Proposition \ref{prop:mgmz3f},
if $X$ has either a $R_n$-point or a $S_n$-point, $n\geq3$,
then the total space $\X$ has multiplicity $n$ at $p$.
Let $\X'\to\X$ be the blow-up of $\X$ at a $R_n$-point [resp.\ $S_n$-point] $p$.
By Proposition \ref{prop:hypsect}, the exceptional divisor $E$ is a Zappatic surface of degree $n$ in $\Pp^{n+1}$
such that all of its irreducible components are rational normal surfaces
meeting along lines and $E$ has at most $R_m$-points, $m\le n$
[resp.\ $S_m$-points, $m\le n$] as Zappatic singularities.
Let $X'$ be the proper transform of $X$.
The curve $\Gamma=E\cap X'$ is a stick curve
$C_{R_n}$ [resp.\ $C_{S_n}$] which, being nodal, does not contain
any Zappatic singularity of $E$.
The new central fibre $E\cup X'$ has either $E_3$- or $E_4$-points
at the double points of $\Gamma$, depending on whether
$E$ is smooth or has a double point there.
These points are accordingly either smooth or double points for $\X'$.

The fact that the process in $(1)$ is repeated finitely many times follows
e.g.\ from Proposition 3.4.13 in \cite{Kollar}.
If $\X^G \to\X$ is the composition of all the blow-ups done in $(1)$,
then $\pi^G: \X^G \to\D$ is a degeneration whose central fibre is a
Zappatic surface with only $E_n$-points, $n\ge3$, as Zappatic singularities. Thus $\pi^G: \X^G \to\D$ is a Gorenstein reduction of $\pi$ (this is clear for the
double points, for the $E_n$-points of the central fibre, see Corollary \ref{cor:Gor} and Remark
\ref{rem:Gor}).

\bigskip
\noindent (2) Our aim is now to prove that, given $\pi : \X \to \D$ as in $(2)$,
the normal crossing reduction algorithm \ref{ssr}  resolves the singularities of the total space $\X$
and produces a central fibre whose support has global normal crossings.

In general, the degeneration $\pi: \X \to \D$ we start with here in $(2)$
will be the output of the Gorenstein reduction algorithm \ref{gr}
applied to an arbitrary degeneration of smooth surfaces to a good
Zappatic one.

\bigskip \noindent
(Step 1) By Proposition \ref{prop:mgmz3f},
if $X$ has an $E_n$-point $p$, $n\ge3$, then
either $\X$ has multiplicity $n$ at $p$, or $n\leq4$
and $\X$ has at most a double point at $p$.
In this step we consider only the former possibility,
since the other cases are considered in the next steps.
Let $\X'\to\X$ be the blow-up of $\X$ at $p$.
By Proposition \ref{prop:hypsect},
the exceptional divisor $E$ is a Gorenstein surface of degree $n$
in $\Pp^{n}$ which is one of the following:
\begin{enumerate}
\item[(I)] an irreducible del Pezzo surface, possible only if $n\leq 6$;
\item[(II)] a union $F=F_1\cup F_2$ of two irreducible components $F_1$ and $F_2$
such that $F_1\cap F_2$ is a (possibly reducible) conic;
the surface $F_i$, $i=1,2$, is either a smooth rational normal cubic scroll,
or a quadric, or a plane;
\item[(III)] a Zappatic surface, whose $m\leq n$ irreducible components
meet along lines and are either planes or smooth quadrics;
moreover $E$ has a unique Zappatic singularity, which is an $E_m$-point.
\end{enumerate}
In case (I), the del Pezzo surface $E$ has at most isolated rational double points.

In case (II), the surface $E$ is Zappatic unless either the conic is reducible
or one of the two components is a quadric cone.
Note that, if $F_1\cap F_2$ is a conic with a double point $p'$, then
$F_1$ and $F_2$ are tangent at $p'$ and $E$ has not normal crossings.

Let $X'$ be the proper transform of $X$.
The curve $\Gamma=E\cap X'$ is a stick curve $C_{E_n}$.
In case (II), if an irreducible component of $E$ is a quadric cone,
the vertex of the cone is a double point of $\Gamma$ and
$\X'$ also has a double point there.
In case (III), the curve $\Gamma$, being nodal, does not contain the $E_m$-point of $E$.
As in $(1)$, one sees that the singular points of $\Gamma$
are either smooth or double points for $\X'$.

In cases (I) and (II), we have eliminated the original Zappatic $E_n$
singularity; in case (III), we have a single $E_m$ ($m \leq n$)
point to still consider.  Whatever extra double points have
been introduced, will be handled in later steps.

As in $(1)$, also Step 1 is repeated
finitely many times e.g.\ by Proposition 3.4.13 in \cite{Kollar}.

\bigskip
\noindent
(Step 2) Now the total space $\X$
of the degeneration has at most double points.
Suppose that $\X$ is singular in dimension one and
let $\gamma$ be an irreducible curve which is double for $\X$.
Then $\gamma$ lies in the intersection
of two irreducible components $X_1$ and $X_2$ of $X$.
By Definition \ref{def:Zappsurf} of Zappatic surface
and the previous steps, one has that $\gamma$ is smooth and the intersection of $X_1$ and $X_2$ is transversal
at the general point of $\gamma$.

Now let $\X'\to\X$ be the blow-up of $\X$ along $\gamma$.
Let $E$ be the exceptional divisor
and $X'_i$, $i=1,2$, be the proper transform of $X_i$ in $\X'$.
Let $p$ be the general point of $\gamma$.
Note that there are effective Cartier divisors of $X$ through $p$
having a node at $p$. Therefore there are effective Cartier divisors
of $\X$ through $p$ having at $p$ a double point of type $A_k$, for some $k\geq 1$.
Since the exceptional divisor of a minimal resolution of such a point
does not contain multiple components, we see that $E$ must be reduced.
Then $E$ is a conic bundle and $\gamma_i=E\cap X'_i$, $i=1,2$,
is a section of $E$ isomorphic to $\gamma$.

Let $C$ be the general ruling of $E$.
If $C$ is irreducible, then $E$ is irreducible and has at most isolated
double points.
We remark moreover that $\gamma_1$ and $\gamma_2$ are generically smooth
for the total space $\X'$, since they are generically smooth for $E$,
which is a Cartier divisor of $\X'$.
In this case, we got rid of the double curve.

Let $C=r_1\cup r_2$ be reducible into two distinct lines.
We may assume that $r_i\cap\gamma_i$, $i=1,2$, is a point
whereas $r_i\cap\gamma_{3-i}=\emptyset$.
This implies that $E$ is reducible;
one component meets $X'_1$ and the other meets $X'_2$.
Hence we may write $E=F_1\cup F_2$, where $F_i$ meets
generically transversally $X'_i$
along $\gamma_i$, $i=1,2$.
It may happen that $F_1$ and $F_2$ meet, generically transversally,
along finitely many fibres of their rulings; away from these,
they meet along the curve $\gamma'$, whose general point is $r_1\cap r_2$.

We note that $\gamma'$, being isomorphic to $\gamma$, is smooth.
Moreover, a local computation shows that $F_1$ and $F_2$
meet transversally at a general point of $\gamma'$.
If the general point of $\gamma'$ is smooth for $\X'$, we have nothing to do
with $\gamma'$, otherwise we go on blowing-up $\X'$ along $\gamma'$.
As usual, after finitely many blow-ups we get rid of all the curves
which are double for the total space.

\bigskip
\noindent
(Step 3)
Now the total space $\X$ of the degeneration has at most isolated double points.
Let $X_\red$ be the support of the central fibre $X$.
Note that, the first time one reaches this step, one has that $X_{\red}=X$,
which implies that $X_\red$ is Cartier.  In what follows, we only require
that in a neighborhood of the singular points where we apply this step,
the reduced set of components is Cartier.

By the discussion of the previous steps,
one sees that
a double point $p$ of $\X$ can be of the following types
(cf.\ Figure \ref{fig:doublepoints}):
\begin{enumerate}
\item[$(a)$] an isolated double point of $X_\red$;
\item[$(b)$] a point of a
double curve of $X_\red$;
\item[$(c)$] an $E_3$-point of $X_\red$;
\item[$(d)$] an $E_4$-point of $X_\red$;
\item[$(e)$] a quadruple point of $X_\red$
which lies in the intersection of three irreducible components
$X_1$, $X_2$ and $X_3$ of $X_\red$; two of them, say $X_2$ and
$X_3$, are smooth at $p$, whereas $X_1$ has a rational double
point of type $A_k$, $k\ge1$, at $p$. In this case, $X_2\cup X_3$
and $X_1$ are both complete intersection of $\X$ locally at $p$.
\end{enumerate}

\begin{figure}[ht]
\[
\begin{array}{ccccc}
\;\;
\raisebox{12pt}{%
\begin{xy}
0; <40pt,0pt>: 
(0,0)*_!/8pt/{p} *{\bullet} *\cir<6pt>{dr^ur} ,
(0,0);(0.13,-0.075)**@{-} ,
(0,0);(-0.13,-0.075)**@{-} ,
\end{xy}%
}
\;\; & \;\;
\begin{xy}
0; <40pt,0pt>: 
(0,-0.5);(0,1)**@{-} ?*{\bullet} ?*_!/7pt/{p} ,
\end{xy}
\;\; & \;\;
\begin{xy}
0; <40pt,0pt>: 
{\xypolygon3{~>{}}},
"1";"0"*{\bullet}**@{-} ?*_!/-6pt/{\gamma_1} ,
"2";"0"**@{-} ?*_!/6pt/{\gamma_2} ,
"3";"0"**@{-} ?*_!/-6pt/{\gamma_3} ,
\end{xy}
\;\; & \;\;
\raisebox{10pt}{%
\begin{xy}
0; <40pt,0pt>: 
{\xypolygon4{~>{}}},
"1";"0"*{\bullet}**@{-} ?*_!/7pt/{\gamma_1} ,
"2";"0"**@{-} ?*_!/-6pt/{\gamma_2} ,
"3";"0"**@{-} ?*_!/6pt/{\gamma_3} ,
"4";"0"**@{-} ?*_!/-6pt/{\gamma_4} ,
\end{xy}%
}
\;\; & \;\;
\begin{xy}
0; <40pt,0pt>: 
{\xypolygon3{~>{}}},
"1";"0"*{\bullet}**@{-} ?*_!/-6pt/{\gamma_1} ,
"2";"0"**@{-} ?*_!/6pt/{\gamma_2} ,
"3";"0"**@{-} ?*_!/-6pt/{\gamma_3} ,
"2";"3"**@{}?*{X_1} ,
"0"*\cir<9pt>{dr^ur} ,
\end{xy}
\;\;
\\[5ex]
(a) & (b) & (c)  & (d) & (e)
\end{array}
\]
\vspace{-3ex}
\caption{Types of double points of the total space $\X$.}\label{fig:doublepoints}
\end{figure}
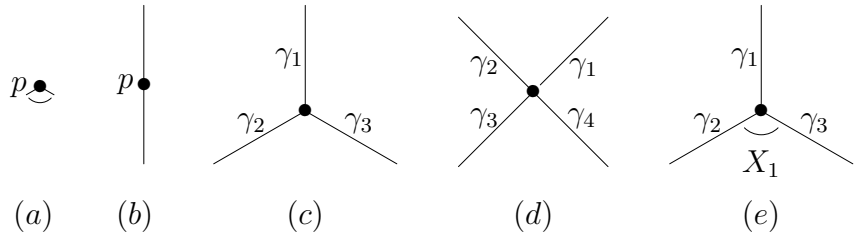

Double points of type (a) may appear either in Step 1, if the exceptional divisor
is a singular del Pezzo surface, or in Step 2, if the exceptional divisor
is a singular conic bundle.
In both cases, they are rational double points for $X_\red$.
By resolving them, one clearly gets as exceptional divisors
only rational surfaces meeting each other (and the proper transform
of the central fibre) along rational curves.

Consider a double point $p$ of type (b), so $p$ lies on a double
curve which is in the intersection of two irreducible components
$X_1$ and $X_2$ of $X_\red$. Let $\X'\to\X$ be the blow-up of $\X$
at $p$ and let $E$ be the exceptional divisor, which is a quadric
surface in $\Pp^3$. Denote by $X'_i$ the proper transform of
$X_i$, $i=1,2$, and by $p'$ the point $p'=E\cap X'_1\cap X'_2$.
Since a general hyperplane section of $X_1\cup X_2$ at $p$ is a
curve with a node at $p$, the quadric $E$ is either:
\begin{enumerate}
\item[(i)] a smooth quadric meeting $X'_i$, $i=1,2$, along a line; or
\item[(ii)] an irreducible quadric cone with vertex $p'$; or
\item[(iii)] the union of two distinct planes meeting along a line $\gamma$ passing through $p'$.
\end{enumerate}
In case (i), we resolved the singularity of the total space at $p$.
In case (ii), the new total space $\X'$ has an isolated double point of type (e) at $p'$.
In case (iii), there are two possibilities: if the line $\gamma$ is a double
curve of $\X'$, then we go back to Step 2, otherwise $\X'$ has an isolated
double point of type (d).

Let $p$ be a double point of type (c).
According to Proposition \ref{prop:mgmz3f},
the embedding dimension of $\X$ at $p$ is 4 and the central fibre is locally analytically near $p$
a hyperplane section of $\X$. Since the multiplicity of the singularity of the threefold is two,
and the multiplicity of the central fibre at this point is three,
the locally analytic hyperplane section must contain a component
of the tangent cone of the threefold singularity.  This tangent
cone is therefore a quadric which has rank at most two:
it is either two distinct hyperplanes
or a double hyperplane (i.e.\ a hyperplane counted twice).
In fact a local computation shows that the latter cannot happen.
In the former case, when one blows up $\X$ at $p$,
one introduces two planes in the new central fibre.
One of these planes meets the proper transforms of the three
components each in a line, forming a triangle in that plane;
this plane is double in the new central fibre.
(Note that at this point we introduce a non-reduced component of the
central fibre; but the rest of the algorithm does not involve this
multiple component.)
The other of the planes, which is simple in the new central fibre,
meets each of the proper transforms at a single distinct point,
which is still an ordinary double point of the total space.
Three more blow-ups, one each at these double points, locally resolve the
total space.  (This analysis follows from a local computation.)

Consider now a double point $p$ of type (d).
By Proposition \ref{prop:mgmz3f},
locally the tangent cone of $\X$ at $p$ is a quadric cone in $\Pp^3$
and the tangent cone $T$ of $X_\red$ at $p$ is obtained by cutting it
with another quadric cone in $\Pp^3$, hence
$T$ is a cone in $\Pp^4$ over a reduced, projectively normal curve of degree 4
and arithmetic genus $1$ which spans a $\Pp^3$.
Let $\X'\to\X$ be the blow-up of $\X$ at $p$
and let $E$ be the exceptional divisor.
Then $E$ is a quadric meeting the proper transform of $X$
along a stick curve $C_{E_4}$, therefore $E$ is either
\begin{enumerate}
\item[(i)] a smooth quadric; or
\item[(ii)] the union of two distinct planes meeting along a line $\gamma$.
\end{enumerate}
In case (i), we resolved the singularity of $\X$ at $p$.
In case (ii), there are two possibilities: if the line $\gamma$ is a double
curve of $\X'$, then we go back to Step 2, otherwise $\X'$ has again two isolated
double points of type (d) at the intersection of $\gamma$ with the proper
transform of $X_\red$.

Here we have created double components of the central fibre, namely
the exceptional divisor is counted twice.
However this exceptional divisor is a Cartier divisor, and therefore
$X_\red$ is also a Cartier divisor locally near this exceptional divisor.

Finally let $p=X_1\cap X_2\cap X_3$ be a double point of type (e).
As in the case of type (d),
locally the tangent cone of $\X$ at $p$ is a quadric cone in $\Pp^3$,
whereas the tangent cone $T$ of $X_\red$ at $p$ is a cone in $\Pp^4$
over a reduced, projectively normal curve of degree 4
and arithmetic genus $1$ which spans a $\Pp^3$.
Let $\X'\to\X$ be the blow-up of $\X$ at $p$
and let $E$ the exceptional divisor.
Denote by $p'$ the intersection of $E$ with the proper transform of $X_2$ and $X_3$.
Then $E$ is a quadric meeting the proper transform of $X_\red$
along the union of two lines and a conic spanning a $\Pp^3$,
therefore $E$ is either
\begin{enumerate}
\item[(i)] a smooth quadric; or
\item[(ii)] a quadric cone with vertex at $p'$; or
\item[(iii)] a pair of planes.
\end{enumerate}
In case (i), we resolved the singularity of the total space at $p$.
In case (ii), the total space $\X'$ has at $p'$ again a point of type (e).
More precisely, if $p$ is a rational double point of type $A_k$,
then $p'$ is a rational double point of type $A_{k-1}$ for $E$.
In case (iii), the line of intersection of the two planes may be
singular for the new total space; if so, we return to Step 2.
If not, there are again isolated double points of type (d) and we
iterate this step again.

As in the case of type (d), the reduced central fibre remains Cartier
in a neighborhood of the new exceptional locus.

It is clear that, after having repeated finitely many times
Steps 2 and 3, one resolves the singularities of the total space
at the double points of these five types (a)-(e).

We remark that we can proceed, in Step 3, by first resolving all of the
points of type (c), and that such points are not created in the resolutions
of points of type (d) and (e).  In fact they are not created in any later
step of the algorithm.
Indeed, anytime three components $X_1$, $X_2$, and $X_3$ concur at a point
as in type (c) where at least one of the three surfaces has been created
by blowing-up, we claim that exactly one of the three surfaces has been
created by blowing-up (i.e., is an exceptional divisor). Since such an
exceptional divisor is locally Cartier and smooth at the point,
then the total space is smooth at the point and therefore the point
cannot be of type (c).
To prove the claim, note that the only other possibility is that
two of the three components, say $X_2$ and $X_3$ belong to an exceptional divisor.
By blowing them down, then $X_1$ acquires a singular point which is worse
than an ordinary double point, which is impossible.

\bigskip
\noindent
(Step 4)
Let $p$ be an isolated double point of the central fibre $X$
which is a smooth point of $\X$.
According to the previous steps, $p$ is either a rational
double point of a del Pezzo surface
or the singular point of a reduced fibre of a conic bundle.
In both cases, the singularity of $X$ at $p$ is resolved
by finitely many blow-ups.
Since $p$ is a smooth point of $\X$, the exceptional divisor
of each blow-up is a plane.

\bigskip
\noindent
(Step 5)
Following the previous steps, one sees that
the support of the central fibre $X$ has global normal crossings,
except at the points $p$, where
two components $X_1$ and $X_2$ of $X$ meet along a curve
with a node at $p$.
Note that $X_1$ and $X_2$ are indeed tangent at $p$.

If one blows-up $\X$ at $p$, the exceptional divisor $E$ is a plane
meeting the proper transform $X'_i$ of $X_i$, $i=1,2$, along a line $\gamma$,
which is a $(-1)$-curve both on $X'_1$ and $X'_2$.
The support of the new central fibre has not yet normal crossings.
However a further blow-up along $\gamma$ produces the normal crossing reduction.
\end{proof}

\section{Combinatorial computation of $K^2$}\label{S:5}

The results contained in Sections \ \ref{S:MS} and \ \ref{S:resolution} will be used in this
section to prove combinatorial formulas for $K^2 = K^2_{\X_t}$, where $\X_t$ is a smooth surface which degenerates to a good Zappatic surface $\X_0 =X=\bigcup_{i=1}^v X_i$, i.e.\ $\X_t$ is
the general fibre of a degeneration of surfaces whose central
fibre is good Zappatic (cf.\ Notation \ref{not:zappdeg}).

Indeed, by using the combinatorial data associated to $X$ and
$G_X$ (cf.\ Definition \ref{def:dualgraph} and Notation
\ref{def:vefg}), we shall prove the following main result:

\begin{theorem}\label{thm:k2Gmain} (cf. \cite[Theorem 6.1]{CCFMk2}) 
Let $\X\to\D$ be a degeneration of surfaces whose central fibre
is a good Zappatic surface $X=\X_0=\bigcup_{i=1}^v X_i$.
Let $C_{ij} = X_i \cap X_j$ be a double curve of $X$, which is considered as a
curve on $X_i$, for $1 \leq i \neq j \leq v$.

If $K^2 := K^2_{\X_t}$, for $t\ne0$, then (cf.\ Notation
\ref{def:vefg}):
\begin{equation}\label{eq:K2Gbounds}
K^2  =  \sum_{i=1}^v \left( K_{X_i}^2 + \sum_{j \ne i} (4 g_{ij} -
C_{ij}^2)\right)-8e+ \sum_{n \ge 3} 2n f_n+r_3+ k,
\end{equation}
where $k$ depends only on the presence of $R_n$- and $S_n$-points,
for $n\ge 4$, and precisely:
\begin{equation} \label{eq:K2Gbounds2}
\sum_{n \geq 4} (n-2)(r_n +s_n) \le k
 \le \sum_{n \ge 4} \left( (2n-5) r_n + \binom{n-1}{2} s_n \right).
\end{equation}
\end{theorem}

In case $\X$ is an embedded degeneration and $X$ is also planar, we have the following:

\begin{corollary}\label{cor:k2Gmain}
Let $\X\to\D$ be an embedded degeneration of surfaces whose central fibre
is a good planar Zappatic surface $X=\X_0=\bigcup_{i=1}^v \Pi_i$.
Then:
\begin{equation}\label{eq:pK2Gbounds}
K^2 = 9v-10e+\sum_{n \ge 3} 2n f_n+r_3+k
\end{equation}
where $k$ is as in \eqref{eq:K2Gbounds2} and depends only on the
presence of $R_n$- and $S_n$-points, for $n\ge 4$.
\end{corollary}

\begin{proof}
Clearly $g_{ij}=0$, for each $1 \leq i \neq j \leq v$, whereas
$C_{ij}^2=1$, for each pair $(i,j)$ s.t. $e_{ij} \in E$, otherwise
$C^2_{ij}=0$.
\end{proof}

The proof of Theorem \ref{thm:k2Gmain} will be done in several
steps. The first one is to compute $K^2$ when $X$ has only
$E_n$-points. In this case, and only in this case, $K_X$ is a
Cartier divisor (cf.\ Remark \ref{rem:pg}).

\begin{theorem}\label{thm:k2Gnc} (cf. \cite[Theorem 6.6]{CCFMk2}) 
Under the assumptions of Theorem \ref{thm:k2Gmain}, if
$X=\bigcup_{i=1}^v X_i$ has only $E_n$-points, for $n\ge3$, then:
\begin{equation}\label{eq:k2Gnc}
K^2  =  \sum_{i=1}^v \left( K_{X_i}^2 + \sum_{j \ne i} (4 g_{ij} -
C_{ij}^2) \right)-8e+ \sum_{n \ge 3} 2n f_n.
\end{equation}
\end{theorem}

\begin{proof} Recall that, in this case, the total space $\X$ is Gorenstein (cf.\ Remark \ref{rem:Gor}).
Thus, $K_{\X}$ is a Cartier divisor on $\X$. Therefore $K_X$ is also
Cartier and it makes sense to consider $K_X^2$ and the adjunction
formula states $ K_X = (K_{\X}+X)_{|X}. $

We claim that
\begin{equation}\label{eq:k|Xi}
{K_X}_{|X_i} =(K_{\X}+ X)_{|X_i}  = K_{X_i} + C_i,
\end{equation}
where $C_i = \sum_{j\neq i}C_{ij}$ is the union of the double
curves of $X$ lying on the irreducible component $X_i$, for each
$1 \leq i \leq v$. Since $\Oc_X(K_X)$ is invertible, it suffices
to prove \eqref{eq:k|Xi} off the $E_n$-points. In other words, we
can consider the surfaces $X_i$ as if they were Cartier divisors
on $\X$. Then, we have:
\begin{equation}\label{eq:k|Xibis}
{K_X}_{|X_i} =(K_{\X}+ X)_{|X_i}  = \big(K_{\X} + X_i + \sum_{j
\neq i} X_j\big)_{|X_i} = K_{X_i} + C_i,
\end{equation}as we had to show.
Furthermore:
\begin{align}
K^2 & = (K_{\X}+ \X_t)^2 \cdot \X_t= (K_{\X}+X)^2\cdot X
  = (K_{\X}+X)^2 \cdot \sum_{i=1}^v X_i
  =  \sum_{i=1}^v \left((K_{\X}+X)_{|X_i}\right)^2 =    \notag\\
& = \sum_{i=1}^v (K_{X_i}^2 + 2 C_i K_{X_i} + C^2_i)
  = \sum_{i=1}^v K_{X_i}^2 + \sum_{i=1}^v C_i K_{X_i}
    + \sum_{i=1}^v C_i (C_i + K_{X_i}) =    \notag\\
& =  \sum_{i=1}^v K_{X_i}^2 + \sum_{i=1}^v (\sum_{j \ne i} C_{ij})
K_{X_i}
    + \sum_{i=1}^v 2(p_a(C_i) - 1). \label{eq:k2vai}
\end{align}

As in Notation \ref{def:vefg}, $C_{ij} =
\sum_{t=1}^{h_{ij}}C_{ij}^t$ is the sum of its disjoint, smooth,
irreducible components, where $h_{ij}$ is the number of these
components. Thus, $$C_{ij} K_{X_i} = \sum_{t=1}^{h_{ij}} (C_{ij}^t
K_{X_i}),$$for each $1 \leq i \neq j \leq v$. If we denote by
$g_{ij}^t$ the geometric genus of the smooth, irreducible curve
$C_{ij}^t$, by the Adjunction Formula on each $C_{ij}^t$, we have
the following intersection number on the surface $X_i$:
\[
C_{ij} K_{X_i} = \sum_{t=1}^{h_{ij}} (2g_{ij}^t - 2 -
(C_{ij}^t)^2) = 2 g_{ij} - 2 h_{ij} - C_{ij}^2,
\]
where the last equality follows from the definition of geometric genus of $C_{ij}$ and the
fact that $C_{ij}^s C_{ij}^t = 0$, for any $1 \leq t \neq s \leq
h_{ij}$.

Therefore, by the distributivity of the intersection form and by
\eqref{eq:k2vai}, we get:
\begin{align}\label{eq:k2G}
K^2 =  \sum_{i=1}^v K_{X_i}^2 + \sum_{i=1}^v (\sum_{j \ne i} (2
g_{ij}- 2 h_{ij}) - C_{ij}^2)
    + \sum_{i=1}^v 2(p_a(C_i) - 1).
\end{align}

For each index $i$, consider now the normalization $\nu_i :
\tilde{C}_i \to C_i$ of the curve $C_i$ lying on $X_i$; this
determines the short exact sequence:
\begin{equation}\label{eq:norm}
0 \to \Oc_{C_i} \to (\nu_i)_* (\Oc_{\tilde{C}_i}) \to
\underline{t}_i \to 0,
\end{equation}
where $\underline{t}_i$ is a sky-scraper sheaf supported on
$\Sing(C_i)$, as a curve in $X_i$. By using Notation
\ref{def:vefg}, the long exact sequence in cohomology induced by
\eqref{eq:norm} gives that:
\[
\chi(\Oc_{C_i}) + h^0(\underline{t}_i)= \sum_{j \neq i}
\sum_{t=1}^{h_{ij}} \chi (\Oc_{C_{ij}^t})=
  \sum_{j \neq i} (h_{ij} - g_{ij}).
\]
Since $\chi (\Oc_{C_i}) = 1 -p_a(C_i)$, we get
\begin{equation}\label{eq:norm2}
p_a(C_i) -1=  \sum_{j \neq i} (g_{ij} - h_{ij}) +
h^0(\underline{t}_i), \; 1 \leq i \leq v.
\end{equation}

By plugging Formula \eqref{eq:norm2} in \eqref{eq:k2G}, we get:
\begin{align}\label{eq:norm2bis}
K^2 =  \sum_{i=1}^v \left(K_{X_i}^2 + \sum_{j \ne i} (4 g_{ij} -
C_{ij}^2)\right)
    - 8 e +  2 \sum_{i=1}^v h^0(\underline{t}_i) .
\end{align}

To complete the proof, we need to compute $h^0(\underline{t}_i)$.
By definition of $\underline{t_i}$, this computation is a local
problem. Suppose that $p$ is an $E_n$-point of $X$ lying on $X_i$,
for some $i$. By the very definition of $E_n$-point (cf.\
Definition \ref{def:zappaticsing} and Example \ref{ex:zngraphs}),
$p$ is a node for the curve $C_i\subset X_i$; therefore
$h^0(\underline{t}_{i|p})= 1$. The same holds on each of the other
$n-1$ curves $C_j \subset X_j$, $1 \leq j \neq i \leq n$,
concurring at the $E_n$-point $p$. Therefore, by
\eqref{eq:norm2bis}, we get \eqref{eq:k2Gnc}.
\end{proof}

\begin{proof}[Proof of Theorem \ref{thm:k2Gmain}]
The previous argument proves that, in this more general case, one
has:
\begin{equation}\label{eq:k2vaic}
K^2   =  \sum_{i=1}^v  \left(K_{X_i}^2 + \sum_{j \ne i} (4 g_{ij}
- C_{ij}^2)\right) -8e + 2 \sum_{i=1}^v h^0(\underline{t}_i) - c
\end{equation}where $c$ is a positive correction term which
depends only on the points where $\X$ is not Gorenstein, i.e.\ at
the $R_n$- and $S_n$-points of its central fibre $X$.

To prove the statement, we have to compute:

\begin{itemize}
\item[(i)] the contribution of $h^0(\underline{t}_i)$ given by the
$R_n$- and the $S_n$-points of $X$, for each $1 \leq i \leq v$;
\item[(ii)] the correction term $c$.
\end{itemize}

For  (i), suppose first that $p$ is a $R_n$-point of
$X$ and let $C_i$ be one of the curves passing through $p$. By
definition (cf.\ Example \ref{ex:tngraphs}), the point $p$ is
either a smooth point or a node for $C_i\subset X_i$. In the first
case we have $h^0(\underline{t}_{i|p})= 0$ whereas, in the latter,
$h^0(\underline{t}_{i|p})= 1$. More precisely, among the $n$
indexes involved in the $R_n$-point there are exactly two indexes,
say $i_1$ and $i_n$, such that $C_{i_j}$ is smooth at $p$, for $j
=1$ and $j=n$, and $n-2$ indexes such that $C_{i_j}$ has a node at
$p$, for $2 \leq j \leq n-1$. On the other hand, if we assume that
$p$ is a $S_n$-point, then $p$ is an ordinary $(n-1)$-tuple point
for only one of the curves concurring at $p$, say $C_i\subset
X_i$, and a simple point for all the other curves $C_j\subset
X_j$, $1 \leq j \neq i \leq n$. Recall that an ordinary
$(n-1)$-tuple point contributes $\binom{n-1}{2}$ to
$h^0(\underline{t}_i)$.

Therefore, from \eqref{eq:k2vaic}, we have:
\[
K^2  =  \sum_{i=1}^v  \left(K_{X_i}^2 + \sum_{j \ne i} (4 g_{ij} -
C_{ij}^2)\right) -8e + \sum_{n\ge3} 2nf_n + \sum_{n\ge3} 2(n-2)r_n
+ \sum_{n\ge4} (n-1)(n-2)s_n - c.
\]

In order to compute the correction term $c$, we have to perform a
partial resolution of $\X$ at the $R_n$- and $S_n$-points of $X$,
which makes the total space Gorenstein; i.e.\ we have to consider a Gorenstein reduction
of the degeneration  $\X \to \D$.This will give us
\eqref{eq:K2Gbounds}, i.e.\ $$K^2  =  \sum_{i=1}^v  \left(K_{X_i}^2
+ \sum_{j \ne i} (4 g_{ij} - C_{ij}^2)\right)- 8 e +\sum_{n \ge 3}
2n f_n+r_3+k,$$where
\[
k:= \sum_{n\ge3} 2(n-2)r_n  \; - r_3 \; + \sum_{n\ge4}
(n-1)(n-2)s_n  \; - c.
\]

From Algorithm \ref{gr} and Proposition \ref{prop:ssr} - $(1)$, we know that
$\X \to \D$ admits a Gorenstein reduction; we now consider a detailed analysis of
this Gorenstein reduction in order to compute the contribution $c$.
It is clear that the contribution to $c$ of each such point is
purely local. In other words,
\[
c=\sum_x c_x
\]
where $x$ varies in the set of $R_n$- and $S_n$-points of $X$ and
where $c_x$ is the contribution at $x$ to the computation of $K^2$
as above.

In the next Proposition \ref{prop:cx}, we shall compute such local
contributions. This result, together with Theorem \ref{thm:k2Gnc},
will conclude the proof.
\end{proof}

\begin{proposition}\label{prop:cx}
In the hypothesis of Theorem \ref{thm:k2Gmain},
if $x\in X$ is a $R_n$-point then:
\[
n-2 \ge c_x \ge 1,
\]
whereas if $x\in X$ is a $S_n$-point then:
\[
(n-2)^2 \ge c_x \ge \binom{n-1}{2}.
\]
\end{proposition}

\begin{proof}
Since the problem is local, we may (and will) assume that $\X$ is
Gorenstein, except at a point $x$, and that each irreducible
component $X_i $ of $X$ passing through $x$ is a plane, denoted by
$\Pi_i$.

First we will deal with the case $n=3$.

\begin{claim}\label{c:c3}
If $x$ is a $R_3$-point, then
\[
c_x=1.
\]
\end{claim}

\begin{proof}[Proof of the claim]
From Proposition \ref{prop:ssr}-$(1)$, let us blow-up the point $x\in\X$ as in Corollary
\ref{cor:hypsect}.

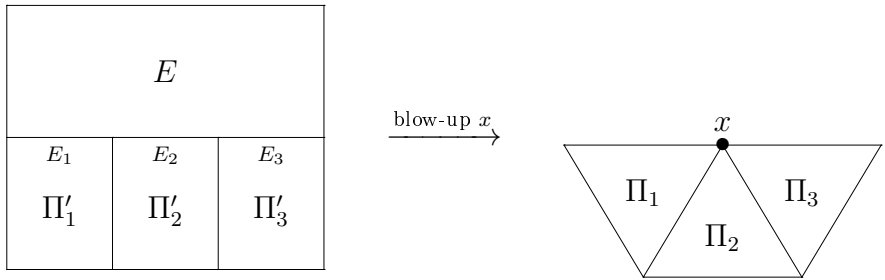
\begin{figure}[ht]
  \[
    \xymatrix@C=40pt@R=50pt{%
      *=0{} \ar@{-}[rrr] \ar@{-}[dd] & *=0{} & *=0{} & *=0{} \ar@{-}[dd] \\
      *=0{} \ar@{-}[r]_{E_1} & *=0{} \ar@{-}[r]_{E_2} \ar@{-}[d] & *=0{} \ar@{-
}[r]_{E_3} \ar@{-}[d] & *=0{} \\
      *=0{} \ar@{-}[rrr] & *=0{} & *=0{} & *=0{}
      \save "1,1"."2,4"!C*\txt{$E$} \restore
      \save "2,1"."3,2"!C-<0ex,.5ex>*\txt{$\Pi'_1$} \restore
      \save "2,2"."3,3"!C-<0ex,.5ex>*\txt{$\Pi'_2$} \restore
      \save "2,3"."3,4"!C-<0ex,.5ex>*\txt{$\Pi'_3$} \restore
    }
    \qquad
    \raisebox{-50pt}{$\xrightarrow{\text{blow-up }x}$}
    \qquad
    \raisebox{-50pt}{%
    \xymatrix@C=30pt@R=50pt{%
      *=0{} \ar@{-}[rrrr] 
        \ar@{-}[dr] & *=0{} & *=0{\bullet} \ar@{-}[dr] \ar@{-}[dl] & *=0{} &
*=0{} \ar@{-}[dl] \\
      *=0{} & *=0{} \ar@{-}[rr] & *=0{} & *=0{} & *=0{}
      \save "1,3"+<0ex,1.5ex> *\txt{$x$} \restore
      \save "1,1"."2,3"!C+<0pt,7pt> *\txt{$\Pi_1$} \restore
      \save "1,2"."2,4"!C-<0pt,10pt> *\txt{$\Pi_2$} \restore
      \save "1,3"."2,5"!C+<0pt,7pt> *\txt{$\Pi_3$} \restore
    }
    }
  \]
\caption{Blowing-up a $R_3$-point $x$.}\label{fig:blowR3}
\end{figure}

We get a new total space $\X'$. We denote by $E$ the exceptional
divisor, by $\Pi'_i$ the proper transform of $\Pi_i$ and by
$X'=\cup\Pi'_i$ the proper transform of $X$, as in Figure
\ref{fig:blowR3}. We remark that the three planes $\Pi_i$,
$i=1,2,3$, concurring at $x$, are blown-up in this process,
whereas the remaining planes stay untouched. We call $E_i$ the
exceptional divisor on the blown-up plane $\Pi_i$. Let
$\Gamma=E_1+E_2+E_3$ be the intersection curve of $E$ and $X'$. By
Corollary \ref{cor:hypsect}, $E$ is a non-degenerate surface of
degree 3 in $\Pp^4$, with $\Gamma$ as a hyperplane section.

Suppose first that $E$ is irreducible. Then $\X'$ is Gorenstein
and by adjunction:
\begin{equation}\label{eq:k2=}
K^2=(K_{X'}+\Gamma)^2+(K_E+\Gamma)^2.
\end{equation}
Since $E$ is a rational normal cubic scroll in $\Pp^4$, then:
\begin{equation}\label{eq:=1}
(K_E+\Gamma)^2=1,
\end{equation}
whereas the other term is:
\[
(K_{X'}+\Gamma)^2=\sum_i (K_{X'|\Pi'_i} + \Gamma_{\Pi'_i})^2
=\sum_{i=1}^3 (K_{X'|\Pi'_i} + E_i)^2+\sum_{j\ge4}
K_{X'|\Pi'_j}^2.
\]
Reasoning as in the proof of Theorem \ref{thm:k2Gnc}, one sees
that
\[
\sum_{j\ge4} K_{X'|\Pi'_j}^2=\sum_{j\ge4}(w_j-3)^2.
\]
On the other hand, one has
\[
(K_{X'|\Pi'_i}+E_i)^2=(w_i-3)^2-1 ,\quad  i=1,3, \qquad
(K_{X'|\Pi'_2}+E_2)^2=(w_2-3)^2.
\]
Putting all together, it follows that $c_x=1$.

Suppose now that $E$ is reducible and $X'$ is still Gorenstein. In
this case $E$ is as described in Proposition \ref{prop:hypsect}
(ii), b, and in Corollary \ref{cor:hypsect} and the proof proceeds
as above, once one remarks that \eqref{eq:=1} holds. This can be
left to the reader to verify (see Figure \ref{fig:Ered}).

\begin{figure}[ht]
  \[
  \xymatrix@C=60pt@R=60pt{%
  *=0{} & *=0{} \ar@{-}[d]_{1}^{0} \ar@{-}[dr]_{0} & *=0{} \\
  *=0{} \ar@{-}[ur]_{1} \ar@{-}[r]_{E_1}^{1} & *=0{} & *=0{}
  \save "2,2"."2,3"!C="a" \restore
  \ar@{-} "2,2";"a"-<0pt,10pt> _{E_2}^{0}  \ar@{-} "a"-<0pt,10pt>;"2,3"
_{E_3}^{0}
  }
  \]
  \caption{$E$ splits in a plane and a quadric.}
  \label{fig:Ered}
\end{figure}
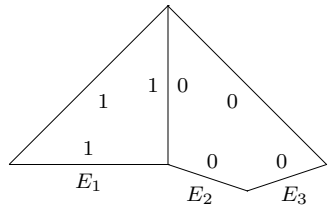

Suppose that $E$ is reducible and $X'$ is not Gorenstein. This
means that $E$ consists of a cone over a $C_{R_3}$ with vertex
$x'$, hence $x'$ is again a $R_3$-point. Therefore, as in Proposition \ref{prop:ssr}-$(1)$,
we have to repeat the process by blowing-up $x'$. After finitely many steps
this procedure stops (cf.\ e.g.\ Proposition 3.4.13 in \cite{Kollar}).
In order to conclude the proof in this case,
one has simply to remark that no contribution to $K^2$ comes from
the surfaces created in the intermediate steps.

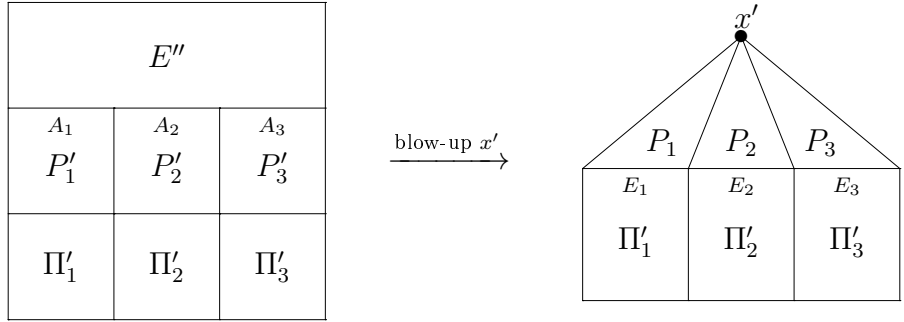
\begin{figure}[ht]
  \[
    \xymatrix@C=40pt@R=40pt{%
      *=0{} \ar@{-}[rrr] \ar@{-}[ddd] & *=0{} & *=0{} & *=0{} \ar@{-}[ddd] \\
      *=0{} \ar@{-}[r]_{A_1} & *=0{} \ar@{-}[r]_{A_2} \ar@{-}[d] & *=0{} \ar@{-
}[r]_{A_3} \ar@{-}[d] & *=0{} \\
      *=0{} \ar@{-}[r] & *=0{} \ar@{-}[r] \ar@{-}[d] & *=0{} \ar@{-}[r] \ar@{-
}[d] & *=0{} \\
      *=0{} \ar@{-}[rrr] & *=0{} & *=0{} & *=0{}
      \save "1,1"."2,4"!C*\txt{$E''$} \restore
      \save "2,1"."3,2"!C-<0ex,.5ex>*\txt{$P'_1$} \restore
      \save "2,2"."3,3"!C-<0ex,.5ex>*\txt{$P'_2$} \restore
      \save "2,3"."3,4"!C-<0ex,.5ex>*\txt{$P'_3$} \restore
      \save "3,1"."4,2"!C*\txt{$\Pi'_1$} \restore
      \save "3,2"."4,3"!C*\txt{$\Pi'_2$} \restore
      \save "3,3"."4,4"!C*\txt{$\Pi'_3$} \restore
    }
    \qquad
    \raisebox{-60pt}{$\xrightarrow{\text{blow-up }x'}$}
    \qquad
    \raisebox{-10pt}{
    \xymatrix@C=40pt@R=50pt{%
      *=0{} & *=0{} & *=0{} & *=0{} \\
      *=0{} \ar@{-}[d] \ar@{-}[r]_{E_1}
        & *=0{} \ar@{-}[r]_{E_2} \ar@{-}[d]
        & *=0{} \ar@{-}[r]_{E_3} \ar@{-}[d]
        & *=0{} \ar@{-}[d] \\
      *=0{} \ar@{-}[rrr] & *=0{} & *=0{} & *=0{}
      \save "1,1"."1,4"!C="x" \restore
      \save "x" *\txt{$\bullet$} \restore
      \save "x"+<.4ex,1.5ex> *\txt{$x'$} \restore
      \ar@{-} "x";"2,1" \ar@{-} "x";"2,2" \ar@{-} "x";"2,3" \ar@{-} "x";"2,4"
      \save "2,1"."3,2"!C-<0ex,.5ex>*\txt{$\Pi'_1$} \restore
      \save "2,2"."3,3"!C-<0ex,.5ex>*\txt{$\Pi'_2$} \restore
      \save "2,3"."3,4"!C-<0ex,.5ex>*\txt{$\Pi'_3$} \restore
      \save "1,1"."2,2"!C+< 10pt,-15pt>*\txt{$P_1$} \restore
      \save "1,2"."2,3"!C+<  0pt,-15pt>*\txt{$P_2$} \restore
      \save "1,3"."2,4"!C+<-10pt,-15pt>*\txt{$P_3$} \restore
    }
   }
  \]
\caption{blowing-up a $R_3$-point $x'$ infinitely near to the
$R_3$-point $x$}\label{fig:blowR3twice}
\end{figure}

To see this, it suffices to make this computation when only two
blow-ups are needed. This is the situation showed in Figure
\ref{fig:blowR3twice} where:
\begin{itemize}
\item  $\X''\to \X'$ is the blow-up at $x'$, \item $X'=\sum
\Pi'_i$ the proper transform of $X'$ on $\X''$, \item
$E'=P'_1+P'_2+P'_3$ is the strict transform of $E=P_1+P_2+P_3$ on
$\X''$, \item $E''$ is the exceptional divisor of the blow-up.
\end{itemize}
We remark that $P'_i$, $i=1,2,3$, is the blow-up of the plane
$P_i$. We denote by $\Lambda_i$ the pullback to $P'_i$ of a line,
and by $A_i$ the exceptional divisor of $P'_i$. Then their
contributions to the computation of $K^2$ are:
\begin{align*}
&
(K_{P'_i}+\Lambda_i+(\Lambda_i-A_i)+A_i)^2=(-\Lambda_i+A_i)^2=0,\quad
i=1,3,
\\
& (K_{P'_2}+\Lambda_2+2(\Lambda_2-A_2)+A_2)^2=0.
\end{align*}
This concludes the proof of Claim \ref{c:c3}.
\end{proof}

Consider now the case that $n=4$ and $x$ is a $R_4$-point.

\begin{claim}\label{c:c4}
If $x$ is a $R_4$-point, then
\[
2 \ge c_x \ge 1.
\]
\end{claim}

\begin{proof}[Proof of the claim]
As before, we blow-up the point $x\in\X$; let $\X'$ be the new
total space and let $E$ be the exceptional divisor. By Corollary
\ref{cor:hypsect}, $E$ is a non-degenerate surface of minimal
degree in $\Pp^5$ with $\Gamma=E_1+E_2+E_3+E_4$ as a hyperplane
section. By Proposition \ref{prop:hypsect}, $E$ is reducible and
the following cases may occur:

\begin{itemize}

\item[(i)] $E$ has global normal crossings, in which case $E$
consists of two quadrics $Q_1,Q_2$ meeting along a line (see
Figure \ref{fig:c-r4});

\begin{figure}[ht]
  \[
  \begin{xy}
  0; <50pt,0pt>: 
  (0,.7)="p";(.7,0)**@{-}="dml";(2.7,0)**@{-}="dmr";
  (3.4,.7)**@{-};(3.4,1.7)**@{-}="ur";(0,1.7)**@{-}="ul";"p"**@{-} ,
  (1.7,0);(1.7,1)**@{-};(1.7,1.7)**@{-} ,
  "ul";(.7,1)**@{-}="ml" ?+<0pt,-6pt>*{E_1};
  "ml";(1.7,1)**@{-}="m" ?+<0pt,-5pt>*{E_2};
  "m";(2.7,1)**@{-}="mr" ?+<0pt,-5pt>*{E_3};
  "mr";"ur"**@{-}        ?+<1pt,-8pt>*{E_4} ,
  "dml";"ml"**@{-} ,
  "dmr";"mr"**@{-} ,
  "dml"."m"!C+<0pt,-3pt>*=0{\dt \Pi'_2} ,
  "dmr"."m"!C+<0pt,-3pt>*=0{\dt \Pi'_3} ,
  "ul"."m"!C+(.2,0)*=0{\dt Q_1} ,
  "ur"."m"!C-(.2,0)*=0{\dt Q_2} ,
  "ul"."dml"!C+<0pt,-3pt>*=0{\dt \Pi'_1} ,
  "ur"."dmr"!C+<2pt,-3pt>*=0{\dt \Pi'_4}
  \end{xy}
  \]
  \caption{The exceptional divisor $E$ has global normal
crossings.}\label{fig:c-r4}
\end{figure}
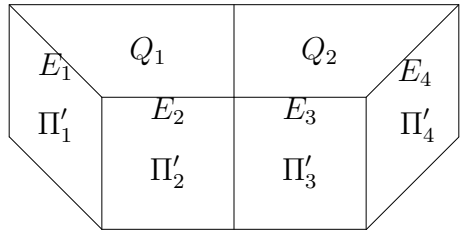

\item[(ii)] $E$ has one $R_3$-point $x'$, in which case $E$
consists of a quadric $Q$ and two planes $P_1,P_2$ (see Figure
\ref{fig:c-r4ii});

\begin{figure}[ht]
  \[
  \quad
  \begin{array}{cc}
  \begin{xy}
  0; <50pt,0pt>: 
  (0,1)="p";(.9,.5)**@{-}="dml";(1.8,0)**@{-}="dm";
  (2.7,.5)**@{-}="dmr";(3.6,1)**@{-}="pp";(3.6,2)**@{-}="ur";
  (1.8,2.5)**@{-}="uu";(0,2)**@{-}="ul";"p"**@{-} ,
  "ul";(.9,1.5)**@{-}="ml" ?+<0pt,-6pt>*{E_1}="e1" ;
  "ml";(1.8,1)**@{-}="m"   ?+<0pt,-6pt>*{E_2} ;
  "m";(2.7,1.5)**@{-}="mr" ?+<1pt,-8pt>*{E_3} ;
  "mr";"ur"**@{-}          ?+<1pt,-8pt>*{E_4}="e4" ,
  "dml";"ml"**@{-};"uu"**@{-} ?+<4pt,-3pt>*{A_1} ,
  "dm";"m"**@{-};
  "dmr";"mr"**@{-};"uu"**@{-} ?-<4pt,3pt>*{A_2} ,
  "m"."uu"!C-<0pt,5pt>*=0{\dt Q} ,
  "dml"."m"!C+<0pt,-3pt>*=0{\dt \Pi'_2} ,
  "dmr"."m"!C+<2pt,-3pt>*=0{\dt \Pi'_3} ,
  "ul"."dml"!C+<0pt,-3pt>*=0{\dt \Pi'_1} ,
  "ur"."dmr"!C+<2pt,-3pt>*=0{\dt \Pi'_4} ,
  "e1"."uu"!C-<12pt,6pt>*=0{\dt P_1} ,
  "e4"."uu"!C+<12pt,-6pt>*=0{\dt P_2} ,
  "uu"*=0{\bullet}; "uu"+<.4ex,1.5ex>*=0{\dt x'}
  \end{xy}
  \quad & \quad
  \raisebox{15pt}{%
  \begin{xy}
  0; <50pt,0pt>: 
  (0,.7)="p";(.7,0)**@{-}="dml";(2.7,0)**@{-}="dmr";
  (3.4,.7)**@{-};(3.4,1.7)**@{-}="ur";(1.7,2.2)**@{-}="uu";
  (0,1.7)**@{-}="ul";"p"**@{-} ,
  (1.7,0);(1.7,1)**@{-};"uu"**@{-} ?-<6pt,0pt>*{A_1} ,
  "ul";(.7,1)**@{-}="ml" ?+<0pt,-6pt>*{E_1} ;
  "ml";(1.7,1)**@{-}="m" ?+<0pt,-5pt>*{E_2} ;
  "m";(2.7,1)**@{-}="mr" ?+<0pt,-5pt>*{E_3}="e3" ;
  "mr";"ur"**@{-}        ?+<1pt,-8pt>*{E_4}="e4" ,
  "dml";"ml"**@{-} ,
  "dmr";"mr"**@{-} ,
  "mr";"uu"**@{-} ?+<6pt,3pt>*{A_2} ,
  "dml"."m"!C+<0pt,-3pt>*=0{\dt \Pi'_2} ,
  "dmr"."m"!C+<0pt,-3pt>*=0{\dt \Pi'_3} ,
  "ml"."uu"!C-<6pt,6pt>*=0{\dt Q} ,
  "e4"."uu"!C+<13pt,-6pt>*=0{\dt P_2} ,
  "e3"."uu"!C+<2pt,-6pt>*=0{\dt P_1} ,
  "ul"."dml"!C+<0pt,-3pt>*=0{\dt \Pi'_1} ,
  "ur"."dmr"!C+<2pt,-3pt>*=0{\dt \Pi'_4} ,
  "uu"*=0{\bullet}; "uu"+<.4ex,1.5ex>*=0{\dt x'}
  \end{xy}
  }
  \quad \\
  a) \text{ The quadric in the middle}
  &
  b) \text{ The quadric on one side}
  \end{array}
  \]
  \caption{$E$ consists of a quadric and two planes and has a $R_3$-point
$x'$.}\label{fig:c-r4ii}
\end{figure}
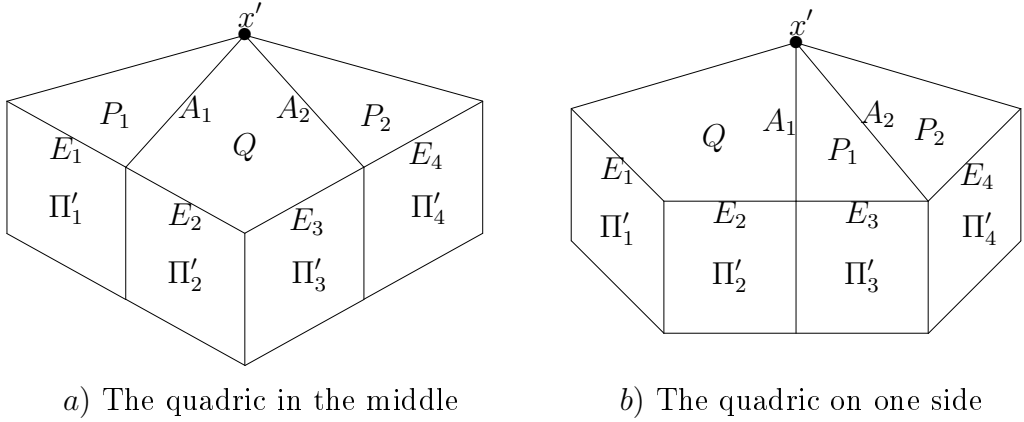

\item[(iii)] $E$ has two $R_3$-points $x',x''$, in which case $E$
consists of four planes $P_1,\ldots,P_4$, i.e.\ a planar Zappatic
surface whose associated graph is the tree $R_4$ (see Figure
\ref{fig:c-r4iii});

\begin{figure}[ht]
  \[
  \begin{xy}
  0; <50pt,0pt>: 
  (0,.7)="p";(.7,0)**@{-}="dml";(2.7,0)**@{-}="dmr";
  (3.4,.7)**@{-};(3.4,1.7)**@{-}="ur" ,
  (1.7,2.2)="uu"; (0,1.7)**@{-}="ul";"p"**@{-} ,
  (1.7,0);(1.7,1)**@{-};"uu"**@{-} ?-<6pt,3pt>*{\st A_2} ,
  "ul";(.7,1)**@{-}="ml" ?+<0pt,-6pt>*{\st E_1}="e1" ;
  "ml";(1.7,1)**@{-}="m" ?+<0pt,-5pt>*{\st E_2}="e2" ;
  "m";(2.7,1)**@{-}="mr" ?+<0pt,-5pt>*{\st E_3}="e3" ;
  "mr";"ur"**@{-}        ?+<1pt,-8pt>*{\st E_4}="e4" ,
  "dml";"ml"**@{-} ,
  "dmr";"mr"**@{-} ,
  "ml";"uu"**@{-} ?+<-6pt,3pt>*{\st A_1} ,
  "dml"."m"!C+<0pt,-3pt>*=0{\dt \Pi'_2} ,
  "dmr"."m"!C+<0pt,-3pt>*=0{\dt \Pi'_3} ,
  "e1"."uu"!C-<13pt,6pt>*=0{\dt P_1} ,
  "e2"."uu"!C+<-2pt,-15pt>*=0{\dt P_2} ,
  "e3"."uu"!C+<2pt,-15pt>*=0{\dt P_3} ,
  "ul"."dml"!C+<0pt,-3pt>*=0{\dt \Pi'_1} ,
  "ur"."dmr"!C+<2pt,-3pt>*=0{\dt \Pi'_4} ,
  "uu"*=0{\bullet}; "uu"+<.4ex,1.5ex>*=0{\dt x'} ,
  (1.7,1.8)*=0{\bullet}="uv" ; "uv"+<1.5ex,1.5ex>*=0{\dt x''} ,
  "e4"."uv"!C+<10pt,-3pt>*=0{\dt P_4} ,
  "ur";"uv"**@{-} ?!{"mr";"uu"}="ii" ,
  "uv";"mr"**@{-} ?+<-14pt,7pt>*{\st A_3} ,
  "uu";"ii"**@{-};"mr"**@{.}
  \end{xy}
  \]
  \caption{$E$ consists of four planes and has two $R_3$-points
$x',x''$.}\label{fig:c-r4iii}
\end{figure}
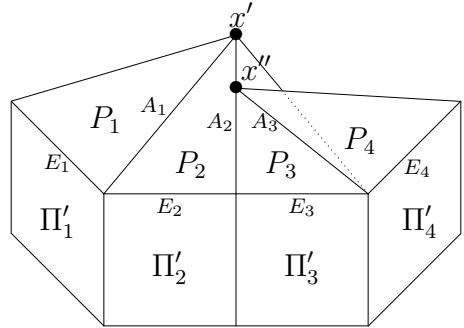

\item[(iv)] $E$ has one $R_4$-point $x'$, in which case $E$
consists of four planes, i.e.\ a planar Zappatic surface whose
associated graph is an open 4-face (cf.\ Figures \ref{fig:R3E3R4},
\ref{fig:graphR3E3R4} and \ref{fig:blowR4twice}).
\end{itemize}

In case (i), $\X'$ is Gorenstein and we can compute $K^2$ as we
did in the proof of Claim \ref{c:c3}. Formula \eqref{eq:k2=} still
holds and one has $ (K_E+\Gamma)^2=0, $ whereas:
\begin{equation}\label{eq:k2x'}
(K_{X'}+\Gamma)^2=\sum_i (K_{X'|\Pi'_i} + \Gamma_{\Pi'_i})^2
=\sum_{i=1}^4 (K_{X'|\Pi'_i} + E_i)^2+\sum_{j\ge4} K_{X'|\Pi'_j}^2
=\sum_{j\ge1} (w_j-3)^2-2,
\end{equation}
because the computations on the blown-up planes
$\Pi'_1,\ldots,\Pi'_4$ give:
\[
(K_{X'|\Pi'_i}+E_i)^2=(w_i-3)^2-1 ,\quad  i=1,4, \qquad
(K_{X'|\Pi'_i}+E_i)^2=(w_i-3)^2,\quad  i=2,3.
\]
This proves that $c_x=2$ in this case.

In case (ii), there are two possibilities corresponding to cases
(a) and (b) of
 Figure
\ref{fig:c-r4ii}. Let us first consider the former possibility. By
Claim \ref{c:c3}, in order to compute $K^2$ we have to add up
three quantities:
\begin{itemize}

\item the contribution of $(K_{X'}+\Gamma)^2$, which is computed
in \eqref{eq:k2x'};

\item the contribution to $K^2$ of $E$, as if $E$ had only global
normal crossings, i.e.:
\[
(K_{P_1}+A_1+E_1)^2+(K_{P_2}+A_2+E_4)^2+(K_{Q}+A_1+A_2+E_2+E_3)^2=2
\]

\item the contribution of the $R_3$-point $x'$, which is
$c_{x'}=1$ by Claim \ref{c:c3}.

\end{itemize}
Putting all this together, it follows that $c_x=1$ in this case.
Consider now the latter possibility, i.e.\ suppose that the
quadric meets only one plane. We can compute the three
contributions to $K^2$ as above: the contribution of
$(K_{X'}+\Gamma)^2$ and of the $R_3$-point $x'$ do not change,
whereas the contribution to $K^2$ of $E$, as if $E$ had only
global normal crossings, is:
\[
(K_{Q}+A_1+E_1+E_2)^2+(K_{P_1}+A_1+A_2+E_3)^2+(K_{P_4}+A_3+E_4)^2=1,
\]
therefore we find that $c_x=2$, which concludes the proof for case
(ii).

In case (iii), we use the same strategy as in case (ii), namely we
add up $(K_{X'}+\Gamma)^2$, the contribution to $K^2$ of $E$, as
if $E$ had only global normal crossings, which turns out to be 2,
and then subtract 2, because of the contribution of the two
$R_3$-points $x',x''$. Summing up, one finds $c_x=2$ in this case.

In case (iv), we have to repeat the process by blowing-up $x'$,
see Figure \ref{fig:blowR4twice}. After finitely many steps
(cf.\ e.g.\ Proposition 3.4.13 in \cite{Kollar}), this
procedure stops in the sense that the exceptional divisor will be
as in case (i), (ii) or (iii).

\begin{figure}[ht]
  \[
\raisebox{120pt}{
    \xymatrix@C=40pt@R=40pt{%
      *=0{} \ar@{-}[rrrr] \ar@{-}[ddd] & *=0{} & *=0{} & *=0{} & *=0{} \ar@{-
}[ddd] \\
      *=0{} \ar@{-}[r] & *=0{} \ar@{-}[r] \ar@{-}[d]
        & *=0{} \ar@{-}[r] \ar@{-}[d] & *=0{} \ar@{-}[r] \ar@{-}[d] & *=0{} \\
      *=0{} \ar@{-}[r] & *=0{} \ar@{-}[r] \ar@{-}[d] & *=0{} \ar@{-}[r] \ar@{-
}[d]
        & *=0{} \ar@{-}[r] \ar@{-}[d] & *=0{} \\
      *=0{} \ar@{-}[rrrr] & *=0{} & *=0{} & *=0{} & *=0{}
      \save "1,1"."2,5"!C*\txt{$E''$} \restore
      \save "2,1"."3,2"!C*\txt{$P'_1$} \restore
      \save "2,2"."3,3"!C*\txt{$P'_2$} \restore
      \save "2,3"."3,4"!C*\txt{$P'_3$} \restore
      \save "2,4"."3,5"!C*\txt{$P'_4$} \restore
      \save "3,1"."4,2"!C*\txt{$\Pi'_1$} \restore
      \save "3,2"."4,3"!C*\txt{$\Pi'_2$} \restore
      \save "3,3"."4,4"!C*\txt{$\Pi'_3$} \restore
      \save "3,4"."4,5"!C*\txt{$\Pi'_4$} \restore
    }
}
    \qquad
    \raisebox{60pt}{$
      \xrightarrow{\text{blow-up }x'}
    $}
    \qquad
    \begin{xy}
    0; <50pt,0pt>: 
    (0,.7)="p";(.7,0)**@{-}="dml";(2.7,0)**@{-}="dmr";
    (3.4,.7)**@{-};(3.4,1.7)**@{-}="ur" ,
    (1.7,2.2)="uu"; (0,1.7)**@{-}="ul";"p"**@{-} ,
    (1.7,0);(1.7,1)**@{-};"uu"**@{-}  ,
    "ul";(.7,1)**@{-}="ml" ?+<0pt,-6pt>*{\st E_1} ;
    "ml";(1.7,1)**@{-}="m" ?+<0pt,-5pt>*{\st E_2} ;
    "m";(2.7,1)**@{-}="mr" ?+<0pt,-5pt>*{\st E_3} ;
    "mr";"ur"**@{-} ?+<1pt,-8pt>*{\st E_4} ,
    "dml";"ml"**@{-} ,
    "dmr";"mr"**@{-} ,
    "mr";"uu"**@{-} ,
    "ml";"uu"**@{-}  ,
    "ur";"uu"**@{-}  ,
    "ul"."dml"!C+<0pt,-3pt>*=0{\dt \Pi'_1} ,
    "dml"."m"!C+<0pt,-3pt>*=0{\dt \Pi'_2} ,
    "dmr"."m"!C+<0pt,-3pt>*=0{\dt \Pi'_3} ,
    "ur"."dmr"!C+<2pt,-3pt>*=0{\dt \Pi'_4} ,
    "ul";"ml"**@{} ?*{}; "uu"**@{} ?(.3)*=0{\dt P_1} ,
    "ml";"m"**@{} ?*{}; "uu"**@{} ?(.3)*=0{\dt P_2} ,
    "m";"mr"**@{} ?*{}; "uu"**@{} ?(.3)*=0{\dt P_3} ,
    "mr";"ur"**@{} ?*{}; "uu"**@{} ?(.3)*=0{\dt P_4} ,
    "uu"*=0{\bullet}; "uu"+<.4ex,1.5ex>*=0{\dt x'} ,
    \end{xy}
  \]
\caption{Blowing-up a $R_4$-point $x'$ infinitely near to
$x$.}\label{fig:blowR4twice}
\end{figure}
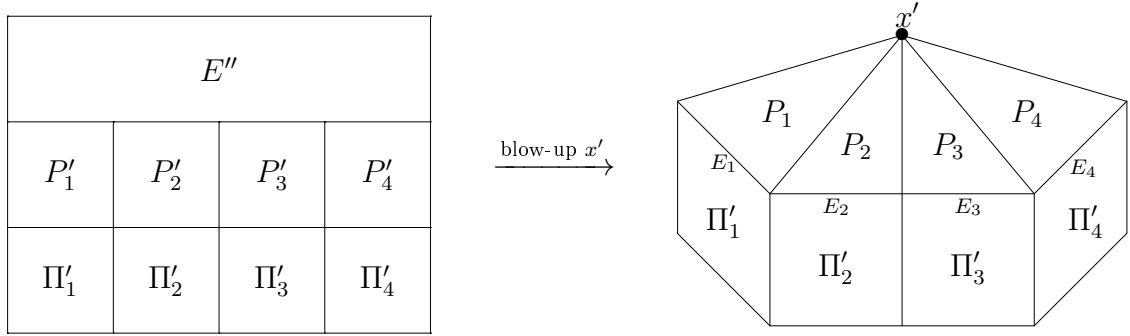

\noindent In order to conclude the proof of Claim \ref{c:c4}, one
has to remark that no contribution to $K^2$ comes from the
surfaces created in the intermediate steps (the blown-up planes
$P'_i$ in Figure \ref{fig:blowR4twice}). This can be done exactly
in the same way as we did in the proof of Claim \ref{c:c3}.
\end{proof}

\begin{remark}\label{rem:R4expli}
The proof of Claim \ref{c:c4} is purely combinatorial. However
there is a nice geometric motivation for the two cases $c_x=2$ and
$c_x=1$, when $x$ is a $R_4$-point, which resides in the fact that
the local deformation space of a $R_4$-point is reducible. This
corresponds to the fact that the cone over $C_{R_4}$ can be
smoothed in both a Veronese surface and a rational normal quartic
scroll, which have $K^2=9$ and $K^2=8$, respectively.
\end{remark}

Consider now the case that $x$ is a $R_n$-point.

\begin{claim}\label{c:cn+}
If $x$ is a $R_n$-point, then
\begin{equation}\label{eq:cn+}
n-2 \ge c_x \ge 1.
\end{equation}
\end{claim}

\begin{proof}[Proof of the claim]
The claim for $n=3,4$ has already been proved, so we assume $n \ge
5$ and proceed by induction on $n$. As usual, we blow-up the point
$x\in\X$.

By Corollary \ref{cor:hypsect}, the exceptional divisor
$E$ is a non-degenerate surface of minimal degree in $\Pp^{n+1}$
with $\Gamma=E_1+\ldots+E_n$ as a hyperplane section. By
Proposition \ref{prop:hypsect}, $E$ is reducible and the following
cases may occur:

\begin{itemize}

\item[(i)] $E$ consists of $\nu\ge3$ irreducible components
$P_1,\ldots,P_\nu$, which are either planes or smooth quadrics,
and $E$ has $h$ Zappatic singular points $x_1,\ldots,x_h$ of type
$R_{m_1},\ldots,R_{m_h}$ such that $m_i < n$, $i=1,\ldots,h$;

\item[(ii)] $E$ has one $R_n$-point $x'$, in which case $E$
consists of $n$ planes, i.e.\ a planar Zappatic surface whose
associated graph is an open $n$-face.

\end{itemize}

In case (ii), one has to repeat the process by blowing-up $x'$.
After finitely many steps (cf.\ e.g.\ Proposition 3.4.13
in \cite{Kollar}), the exceptional divisor will
necessarily be as in case (i). We remark that no contribution to
$K^2$ comes from the surfaces created in the intermediate steps,
as one can prove exactly in the same way as we did in the proof of
Claim \ref{c:c3}.

Thus, it suffices to prove the statement for the case (i). Notice
that $\X'$ is not Gorenstein, nonetheless we can compute $K^2$
since we know (the upper and lower bounds of) the contribution of
$x_i$ by induction. We can indeed proceed as in case (ii) of the
proof of Claim \ref{c:c4}, namely, we have to add up three
quantities:

\begin{itemize}

\item the contribution of $(K_{X'}+\Gamma)^2$;

\item the contribution to $K^2$ of $E$, as if $E$ had only global
normal crossings;

\item the contributions of the points $x_i$ which is known  by
induction.
\end{itemize}

Let us compute these contributions. As for the first one, one has:
\[
(K_{X'}+\Gamma)^2=\sum_{i=1}^n (K_{X'|\Pi'_i} + E_i)^2+\sum_{j\ge
n} K_{X'|\Pi'_j}^2 =\sum_{j\ge1} (w_j-3)^2-2,
\]
since the computations on the blown-up planes
$\Pi'_1,\ldots,\Pi'_n$ give:
\begin{align*}
&(K_{X'|\Pi'_i}+E_i)^2=(w_i-3)^2-1,\quad  i=1,n, \\
&(K_{X'|\Pi'_i}+E_i)^2=(w_i-3)^2,\quad  2\le i \le n-1.
\end{align*}

In order to compute the second contribution, one has to introduce
some notation, precisely we let:

\begin{itemize}

\item $P_1,\ldots,P_\nu$ be the irreducible components of $E$,
which are either planes or smooth quadrics, ordered in such a way
that the intersections in codimension one are as follows: $P_i$
meets $P_{i+1}$, $i=1,\ldots,\nu-1$, along a line;

\item $A_i$ be the line which is the intersection of $P_i$ and
$P_{i+1}$;

\item $\eps_i=\deg(P_i)-1$, which is 0 if $P_i$ is a plane and 1
if $P_i$ is a quadric;

\item $j(i)=i+\sum_{k=1}^{i-1} \eps_j$. With this notation, if
$P_i$ is a plane, it meets the blown-up plane $\Pi'_{j(i)}$ along
$E_{j(i)}$, whereas if $P_i$ is a quadric, it meets the blown-up
planes $\Pi'_{j(i)}$ and $\Pi'_{j(i)+1}$ along $E_{j(i)}$  and
$E_{j(i)+1}$, respectively.

\end{itemize}

Then the contribution to $K^2$ of $E$, as if $E$ had only global
normal crossings, is:
\begin{align*}
&(K_{P_1}+A_1+E_1+\eps_1 E_2)^2 + (K_{P_\nu}+A_{\nu-1}+\eps_\nu
E_{n-1}+E_n)^2+
\\
&+\sum_{i=2}^{\nu-1}(K_{P_i}+A_{i-1}+A_i+E_{j(i)}+\eps_i
E_{j(i)+1})^2=2-\eps_1- \eps_\nu.
\end{align*}

Finally, by induction, the contribution $\sum_{i=1}^h c_{x_i}$ of
the points $x_i$ is such that:
\[
\nu -2 = \sum_{i=1}^h (m_i-2) \ge \sum_{i=1}^h c_{x_i} \ge
\sum_{i=1}^h 1=h,
\]
where the first equality is just \eqref{eq:nu-2}.

Putting all this together, it follows that:
\[
c_x = \eps_1+\eps_\nu+\sum_{i=1}^h c_{x_i},
\]
hence an upper bound for $c_x$ is
\[
c_x \le \eps_1+\eps_\nu+\nu-2 \le n-2,
\]
because $n=\nu+\sum_{i=1}^\nu \eps_i$, whereas a lower bound is
\begin{equation}\label{eq:cx>=1}
c_x \ge \eps_1+\eps_\nu+ h \ge h \ge 1,
\end{equation}
which concludes the proof of Claim \ref{c:cn+}.
\end{proof}

\begin{remark}
If $c_x=1$, then in \eqref{eq:cx>=1} all inequalities must be
equalities, thus $h=1$ and $\eps_1=\eps_\nu=0$. This means that
there is only one point $x_1$ infinitely near to $x$, of type
$R_{\nu}$, and that the external irreducible components of $E$,
i.e.\ $P_1$ and $P_\nu$, are planes. There is no combinatorial
obstruction to this situation.

For example, let $x$ be a $R_n$-point such that the exceptional
divisor $E$ consists of $\nu=n-1$ irreducible components, namely
$n-2$ planes and a quadric adjacent to two planes, forming
a $R_{n-1}$-point $x'$. By the proof of Claim \ref{c:c4} (case (ii),
former possibility), it follows that $c_x=c_{x'}$. Since, as we
saw, the contribution of a $R_4$-point can be 1, by induction we
may have that also a $R_n$-point contributes by 1.

From the proof of Claim \ref{c:cn+}, it follows that the upper
bound $c_x=n-2$ is attained when for example the exceptional
divisor $E$ consists of $n$ planes forming $n-2$ points of type
$R_3$.

More generally, one can see that there is no combinatorial
obstruction for $c_x$ to attain any possible value between the
upper and lower bounds in \eqref{eq:cn+}.
\end{remark}

Finally, consider the case that $x$ is of type $S_n$.

\begin{claim}\label{c:sn+}
If $x$ is a $S_n$-point, then
\begin{equation}\label{eq:cx>^2}
(n-2)^2 \ge c_x \ge \binom{n-1}{2}.
\end{equation}
\end{claim}

\begin{proof}
We remark that we do not need to take care
of 1-dimensional singularities of the total space of the degeneration,
as we have already noted in Claim \ref{c:cn+}.

Notice that $S_3=R_3$ and, for $n=3$, Formula \eqref{eq:cx>^2}
trivially follows from Claim \ref{c:c3}. So we assume $n\ge4$.
Blow-up $x$, as usual; let $\X'$ be the new total space and $E$
the exceptional divisor.
By Proposition \ref{prop:hypsect}, three
cases may occur: either

\begin{itemize}

\item[(i)] $E$ has global normal crossings, i.e.\ $E$ is the union
of a smooth rational normal scroll $X_1=S(1,d-1)$ of degree $d$,
$2\le d \le n$, and of $n-d$ disjoint planes $P_1,\ldots,P_{n-d}$,
each meeting $X_1$ along different lines of the same ruling; or

\item[(ii)] $E$ is a union of $n$ planes $P_1,\ldots,P_n$ with $h$
Zappatic singular points $x_1,\ldots,x_h$ of type
$S_{m_1},\ldots,S_{m_h}$ such that $3\le m_i < n$, $i=1,\ldots,h$,
and \eqref{eq:n-1su2} holds; or

\item[(iii)] $E$ is a union of $n$ planes with one $S_n$-point
$x'$.

\end{itemize}

In case (iii), one has to repeat the process by blowing-up $x'$.
After finitely many steps (cf.\ e.g.\ Proposition 3.4.13 in \cite{Kollar}),
the exceptional divisor will
necessarily be as in cases either (i) or (ii). We remark that no
contribution to $K^2$ comes from the surfaces created in the
intermediate steps. Indeed, by using the same notation of the
$R_n$-case in Claim \ref{c:c3}, if $x$ is a $S_n$ point and if
$\Pi_1$ is the plane corresponding to the vertex of valence $n-1$
in the associated graph, we have (cf.\ Figure \ref{fig:Sntwice}):
\begin{align*}
& (K_{P_1'} + \Lambda_1 + A_1 + (n-1) (\Lambda_1 - A_1))^2 =
(n-3)^2 - (n-3)^2 =
0, \\
& (K_{P_i'} + \Lambda_i + A_i + (\Lambda_i - A_i))^2 = 1-1 = 0,
\;\;  2 \leq i \leq n.
\end{align*}

\begin{figure}[ht]
  \[
  \raisebox{12.5pt}{$
  \begin{xy}
    0; <25pt,0pt>: , 
    (1,0.5)="p2dr";(6,0.5)="p5dl";"p2dr"."p5dl"!C="m";
    "m"+(0,3)="u";
    "p2dr";(1,2)**@{-}="p2ur";(0,3)**@{-}="p2ul";
    (0,1.5)**@{-}="p2dl";"p2dr"**@{-} ,
    "p2dr"."p2ul"!C*{\dt P'_2} ,
    (2.5,0.5)="p3dr";(2.5,2)**@{-}="p3ur";(1.5,3)**@{-}="p3ul";
    (1.5,1.5)**@{-}="p3dl";"p3dr"**@{-} ,
    "p3dr"."p3ul"!C*{\dt P'_3} ,
    (4.5,0.5)="p4dl";(4.5,2)**@{-}="p4ul";(5.5,3)**@{-}="p4ur";
    (5.5,1.5)**@{-}="p4dr";"p4dl"**@{-} ,
    "p4dl"."p4ur"!C*{\dt P'_4} ,
    "p5dl";(6,2)**@{-}="p5ul";(7,3)**@{-}="p5ur";
    (7,1.5)**@{-}="p5dr";"p5dl"**@{-} ,
    "p5dl"."p5ur"!C*{\dt P'_5} ,
    "p3dr"."p4ul"!C*{\dt P'_1}="p1" ,
    "p3ur";"p4ul"**@{-};
    "p3ur"."p4ul"!C;"u"**@{}?*{E''};
    "p2ul";"p5ur"**\crv{"u"};
    "p3ul";"p4ur"**\crv{"u"};
    {"p2ur";"p3ur":"p3ul";"p3dl",x}="i3" ,
    "p2ur";"i3"**@{-};"p3ur"**@{.} ,
    {"p5ul";"p4ul":"p4ur";"p4dr",x}="i4" ,
    "p5ul";"i4"**@{-};"p4ul"**@{.} ,
    {"p2dr";"p3dr":"p3ul";"p3dl",x}="i5" ,
    {"p5dl";"p4dl":"p4ur";"p4dr",x}="i6" ,
    "p2dr";"i5"**@{-};"p3dr"**@{.};"p4dl"**@{-};"i6"**@{.};"p5dl"**@{-} ,
    "p2dr";"p2dr"-"p2dl"**@{-}="p2ddr";"p2dl"-"p2dl"**@{-};"p2dl"**@{-};
    "p2ddr"."p2dl"!C*{\dt \Pi'_2} ,
    "p3dr";"p3dr"-"p2dl"**@{-}="p3ddr";"p3dl"-"p2dl"**@{-};"p3dl"**@{-};
    "p3ddr"."p3dl"!C*{\dt \Pi'_3} ,
    "p4dr";"p4dr"-"p2dl"**@{-}="p4ddr";"p4dl"-"p2dl"**@{-};"p4dl"**@{-};
    "p4ddr"."p4dl"!C*{\dt \Pi'_4} ,
    "p5dr";"p5dr"-"p2dl"**@{-}="p5ddr";"p5dl"-"p2dl"**@{-};"p5dl"**@{-};
    "p5ddr"."p5dl"!C*{\dt \Pi'_5} ,
    "p3ddr"."p4dl"!C*{\dt \Pi'_1} ,
    "p2dr"-"p2dl";"p5dl"-"p2dl"**@{-} ,
  \end{xy}
    $}
    \quad
    \raisebox{50pt}{$
      \xrightarrow{\text{blow-up }x'}
    $}
    \quad
  \begin{xy}
    0; <25pt,0pt>: 
    (1,0.5)="p2dr";(6,0.5)="p5dl"**@{-};"p2dr"."p5dl"!C="m";
    "m"+(0,3)="u"*{\bullet};"u"+<2pt,9pt>*{\dt x'};
    "p2dr";(1,2)**@{-}="p2ur";(0,3)**@{-}="p2ul";
    (0,1.5)**@{-}="p2dl";"p2dr"**@{-} ,
    "p2dr"."p2ul"!C*{\dt \Pi'_2} ,
    (2.5,0.5)="p3dr";(2.5,2)**@{-}="p3ur";(1.5,3)**@{-}="p3ul";
    (1.5,1.5)**@{-}="p3dl";"p3dr"**@{-} ,
    "p3dr"."p3ul"!C*{\dt \Pi'_3} ,
    (4.5,0.5)="p4dl";(4.5,2)**@{-}="p4ul";(5.5,3)**@{-}="p4ur";(
    5.5,1.5)**@{-}="p4dr";"p4dl"**@{-} ,
    "p4dl"."p4ur"!C*{\dt \Pi'_4} ,
    "p5dl";(6,2)**@{-}="p5ul";(7,3)**@{-}="p5ur";
    (7,1.5)**@{-}="p5dr";"p5dl"**@{-} ,
    "p5dl"."p5ur"!C*{\dt \Pi'_5} ,
    "p3dr"."p4ul"!C*{\dt \Pi'_1}="p1" ,
    "p3ul";"u"**@{-};"p3ur"**@{-};"p4ul"**@{-};"u"**@{-};"p4ur"**@{-};
    "p2ul";"u"**@{-};"p5ur"**@{-} ,
    {"p2ur";"u":"p3ul";"p3dl",x}="i1" ,
    "p2ur";"i1"**@{-};"u"**@{.} ,
    {"p5ul";"u":"p4ur";"p4dr",x}="i2" ,
    "p5ul";"i2"**@{-};"u"**@{.} ,
    {"p2ur";"p3ur":"p3ul";"p3dl",x}="i3" ,
    "p2ur";"i3"**@{-};"p3ur"**@{.} ,
    {"p5ul";"p4ul":"p4ur";"p4dr",x}="i4" ,
    "p5ul";"i4"**@{-};"p4ul"**@{.} ,
    "p2ul";"p2ur"**@{}?;"p3ul"**@{}?*{\dt P_2}="p2" ,
    "p5ul";"p5ur"**@{}?;"p4ur"**@{}?*{\dt P_5}="p5" ,
    {"p2";"p5":"p3dr";"p3ur",x}+<-5pt,0pt> *{\dt P_3}="p3" ,
    {"p2";"p5":"p4dl";"p4ul",x}+<5pt,0pt> *{\dt P_4}="p4" ,
    {"m";"u":"p3";"p4",x} *{\dt P_1} ,
  \end{xy}
  \]
  \caption{Blowing-up a $S_5$-point $x'$ infinitely near to a $S_5$-point
$x$}\label{fig:Sntwice}
\end{figure}
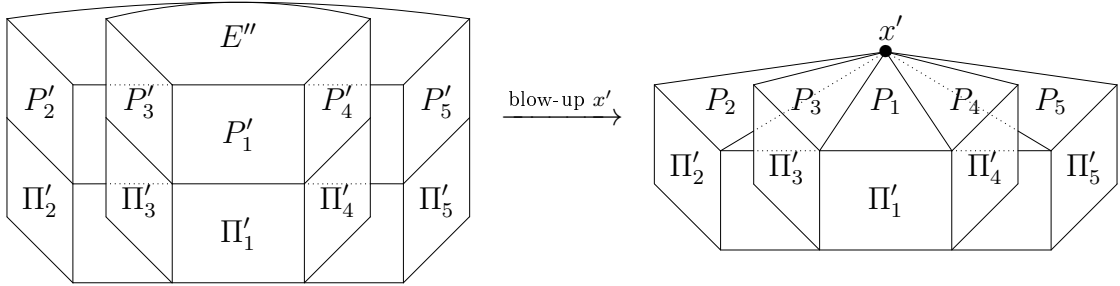

Thus, it suffices to prove the statement for the first two cases
(i) and (ii).

Consider the case (i), namely $E$ has global normal crossings.
Then $\X'$ is Gorenstein and we may compute $K^2$ as in
\eqref{eq:k2=}. The contribution of the blown-up planes
$\Pi'_1,\ldots,\Pi'_n$ (choosing again the indexes in such a way
that $\Pi'_1$ meets $\Pi'_2,\ldots,\Pi'_n$ in a line) is:
\begin{equation}\label{eq:KX'sn}
\begin{aligned}
 (K_{X'|\Pi'_i}+E_i)^2 &= (w_i-3)^2-1,\qquad  i=2,\ldots,n, \\
 (K_{X'|\Pi'_1}+E_1)^2 &= (w_1-3)^2-(n-3)^2,
\end{aligned}
\end{equation}
whereas the contribution of $E$ turns out to be:
\begin{equation}\label{eq:k2=4-n}
(K_{E}+\Gamma)^2=4-n.
\end{equation}
Indeed, one finds that:
\begin{align*}
& \left( (K_{E}+\Gamma)_{|X_1} \right)^2=(-A+(n-d-1)F)^2=d+4-2n, \\
& \left( (K_{E}+\Gamma)_{|P_i} \right)^2=1, \quad i=1,\ldots,n-d,
\end{align*}
where $A$ is the linear directrix of $X_1$ and $F$ is its fibre,
therefore \eqref{eq:k2=4-n} holds. Summing up, it follows that
\begin{equation}\label{eq:cx=^2}
c_x=n-4+(n-1)+(n-3)^2=(n-2)^2,
\end{equation}
which proves \eqref{eq:cx>^2} in this case (i).

In case (ii), $E$ is not Gorenstein, nonetheless we can compute
$K^2$ since we know (the upper and lower bounds of) the
contribution of $x_i$ by induction. We can indeed proceed as in
case (ii) of the proof of Claim \ref{c:c4}, namely, we have to add
up three quantities:

\begin{itemize}

\item the contribution of $(K_{X'}+\Gamma)^2$, which has been
computed in \eqref{eq:KX'sn};

\item the contribution to $K^2$ of $E$, as if $E$ had only global
normal crossings, which is:
\[
\Bigl( K_{P_1}+E_1+\sum_{i=2}^n A_i \Bigr)^2+ \sum_{i=2}^n
(K_{P_i}+E_i+A_i)^2 = (n-3)^2+ n-1,
\]
where $\Pi'_1$ is the blown-up plane meeting all the other
blown-up planes in a line, $E_i$ is the exceptional curve on
$\Pi'_i$ and $A_i$ is the double line intersection of $P_1$ with
$P_i$;

\item the contribution $\sum_{i=1}^h c_{x_i}$ of the points $x_i$,
which by induction, is such that:
\begin{equation}\label{eq:sn}
\sum_{i=1}^h (m_i-2)^2 \ge \sum_{i=1}^h c_{x_i} \ge \sum_{i=1}^h
\binom{m_i-1}{2}=\binom{n-1}{2},
\end{equation}
where the last equality is just \eqref{eq:n-1su2}.

\end{itemize}
Putting all together, one sees that
\[
c_x = \sum_{i=1}^h c_{x_i},
\]
hence \eqref{eq:sn} gives the claimed lower bound, as for the
upper bound:
\begin{align*}
c_x & \le \sum_{i=1}^h (m_i-2)^2
= \sum_{i=1}^h (m_i-1)(m_i-2)-\sum_{i=1}^h (m_i-2) \stackrel{(*)}{=} \\
& \stackrel{(*)}{=} (n-1)(n-2)-\sum_{i=1}^h (m_i-2)\le
(n-1)(n-2)-(n-2)=(n-2)^2,
\end{align*}
where the equality $(*)$ follows from \eqref{eq:n-1su2}. This
completes the proof of Claim \ref{c:sn+}.
\end{proof}

The above Claims \ref{c:cn+} and \ref{c:sn+} prove Proposition
\ref{prop:cx} and, so, Theorem \ref{thm:k2Gmain}.
\end{proof}

\begin{remark}
Notice that the upper bound $c_x=(n-2)^2$ is attained when for
example the exceptional divisor $E$ has global normal crossings
(cf.\ case (i) in Claim \ref{c:sn+}). The lower bound
$c_x=\binom{n-1}{2}$ can be attained if the exceptional divisor
$E$ consists of $n$ planes forming $\binom{n-1}{2}$ points of type
$S_3=R_3$.

Contrary to what happens for the $R_n$-points, not all the values
between the upper and the lower bound are realized by $c_x$,
for a $S_n$-point $x$. Indeed they are not even combinatorially
possible.
For example, consider the case of a $S_6$-point $x$:
the bounds in \eqref{eq:cx>^2} say that $10\leq c_x \leq 16$.
If the exceptional divisor $E$ has global normal crossings,
then $c_x=16$.
Otherwise $E$ is a union of planes and has the following Zappatic singularities:
$q_5 \leq 1$ points of type $S_5$,
$q_4 \leq 3-2q_5$ points of type $S_5$
and $q_3=10-3q_4-6q_5$ points of type $S_3=R_3$.
It follows that
\[
c_x\leq q_3+4q_4+9q_5=10+q_4+3q_5\leq 13+q_5 \leq 14.
\]
Therefore the case $c_x=15$ cannot occur if $x$ is a $S_6$-point.
\end{remark}

\section{The genus of the fibres of degenerations of surfaces to Zappatic ones}\label{S:pgzappdeg}

In this section we want to investigate on the behaviour of the geometric
genus of the smooth fibres of a degeneration of surfaces to a good Zappatic one,
in terms of the $\omega$-genus of the central fibre.

By recalling Definition \ref{def:degen}, the geometric genus
of the general fibre of a semistable degeneration of surfaces
can be computed via the \emph{Clemens-Schmid exact sequence},
cf.\ \cite{Morr}. Clemens-Schmid result implies the following:

\begin{theorem}\label{prop:CS}
Let $X=\bigcup_{i=1}^v X_i$ be the central fibre
of a semistable degeneration of surfaces $\X \to \D$.
Let $G_X$ be the graph associated to $X$ and $\Phi_X$ be the map
introduced in Definition \ref{rem:indices}.
Then, for $t\neq 0$, one has:
\begin{equation}\label{eq:CS}
p_g(\X_t) = h^ 2(G_X, \CC) + \sum_{i=1}^v p_g(X_i) + \dim({\rm coker}
(\Phi_X)). \end{equation}
\end{theorem}

Then Theorem \ref{prop:CS} and our Theorem \ref{thm:4.pgbis}
imply the following:

\begin{corollary}\label{cor:CS}
Let $\X\to\D$ be a semistable degeneration of surfaces,
so that its central fibre $X=\X_0$
is a good Zappatic surface with only $E_3$-points as Zappatic singularities.
Then, for any $t\ne0$, one has:
\[
p_g(\X_t)=p_\omega(X).
\]
\end{corollary}

\begin{remark}
Let $\X \to \D$ be a degeneration of surfaces
with central fibre $X$.
Consider the dualizing sheaf $\omega_{\X}$ of $\X$.
By general properties of dualizing sheaves,
one knows that $\omega_{\X}$ is torsion-free as an $\Oc_{\X}$-module.
Since one has the injection $\Oc_{\D} \hookrightarrow \Oc_{\X}$,
then $\omega_{\X}$ is torsion-free over $\Delta$.
Since $\Delta$ is the spectrum of a DVR, then $\omega_\X$ is free
and therefore flat over $\Delta$.
By semi-continuity, this implies that, for $t\ne0$, $p_g(\X_t)\leq p_\omega(X)$.
The above corollary shows that equality holds for semistable degenerations of surfaces.
\end{remark}

Consider, from now on, a degeneration $\pi:\X\to\D$ of surfaces
with good Zappatic central fibre $X=\X_0$.
Our main purpose in this section is to prove Proposition \ref{prop:ssr2},
where we show that the $\omega$-genus of the central fibre
of a semistable reduction $\tilde\pi:\tilde\X\to\D$ of $\pi$
equals the $\omega$-genus of $X$. As a consequence we will have that the
$\omega$--genus of the fibres of $\pi:\X\to\D$ is constant
(see Theorem \ref{thm:pg} below), exactly as it happens in the normal crossings case, as
we saw in Corollary \ref{cor:CS}.

In order to prove Proposition
\ref{prop:ssr2}, we make use of Proposition \ref{prop:ssr}. Indeed,
let $X=\bigcup_{i=1}^v X_i$ be the central fibre of the original degeneration
and let $\bar X_{\red}=\bigcup_{i=1}^w \bar X_i$ be
the support of the central fibre $\bar X$
of its normal crossing reduction obtained as in Proposition \ref{prop:ssr}, where $w \geq v$.
Next we describe the relation between the graph $G$ associated to $X$
and the one $\bar G$ associated to $\bar X_\red$.
By the proof of Proposition \ref{prop:ssr}, one has that $G$
is a subgraph of $\bar G$ and we may assume that $\bar X_i$ is birational to $X_i$, $i=1,\ldots,v$.

\begin{proposition}\label{XXX} (cf. \cite[Proposition 4.10]{CCFMpg}) 
In the above situation, one has:
\begin{enumerate}
\item[(i)] $p_g(\bar X_i)=0$, $i=v+1,\ldots,w$;
\item[(ii)] $\dim({\rm coker} (\Phi_{\bar X_\red}))=\dim({\rm coker}
(\Phi_{X}))$; \item[(iii)] the graphs $G$ and $\bar G$ have the same Betti numbers.
\end{enumerate}
\end{proposition}

\begin{proof}
Following the discussion of Gorenstein reduction algorithm \ref{gr} and of each Step
of the normal crossing reduction algorithm \ref{ssr},
one sees that each new component $\bar X_i$, $i=v+1,\ldots,w$,
of the central fibre is an exceptional divisor of a blow-up,
which is either a rational or a ruled surface.
This proves (i).

For $i=1,\ldots,v$, the birational morphism $\bar\sigma:\bar\X\to\X$
determines a birational morphism $\bar X_i\to X_i$ which is
the composition of blow-ups at smooth points of $X_i$.
In order to prove (ii), we notice that in algorithms \ref{gr} and  \ref{ssr},
we have added rational double curves (which do not contribute to
the cokernel), new rational components (which also do not contribute to
the cokernel), and irrational ruled surfaces, which are only created by
blowing-up irrational double curves. Focusing on single such irrational double curve, one sees that
it is replaced by a certain number $h$ of irrational ruled surfaces,
and by $h+1$ new double curves.  The map on the $H^1$ level is an isomorphism
between the new surfaces and the new curves.  Hence there is no change in the
dimension of the cokernel. This concludes the proof of (ii).

In order to prove (iii),
let us see what happens at algorithm \ref{gr} and at each step of algorithm \ref{ssr}.

In algorithm \ref{gr}, one blows-up $R_n$- and $S_n$-points of $X=\X_0$.
An example will illustrate the key features of the analysis.
Let $p$ be a $R_4$-point of $X$. After blowing-up $\X$ at $p$,
there are five different possible configurations of the exceptional divisor $E$
(cf.\ the proof of Claim \ref{c:c4}):
\begin{itemize}
\item[(i)] $E$ is the union of two quadrics with normal crossings;

\item[(ii)] $E$ is the union of a quadric and two planes
having a $R_3$-point $p'$, and the quadric is in the middle;

\item[(iii)] $E$ is the union of a quadric and two planes
having a $R_3$-point $p'$, and one of the planes is in the middle;

\item[(iv)] $E$ is the union of four planes having two $R_3$-points $p'$, $p''$;

\item[(v)] $E$ is the union of four planes having a $R_4$-point $p'$.
\end{itemize}
The corresponding associated graphs are illustrated in Figure \ref{fig:r4g},
where the proper transforms of the four components of $X$ concurring at $p$
are the left-hand-side vertices in each graph.
As the pictures show, $G$ is a deformation retract of the new associated graph
(considered as CW-complexes).

\begin{figure}[ht]
\[
\begin{array}{ccccc}
\quad
\includegraphics{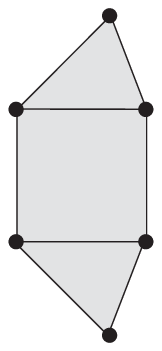}
\quad
&
\quad
\psfrag{x}{$p'$}
\includegraphics{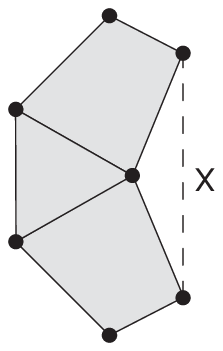}
\quad
&
\quad
\psfrag{X}{$\,p'$}
\includegraphics{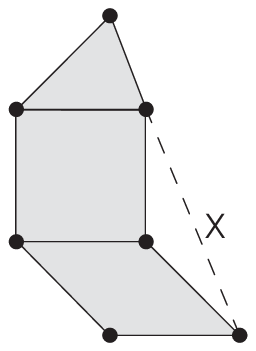}
\quad
&
\quad
\psfrag{x}{$p'$}
\psfrag{Y}{$\,\,p''$}
\includegraphics{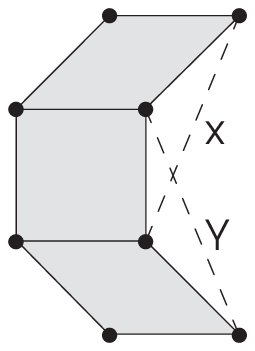}
\quad
&
\quad
\psfrag{x}{$p'$}
\includegraphics{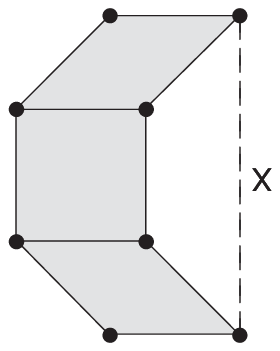}
\quad \\
[3mm]
\text{Case (i)}
& \text{Case (ii)}
& \text{Case (iii)}
& \text{Case (iv)}
& \text{Case (v)}
\end{array}
\]
\caption{After blowing-up a $R_4$-point $p$, there are five possibilities}
\label{fig:r4g}
\end{figure}

Generally, if one blows-up a $R_n$- [resp.\ $S_n$-] point $p$, in
the associated graph to $X$ one builds new 3- and 4- faces (triangles and
quadrangles) over the original chain of length $n$ [resp.\ fork
with $n-1$ teeth] corresponding to the $n$ components of $X$
concurring at $p$. Therefore it is always the case that $G$ is a
deformation retract of the new associated graph.

From this point on there are no more $R_n$ or $S_n$ points ever
appearing in the configuration.  However it may happen that at
intermediate steps of the algorithm, we do not have strict normal
crossings nor Zappatic singularities. If this happens, we still
consider the usual associated graph to the configuration, namely a
vertex for each component, an edge for each connected component of
an intersection between components, and faces for intersections of
three or more components.

Consider Step 1 of algorithm \ref{ssr}.
Each blow-up of an $E_n$-point, where the total space has multiplicity $n$,
has the effect of adding new vertices in the interior of the corresponding $n$-face
and of adding new edges which subdivide the $n$-face.
This does not modify the Betti numbers of the associated graph.

In Step 2 of algorithm \ref{ssr}, the blow-up along a double curve
determines a subdivision of the edge corresponding
to the double curve and a subdivision of the faces adjacent on that edge.

In Step 3 of algorithm \ref{ssr}, the blow-ups at double points of types (a) and (b)
add trees adjacent only to a vertex or an edge,
and again this does not modify the topological properties of the graph.

Resolving a double point of type (c),
one first subdivides the original triangle of vertices $v_1,v_2,v_3$
in three triangles;
then, setting $v_0$ the new vertex,
one adds another vertex $v'_0$ above $v_0$
and three triangles of vertices $v_0,v'_0,v_i$, respectively for $i=1,2,3$.
Clearly the resulting graph retracts back to a subdivision of the original one.

For a double point of type (d),
one subdivides the original quadrangle either in four triangles,
if the exceptional divisor $E$ of the blow-up is a smooth quadric,
or in two triangles and two quadrangles as in Figure \ref{fig:2planes},
if $E$ is the union of two planes.

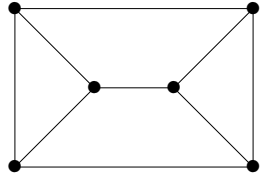
\begin{figure}[ht]
\[
\begin{xy}
0; <30pt,0pt>: 
(0,0)*=0{\bullet};(3,0)*=0{\bullet}**@{-};
(3,2)*=0{\bullet}**@{-};(0,2)*=0{\bullet}**@{-};
(0,0)**@{-};(1,1)*=0{\bullet}**@{-};
(2,1)*=0{\bullet}**@{-};(3,0)**@{-} ,
(0,2);(1,1)**@{-} ,
(2,1);(3,2)**@{-} ,
\end{xy}
\]
\caption{Subdivision of a quadrangle in type (d), case (ii)}
\label{fig:2planes}
\end{figure}

For a double point of type (e),
one subdivides the original triangle either in three triangles,
if the exceptional divisor $E$ of the blow-up is irreducible,
or in a triangle and two quadrangles
as in Figure \ref{fig:2planes2}, if $E$ is reducible.

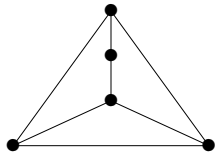
\begin{figure}[ht]
\[
\xymatrix@C=30pt@R=10pt{%
& *=0{\bullet} \ar@{-}[d] \ar@{-}[dddl] \ar@{-}[dddr] & \\
& *=0{\bullet} \ar@{-}[d] & \\
& *=0{\bullet} \ar@{-}[dl] \ar@{-}[dr] & \\
*=0{\bullet} \ar@{-}[rr] & & *=0{\bullet}
}
\]
\caption{Subdivision of a triangle in type (e), case (iii)}
\label{fig:2planes2}
\end{figure}

In all cases, one sees that these modifications,
coming from the resolution of double points of type (c), (d) and (e),
do not change the Betti numbers of the associated graph.

Finally, the blow-ups of Steps 4 and 5
add trees adjacent to a vertex or an edge and again do not modify
the Betti number of the associated graph.
\end{proof}

We are interested not only in $\bar X_\red$ but in $\bar X$ itself.
For each component $i$, let $\mu_i$ be the multiplicity of $\bar X_i$ in $\bar X$.
For the analysis of the semistable reduction, we must understand rather precisely
the components of multiplicity larger than one.

\begin{corollary} (cf. \cite[Corollary 4.11]{CCFMpg}) 
Set $\bar C_{ij}=\bar X_i\cap \bar X_j$ if $\bar X_i$ and $\bar X_j$ meet along a curve,
or $\bar C_{ij}=\emptyset$ otherwise.
If $\mu_i>1$, one has the following possibilities:
\begin{enumerate}
\item[(i)] $\bar X_i$ is a generically ruled surface and the curve
$\sum_{j\ne i} \mu_j \bar C_{ij}$ is generically supported on
a bisection of the ruling.
\item[(ii)] There is a birational morphism $\sigma:\bar X_i\to\Pp^2$ such that
the curve $\sum_{j\ne i} \mu_j \bar C_{ij}$ maps to four distinct lines.
\item[(iii)] $\mu_i=4$ and $\bar X_i$ is a smooth quadric; the curve
$\sum_{j\ne i} \mu_j \bar C_{ij}$ consists of two (multiplicity one) fibres
in one ruling and one double fibre from the other ruling.
\item[(iv)] $\bar X_i$ is a smooth quadric and the curve
$\sum_{j\ne i} \mu_j \bar C_{ij}$ is linearly equivalent
to $\mu_i H$, where $H$ is a plane section of $\bar X_i$.
\item[(v)] There is a birational morphism $\sigma:\bar X_i\to\Pp^2$ such that
the curve $\sum_{j\ne i} \mu_j \bar C_{ij}$ is the total transform
via $\sigma$ of a plane curve of degree $\mu_i$ supported on two distinct lines.
\item[(vi)] $\bar X_i$ is a Hirzebruch surface $\Ff_2$ and
the curve $\sum_{j\ne i} \mu_j \bar C_{ij}$ is
of the form $\mu_i (H+A)$, where $A$ is the $(-2)$-curve
and $H$ is a section of self-intersection 2.
\end{enumerate}
\end{corollary}

\begin{proof}
Following the the Gorenstein reduction algorithm \ref{gr} and the
steps of the normal crossing reduction algorithm \ref{ssr},
one sees that multiple components are not created either in algorithm \ref{gr} or in Step 1 of
algorithm \ref{ssr}. It is possible that a multiple component may be created in Step 2, by blowing-up
a double curve of $X_\red$ which is the intersection of two components that
have multiplicity.  This will create a multiple ruled surface whose
double curve is a bisection, giving case (i).

Multiple components of the central fibre $\bar X$ may arise also
in Step 3 when one blows-up double points of types (c), (d) and (e).
In case (c), two types of multiple components appear.  The first is a plane blown-up
at three collinear points, with multiplicity two; the double curve consists
of the collinearity line, three other general lines, and the three exceptional
divisors counted with multiplicity four; this is case (ii).
The other type of multiple component is a
quadric with multiplicity four, giving case (iii).
This analysis follows from the remark we did at the end of Step 4,
where we showed that the three surfaces coming together to form
this singularity of type (c) each have multiplicity one.

Let $p=X_1\cap X_2\cap X_3\cap X_4$ be a point of type (d),
where $X_1, \ldots, X_4$ are irreducible components of $X_\red$.
One may choose the numbering on the four components
such that $X_1\cup X_2$ and $X_3\cup X_4$
are local complete intersections of $\X$ at $p$,
and moreover the multiplicities satisfy $\mu_1=\mu_2$ and $\mu_3=\mu_4$.
(This is clear at the start, when all multiplicities are one;
and from that point on one proceeds inductively.)
Then the exceptional divisor $E$ appears in the new central fibre
with multiplicity $\mu_1+\mu_3=\mu_2+\mu_4$.
Recall that if $E$ is a smooth quadric, the resolution process stops,
and we have case (iv) above;
while if $E$ is the union of two planes, then both planes
appear with multiplicity $\mu_1+\mu_3$ and we go on inductively;
this gives case (v).

Let now $p=X_1\cap X_2\cap X_3$ be a point of type (e).
As noted above, $X_2\cup X_3$ and $X_1$ are local complete intersections
of $\X$ at $p$.
As above, one may assume that the multiplicities satisfy $\mu_2=\mu_3$.
Then the exceptional divisor $E$ appears in the new central fibre
with multiplicity $\mu_1+\mu_2=\mu_1+\mu_3$.
If $E$ is a smooth quadric, the resolution process stops,
giving case (iv) again.
If $E$ is a quadric cone, then we proceed to blow-up
the vertex of the cone, and therefore the proper transform of $E$
in the final central fibre will be a Hirzebruch surface $\Ff_2$,
which gives the final case (vi).  Finally if $E$ is a pair of planes,
each plane gives rise to a component in case (v).
\end{proof}

Now, by recalling Remark \ref{def:degen2}, we are able to prove the main result of this section:

\begin{proposition}\label{prop:ssr2} (cf. \cite[Proposition 4.12]{CCFMpg}) 
Let $\pi:\X\to\D$ be a degeneration of surfaces
with good Zappatic central fibre $X=\X_0=\bigcup_{i=1}^v X_i$.
Let $\bar\pi:\bar\X\to\D$ be the normal crossing reduction of $\pi$
given by algorithms \ref{gr} and \ref{ssr}
and let $\tilde\pi:\tilde\X\to\D$ be the semistable reduction
of $\bar\pi$ obtained by following the process described in Chapter II of \cite{Kempf}.
Then:
\begin{equation} \label{eq:ssr}
p_\omega(\tilde\X_0)=p_\omega(X).
\end{equation}
\end{proposition}

\begin{proof}
Let $\bar X=\bar\X_0=\sum_{i=1}^w \mu_i \bar X_i$ be the central fibre
of the normal crossing reduction $\bar\pi$.  One has $v\leq w$
and we may assume that $\mu_i=1$ for $1\leq i \leq v$, and that these first
$v$ components are birational to the original components of $X$.
The surface $\bar X$ is a toroidal embedding in $\bar\X$,
in the sense of Definition 1, p.\ 54 of \cite{Kempf}.
To any such a toroidal embedding one can associate a compact
polyhedral complex $\bar\Gamma$ with integral structure
as shown in \cite{Kempf}, pp.\ 71 and 94.
In our present situation, the complex $\bar\Gamma$ is exactly
the associated graph $\bar G$.
The integral structure is recorded by the multiplicities of the components.

By \cite{Kempf}, p.\ 107, there exists a semistable reduction $\tilde\X\to\Delta$
as in Diagram \ref{eq:s}, where the base change $\beta(t)=t^m$ is such that
$m$ is a common multiple of $\mu_1,\ldots,\mu_w$.
Notice that $\tilde\X$ is again a toroidal embedding
of the central fibre $\tilde{X} = \tilde\X_0$.
Denote by $\tilde G$ the associated graph to $\tilde{X}$.
Again by \cite{Kempf}, p.\ 107,
one has that the corresponding polyhedron $\tilde\Gamma$ is a subdivision
of $\bar\Gamma$, in the sense of the definition at p.\ 111 of \cite{Kempf}.
This implies that the CW-complexes $\tilde G$ and $\bar G$ are homeomorphic.
In particular they have the same homology.

Now the central fibre $\tilde{X} = \tilde{\X}_0=\bigcup_{i=1}^u\tilde{X}_i$
is reduced,
with global normal crossings.
One has that $u\ge w$ and, by taking into account the base change,
one may assume that, for $i=1,\ldots,w$, $\tilde X_i$ is birational to
the $\mu_i$-tuple cover of $\bar X_i$,
branched along $\sum_{j\ne i} \mu_j\bar C_{ij}$.

Let us first consider components with $\mu_i=1$.
These include the first $v$ components $\tilde X_i$, $i=1,\ldots,v$,
which correspond to the original components of $X$.
For these components we have
$p_g(\tilde X_i)=p_g(\bar X_i)=p_g(X_i)$, $i=1,\ldots,v$.
There also may be components with $\mu_i=1$ which were introduced
in the normal crossing reduction algorithm.
We have seen in Proposition \ref{XXX}
that all such components have $p_g=0$.
Finally there may be components with $\mu_i=1$ with $i>w$
which have been introduced in the semistable reduction process.
These new surfaces are of two types: they may correspond either to
\begin{enumerate}
\item[(a)] vertices of $\tilde G$ which lie on an edge $\eta$ of $\bar G$; or to
\item[(b)] vertices of $\tilde G$ which lie in the interior of a triangular face of $\bar G$.
\end{enumerate}
We recall that the birational morphism $\tilde{\X}\to\X_\beta$
as in Diagram \ref{eq:s}
is the blow-up of a suitable sheaf of ideals,
cf.\ p.\ 107 of \cite{Kempf}.

Let $\tilde X_j$ be a surface of type (a).
This is an exceptional divisor of such a blow-up with support
on the double curve $\gamma$ of $\bar X$ corresponding to the edge $\eta$.
Then $\tilde X_j$ maps to $\gamma$ with fibres which are rational by the toric
nature of the singularity along $\gamma$.

Suppose that $\tilde X_j$ is of type (b).
Then $\tilde X_j$ is an exceptional divisor appearing in the toric resolution
of a toric singular point.
Therefore $\tilde X_j$ is rational
and moreover it meets the other components along rational curves
(cf., e.g., Section 2.6 in \cite{Fulton}).

Therefore all of these components are rational or ruled,
and hence also have $p_g=0$.

Now let us consider the case $\mu_i>1$.
In this case $\tilde X_i$ is a $\mu_i$-cover
of the surface $\bar X_i$, and such surfaces
were classified in the previous corollary,
along with the double curves which give the branch locus
of the covering.
In each case the cover is easily seen to be rational or ruled.
Hence also for these surfaces one has $p_g=0$.

Since we have shown that the homology of the graphs are the same,
and we have controlled the $p_g$ of the components properly,
the only thing left to prove is that
$\dim({\rm coker}(\Phi_X))=\dim({\rm coker}(\Phi_{\tilde X}))$.

We have already seen that
$\dim({\rm coker} (\Phi_{\bar X_{\red}}))=\dim({\rm coker} (\Phi_{X}))$
in Proposition \ref{XXX}.
The argument here is similar;
it suffices to show that the extra components
$\tilde X_{v+1},\ldots,\tilde X_u$ do not contribute to
$\dim({\rm coker}(\Phi_{\tilde X}))$. These surfaces are either rational or
ruled over a curve $\gamma$. In the rational case, by the proof of
Proposition \ref{ssr} and by the above considerations about toric
resolution of singularities, they meet the other components of $\tilde X$
along rational curves. Hence they do not contribute to
$\dim(\coker\Phi_{\tilde X})$.

In the ruled case, $\tilde X_j$ is a scroll over $\gamma$ and,
by the description of the resolution process, $\tilde X_j$ meets
the other components of $\tilde X$ along curves which are either rational
or isomorphic to $\gamma$.
The same argument as in Proposition \ref{XXX} shows that the
cokernel is unchanged in this case.

Thus the proof is concluded by Theorem \ref{thm:4.pgbis}.
\end{proof}

As a direct consequence, we have the following:

\begin{theorem}\label{thm:pg} (cf. \cite[Theorem 4.14]{CCFMpg}) 
Let $\pi:\X\to\D$ be a degeneration of surfaces with good Zappatic central fibre $X=\X_0$.
Then, for any $t\neq0$, one has:
\[
p_g(\X_t)=p_\omega(X).
\]
\end{theorem}

\begin{proof}
Just consider the semistable reduction $\tilde\pi:\tilde\X\to\D$
as we did before.
One clearly has that $p_g(\X_t)=p_g(\tilde\X_t)$
for $t\neq0$.
Theorem \ref{prop:CS} then implies that $p_g(\tilde\X_t)=p_\omega(\tilde\X_0)$
and finally Proposition \ref{prop:ssr2} concludes that $p_\omega(\tilde\X_0)=p_\omega(X)$.
\end{proof}

\section{The Multiple Point Formula}\label{S:BI}

The aim of this section is to prove a fundamental inequality,
which involves the Zappatic singularities of a given good Zappatic
surface $X$ (see Theorem \ref{thm:BI}), under the hypothesis that
$X$ is the central fibre of a good Zappatic degeneration
as in Definition \ref{def:zappdeg}.
This inequality can be viewed as an extension of the well-known Triple
Point Formula (see Lemma \ref{lem:tpf} and cf.\ \cite{Frie}),
which holds only for semistable degenerations.
As corollaries, we will obtain, among other things, the main result
contained in Zappa's paper \cite{Za3} (cf.\ Section \ref{S:last}).

Let us introduce some notation.

\begin{notation}\label{not:gamma}
Let $X$ be a good Zappatic surface. We denote by:

\begin{itemize}
\item $\gamma=X_1\cap X_2$ the intersection of two irreducible
components $X_1$, $X_2$ of $X$;
\item $F_{\gamma}$ the divisor on $\gamma$ consisting
of the $E_3$-points of $X$ along $\gamma$;
\item $f_n(\gamma)$ the number of $E_n$-points of $X$ along $\gamma$;
in particular, $f_3(\gamma) = \deg(F_{\gamma})$;
\item $r_n(\gamma)$ the number of $R_n$-points of $X$ along $\gamma$;
\item $s_n(\gamma)$ the number of $S_n$-points of $X$ along $\gamma$;
\item $\rho_n(\gamma):=r_n(\gamma)+s_n(\gamma)$,
for $n\ge4$, and $\rho_3(\gamma)=r_3(\gamma)$.
\end{itemize}

If $X$ is the central fibre of a good Zappatic degeneration $\X
\to \D$, we denote by:
\begin{itemize}
\item $D_{\gamma}$ the divisor of $\gamma$ consisting of the
double points of $\X$ along $\gamma$ off the Zappatic singularities of $X$;
\item $d_\gamma = \deg(D_{\gamma})$;
\item $d_\X$ the total number of double points of $\X$
off the Zappatic singularities of $X$.
\end{itemize}
\end{notation}

The main result of this section is the following (cf. \cite[Theorem 7.2]{CCFMk2}):

\begin{theorem}[Multiple Point Formula]\label{thm:BI}
Let $X$ be a surface which is the central fibre of a good Zappatic
degeneration $\X \to \D$. Let $\gamma=X_1\cap X_2$ be the
intersection of two irreducible components $X_1$, $X_2$ of $X$.
Then
\begin{equation}\label{eq:thmBI}
\deg(\N_{\gamma|X_1}) + \deg(\N_{\gamma|X_2}) + f_3(\gamma) -
r_3(\gamma) -\sum_{n\ge4} (\rho_n(\gamma)+f_n(\gamma)) \ge
d_\gamma \ge 0.
\end{equation}
\end{theorem}

In the planar case, one has:

\begin{corollary}\label{cor:BI}
Let $X$ be a surface which is the central fibre of a good, planar
Zappatic degeneration $\X \to \D$. Let $\gamma$ be a double line
of $X$. Then
\begin{equation}\label{eq:corBI1}
2 + f_3(\gamma) - r_3(\gamma) -\sum_{n\ge4}
(\rho_n(\gamma)+f_n(\gamma)) \ge d_\gamma \ge 0.
\end{equation}
Therefore:
\begin{equation}\label{eq:corBI2}
2e+3f_3 - 2 r_3 -\sum_{n\ge4} nf_n -\sum_{n\ge4} (n-1)\rho_n \ge
d_\X \ge 0.
\end{equation}
\end{corollary}

As for Theorem \ref{thm:k2Gmain}, the proof of Theorem
\ref{thm:BI} will be done in several steps, the first of which is
the classical:

\begin{lemma}[Triple Point Formula] \label{lem:tpf}
Let $X$ be a good Zappatic surface with global normal crossings,
which is the central fibre of a good Zappatic degeneration with
smooth total space $\X$. Let $\gamma=X_1\cap X_2$, where $X_1$ and
$X_2$ are irreducible components of $X$. Then:
\begin{equation}\label{eq:tpf1}
\N_{\gamma|X_1}\otimes\N_{\gamma|X_2}\otimes\Oc_\gamma(F_{\gamma})
\cong \Oc_\gamma.
\end{equation}

In particular,
\begin{equation}\label{eq:tpf2}
\deg(\N_{\gamma|X_1}) + \deg(\N_{\gamma|X_2})+ f_3(\gamma) = 0.
\end{equation}

\end{lemma}
\begin{proof}
By Definition \ref{def:zappdeg}, since the total space $\X$ is assumed to be
smooth, the good Zappatic degeneration $\X \to \D$ is
semistable. Let $X = \bigcup_{i=1}^v X_i$. Since $X$ is a Cartier
divisor in $\X$ which is a fibre of the morphism $\X \to \D$, then
$\Oc_X(X) \cong \Oc_X.$
Tensoring by $\Oc_{\gamma}$ gives $\Oc_{\gamma}(X) \cong \Oc_{\gamma}.$
Thus,
\begin{equation}\label{eq:ogamma}
\Oc_{\gamma} \cong \Oc_{\gamma}(X_1) \otimes \Oc_{\gamma}(X_2)
\otimes \Oc_{\gamma}(Y),
\end{equation}
where $Y = \cup_{i=3}^v X_i$. One concludes by observing that in
\eqref{eq:ogamma} one has $\Oc_{\gamma}(X_i) \cong
\N_{\gamma|X_{3-i}}$, $1 \leq i \leq 2$, and $\Oc_{\gamma}(Y)
\cong \Oc_{\gamma}(F_{\gamma})$.
\end{proof}

It is useful to consider the following slightly more general
situation. Let $X$ be a union of surfaces such that $X_{\rm red}$
is a good Zappatic surface with global normal crossings. Then
$X_{\rm red}=\cup_{i=1}^v X_i$ and let $m_i$ be the multiplicity
of $X_i$ in $X$, $i=1,\ldots,v$. Let $\gamma=X_1\cap X_2$ be the
intersection of two irreducible components of $X$. For every point
$p$ of $\gamma$, we define the weight $w(p)$ of $p$ as the
multiplicity $m_i$ of the component $X_i$ such that $p \in
\gamma\cap X_i$.

Of course $w(p)\ne 0$ only for $E_3$-points of $X_{\rm red}$ on
$\gamma$. Then we define the divisor $F_\gamma$ on $\gamma$ as
$$F_\gamma:=\sum_p w(p)p.$$

The same proof of Lemma \ref{lem:tpf} shows the following:

\begin{lemma}[Generalized Triple Point Formula] \label{lem:gtpf}
Let $X$ be a surface such that $X_{\rm red}=\cup_i X_i$ is a good
Zappatic surface with global normal crossings. Let $m_i$ be the
multiplicity of $X_i$ in $X$. Assume that $X$ is the central fibre
of a degeneration $\X \to \D$ with smooth total space $\X$. Let
$\gamma=X_1\cap X_2$, where $X_1$ and $X_2$ are irreducible
components of $X_{\rm red}$. Then:
\begin{equation}\label{eq:gtpf1}
\N_{\gamma|X_1}^{\otimes m_2}\otimes \N_{\gamma|X_2}^{\otimes
m_1}\otimes \Oc_\gamma(F_{\gamma}) \cong \Oc_\gamma.
\end{equation}

In particular,
\begin{equation}\label{eq:gtpf2}
m_2\deg(\N_{\gamma|X_1}) + m_1\deg(\N_{\gamma|X_2})+ \deg
(F_\gamma) = 0.
\end{equation}
\end{lemma}

The second step is given by the following result (cf. \cite[Proposition 7.14]{CCFMk2}):

\begin{proposition}\label{prop:step2}
Let $X$ be a good Zappatic surface with global normal crossings,
which is the central fibre of a good Zappatic degeneration $\X \to
\D$. Let $\gamma=X_1\cap X_2$, where $X_1$ and $X_2$ are
irreducible components of $X$. Then:
\begin{equation}\label{eq:step2-1}
\N_{\gamma|X_1}\otimes\N_{\gamma|X_2}\otimes\Oc_\gamma(F_{\gamma})
\cong \Oc_\gamma (D_{\gamma}).
\end{equation}

In particular,
\begin{equation}\label{eq:step2-2}
\deg(\N_{\gamma|X_1}) + \deg(\N_{\gamma|X_2})+ f_3(\gamma) =
d_\gamma.
\end{equation}
\end{proposition}
\begin{proof}
By the very definition of good Zappatic degeneration, the total
space $\X$ is smooth except for ordinary double points along the
double locus of $X$, which are not the $E_3$-points of $X$. We can
modify the total space $\X$ and make it smooth by blowing-up its
double points.

Since the computations are of local nature, we can focus on the
case of $\X$ having only one double point $p$ on $\gamma$. We
blow-up the point $p$ in $\X$ to get a new total space $\X'$,
which is smooth. Notice that, according to our hypotheses, the
exceptional divisor $E:= E_{\X, p} = \Pp(T_{\X , p})$ is
isomorphic to a smooth quadric in $\Pp^3$ (see Figure
\ref{fig:doppio}).

\begin{figure}[ht]
\[
\begin{xy}
<40pt,0pt>:
(0,0);(3.5,0)**@{-}="d";(3.5,1)**@{-};(0,1)**@{-}="a";(0,0)**@{-}
, (2,0);(2,1)**@{-}="c", (1,0)="b";(1,1)**@{-}?(.7)="p"*{\bullet}
?(0)+<1pt,1pt>*!LD{\st \gamma}, "p"+<3pt,0pt>*!L{\st p} ,
"a"."b"!C*{\dt X_1} , "b"."c"!C*{\dt X_2} , "c"."d"!C*{\dt Y}
\end{xy}
\qquad
\raisebox{20pt}{%
  $\xleftarrow{\text{blow-up $p$}}$
} \qquad
\begin{xy}
<40pt,0pt>: (0,0); (4,0)**@{-}="d"; (4,1)**@{-};
(1.75,1)**@{-}="e"; (1.25,.5)**@{-}="p" ?(.45)*!RD[@!45]{\st 0}
?(.7)*[@!45]!L{\st -1} , "p"; (.75,1)**@{-}="f"
?(.55)*!DL[@!-45]{\st 0} ?(.4)*[@!-45]!R{\st -1} , "f";
(0,1)**@{-}="a"; (0,0)**@{-} , "f"; (1.25,1.5)**@{-}; "e"**@{-} ,
(2.5,0); (2.5,1)**@{-}="c", (1.25,0)="b"; "p" **@{-}
?(0)+<1pt,1pt>*!LD{\st \gamma'} , "e"."f"!C*+!D{\dt E} ,
"a"."b"!C+<-3pt,-2pt>*{\dt X'_1} , "b"."c"!C+<3pt,-2pt>*{\dt X'_2}
, "c"."d"!C*{\dt Y}
\end{xy}
\]
  \caption{Blowing-up an ordinary double point of $\X$}\label{fig:doppio}
\end{figure}
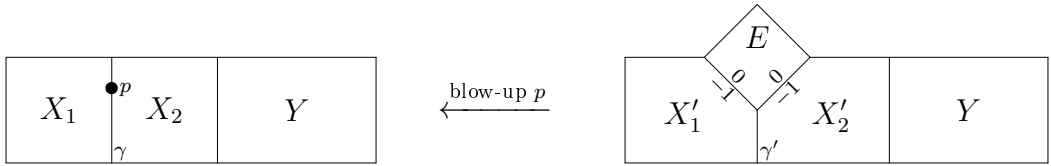

The proper trasform of $X$ is:
$$X'= X'_1+ X'_2+ Y$$
where $X'_1, \; X'_2$ are the proper transforms of $X_1$, $X_2$,
respectively. Let $\gamma'$ be the intersection of $X'_1$ and
$X'_2$, which is clearly isomorphic to $\gamma$. Let $p_1$ be the
intersection of $\gamma' $ with $E$.

Since $\X'$ is smooth, we can apply Lemma \ref{lem:gtpf} to
$\gamma'$. Therefore, by \eqref{eq:gtpf1}, we get $$\Oc_{\gamma'}
\cong \N_{\gamma'|X'_1} \otimes \N_{\gamma'|X'_2} \otimes
\Oc_{\gamma'} (F_{\gamma'}).$$In the isomorphism between $\gamma'$
and $\gamma$, one has: $$\Oc_{\gamma'} (F_{\gamma'} - p_1) \cong
\Oc_{\gamma}(F_{\gamma}), \; \; \N_{\gamma'|X'_i} \cong
\N_{\gamma|X_i}\otimes \Oc_{\gamma} (-p), \;\; 1 \leq i \leq 2.$$
Putting all this together, one has the result.
\end{proof}

Taking into account Lemma \ref{lem:gtpf}, the same proof of
Proposition \ref{prop:step2} gives the following result:

\begin{corollary}\label{cor:step2} (cf. \cite[Corollary 7.17]{CCFMk2})
Let $X$ be a surface such that $X_{\rm red}=\cup_i X_i$ is
a good Zappatic surface with global normal
crossings. Let $m_i$ be the multiplicity of $X_i$ in $X$. Assume
that $X$ is the central fibre of a degeneration $\X \to \D$ with total
space $\X$ having at most ordinary double points outside the
Zappatic singularities of $X_{\rm red}$.

Let $\gamma=X_1\cap X_2$, where $X_1$ and $X_2$ are irreducible
components of $X_{\rm red}$. Then:
\begin{equation}\label{eq:gtpf1step2}
\N_{\gamma|X_1}^{\otimes m_2}\otimes \N_{\gamma|X_2}^{\otimes
m_1}\otimes \Oc_\gamma(F_{\gamma}) \cong
\Oc_\gamma(D_\gamma)^{\otimes(m_1+m_2)}.
\end{equation}
In particular,
\begin{equation}\label{eq:gtpf2step2}
m_2\deg(\N_{\gamma|X_1}) + m_1\deg(\N_{\gamma|X_2})+ \deg
(F_\gamma) = (m_1+m_2)d_\gamma.
\end{equation}
\end{corollary}

Now we can come to the:

\begin{proof}[Proof of Theorem \ref{thm:BI}]
Recall that, by  Definition \ref{def:zappdeg} of Zappatic
degenerations, the total space $\X$ has only isolated
singularities. We want to apply Corollary \ref{cor:step2} after
having resolved the singularities of the total space $\X$ at the
Zappatic singularities of the central fibre $X$, i.e.\ at the
$R_n$-points of $X$, for $n \geq 3$, and at the $E_n$- and
$S_n$-points of $X$, for $n \geq 4$.

Now we briefly describe the resolution process, which will become
even clearer in the second part of the proof, when we will enter
into the details of the proof of Formula \eqref{eq:thmBI}.

Following the blowing-up process of Algorithm \ref{gr} and of Proposition
\ref{ssr} - $(1)$ at the $R_n$- and $S_n$-points
of the central fibre $X$, as described in details in Section \ref{S:5},
one gets a degeneration such that
the total space is Gorenstein, with isolated singularities, and
the central fibre is a Zappatic surface with only $E_n$-points.

The degeneration will not be Zappatic, if the double points of
the total space occurring along the double curves, off the
Zappatic singularities, are not ordinary. According to our
hypotheses, this cannot happen along the proper transform of the
double curves of the original central fibre. All these
non-ordinary double points can be resolved with finitely many
subsequent blow-ups and they will play no role in the computation
of Formula \eqref{eq:thmBI}.

Recall that the total space $\X$ is smooth at the $E_3$-points of
the central fibre, whereas $\X$ has multiplicity either 2 or 4 at
an $E_4$-point of $X$. Thus, we can consider only $E_n$-points
$p\in X$, for $n\ge4$.

By Proposition \ref{prop:mgmz3f}, $p$ is a quasi-minimal singularity for $\X$, unless $n=4$
and $\mult_p(\X)=2$. In the latter case, this singularity is resolved by a sequence
of blowing-ups at isolated double points.

Assume now that $p$ is a quasi-minimal singularity for $\X$. Let us
blow-up $\X$ at $p$ and let $E'$ be the exceptional divisor. Since a hyperplane section of $E'$
is $C_{E_n}$, the possible configurations of $E'$ are described in Proposition \ref{prop:hypsect}, (iii).

If $E'$ is irreducible, that is case (iii.a) of Proposition \ref{prop:hypsect},
then $E'$ has at most isolated rational double points,
where the new total space is either smooth or it has a double point. This
can be resolved by finitely many blowing-ups at analogous double points.

Suppose we are in case (iii.b) of Proposition \ref{prop:hypsect}.
If $E'$ has global normal crossings, then the desingularization process
proceeds exactly as before.

If $E'$ does not have global
normal crossings then, either $E'$ has a component which
is a quadric cone or the two components of $E'$ meet along a
singular conic. In the former case, the new total space
has a double point at the vertex of the cone.
In the latter case, the total
space is either smooth or it has an isolated double point at
the singular point of the conic. In either cases, one resolves
the singularities by a sequence of blowing-ups as before.

Suppose finally we are in case (iii.c) of Proposition \ref{prop:hypsect},
i.e.\ the new central fibre is a Zappatic surface with one point $p'$ of type $E_m$,
with $m \le n$. Then we can proceed by induction on $n$. Note that if an exceptional divisor has an
$E_3$-point $p''$, then $p''$ is either a smooth, or a double, or a triple point for the total space.
In the latter two cases, we go on by blowing-up $p''$. After finitely many blow-ups (by Definition \ref{def:zappdeg},
cf.\ Proposition 3.4.13 in \cite{Kollar}), we get a central fibre which might be non-reduced,
but its support has only global normal
crossings, and the total space has at most ordinary double points
off the $E_3$-points of the reduced part of the central fibre.

Now we are in position to apply Corollary \ref{cor:step2}.
In order to do this, we have to understand the relations
between the invariants of a double curve of the original Zappatic surface $X$
and the invariants appearing in Formula \eqref{eq:gtpf2step2}
for the double curve of the strict transform of $X$.

Since all the computations are of local nature, we may
assume that $X$ has a single Zappatic singularity $p$, which is
not an $E_3$-point. We will prove the theorem in this case. The
general formula will follow by iterating these considerations for
each Zappatic singularity of $\X$.

Let $X_1$, $X_2$ be irreducible components of $X$ containing $p$
and let $\gamma$ be their intersection.
As we saw in the above resolution process, we blow-up $\X$ at
$p$. We obtain a new total space $\X'$, with the exceptional
divisor $E' := E_{\X, p}=\Pp(T_{\X, p})$ and the proper transform
$X'_1, \; X'_2$ of $X_1$, $X_2$. Let $\gamma'$ be the intersection
of  $X'_1, \; X'_2$. We remark that $\gamma' \cong \gamma$ (see
Figure \ref{fig:BI1}).

\begin{figure}[ht]
\[
    \xymatrix@C=40pt@R=50pt{%
      *=0{} \ar@{-}[rrr] \ar@{-}[dd] & *=0{} & *=0{} & *=0{} \ar@{-}[dd] \\
      *=0{} \ar@{-}[r]_{E_1} & *=0{\bullet} \ar@{-}[r]_{E_2} \ar@{-
}[d]_>>>{\gamma'} & *=0{} \ar@{-}[r] \ar@{-}[d] & *=0{} \\
      *=0{} \ar@{-}[rrr] & *=0{} & *=0{} & *=0{}
      \save "1,1"."2,4"!C*\txt{$E'$} \restore
      \save "2,1"."3,2"!C-<0ex,.5ex>*\txt{$X'_1$} \restore
      \save "2,2"."3,3"!C-<0ex,.5ex>*\txt{$X'_2$} \restore
      \save "2,3"."3,4"!C-<0ex,.5ex>*\txt{$Y'$} \restore
      \save "2,2"+<0pt,6pt>*\txt{$p_1$} \restore
    }
    \qquad
    \raisebox{-50pt}{$ %
    \xrightarrow{\text{blow-up }p}
    $}
    \qquad
    \raisebox{-50pt}{%
    \xymatrix@C=30pt@R=50pt{%
      *=0{} \ar@{-}[rrrr] 
        \ar@{-}[dr] & *=0{} & *=0{\bullet} \ar@{-}[dr] \ar@{-}[dl] & *=0{} &
*=0{} \ar@{-}[dl] \\
      *=0{} & *=0{} \ar@{-}[rr] \ar@{}[ur]_<<<<{\gamma} & *=0{} & *=0{} & *=0{}
      \save "1,3"+<0ex,1.5ex> *{\dt p} \restore
      \save "1,1"."2,3"!C+<0pt,7pt> *{\dt X_1} \restore
      \save "1,2"."2,4"!C-<0pt,10pt> *{\dt X_2} \restore
      \save "1,3"."2,5"!C+<0pt,7pt> *{\dt Y} \restore
    }
    }
  \]
\caption{Blowing-up $\X$ at $p$}\label{fig:BI1}
\end{figure}
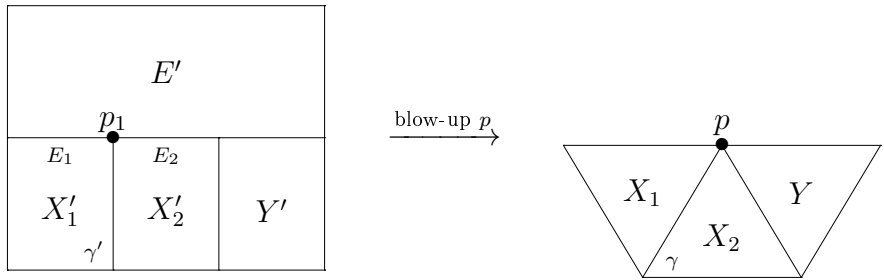

Notice that $\X'$ might have Zappatic singularities off $\gamma'$.
These will not affect our considerations. Therefore, we can assume
that there are no singularities of $\X'$ of this sort. Thus, the
only point of $\X'$ we have to take care of is $p_1 := E' \cap
\gamma'$.

If $p_1$ is smooth for $E'$, then it must be smooth also for $\X'$.
Moreover, if $p_1$ is singular for $E'$, then $p_1$ is a double point of $E'$ as it follows from the
above resolution process and from Proposition  \ref{prop:hypsect}. Therefore, $p_1$
is at most double also for $\X'$;
since $p_1$ is a quasi-minimal, Gorenstein singularity of multiplicity 4 for the central fibre of $\X'$,
then $p_1$ is a double point of $\X'$ by Proposition \ref{prop:mgmz3f}.

Thus there are two cases to be considered: either
\begin{itemize}
\item[(i)] $p_1$ is smooth for both $E'$ and $\X'$, or
\item[(ii)] $p_1$ is a double point for both $E'$ and $\X'$.
\end{itemize}

In case (i), the central fibre of $\X'$ is $\X'_0 = X'_1 \cup X'_2
\cup Y'\cup E'$ and we are in position to use the enumerative information
\eqref{eq:step2-2} from Proposition \ref{prop:step2} which reads:
$$\deg(\N_{\gamma'|X'_1}) + \deg(\N_{\gamma'|X'_2})+
f_3(\gamma') = d_{\gamma'}.$$
Observe that $f_3(\gamma')$ is the
number of $E_3$-points of the central fibre $\X'_0$ of $\X'$ along $\gamma'$, therefore
$$f_3(\gamma') = f_3(\gamma)+1.$$
On the other hand:
$$\deg(\N_{\gamma'|X'_i}) = \deg (\N_{{\gamma}|X_i}) -1, \;\; 1 \leq i \leq 2.$$
Finally,
$$d_{\gamma} = d_{\gamma'}$$
and therefore we have
\begin{equation}\label{eq:+1}
\deg(\N_{\gamma|X_1}) + \deg(\N_{\gamma|X_2})+ f_3(\gamma) - 1 =
d_{\gamma}
\end{equation}
which proves the theorem in this case (i).

\medskip
Consider now case (ii), i.e.\ $p_1$ is a double point for both $E'$ and $\X'$.

If $p_1$ is an ordinary double point for $\X'$,
we blow-up $\X'$ at $p_1$ and we get
a new total space $\X''$.
Let $X''_1$, $X''_2$ be the proper transforms of $X'_1$, $X'_2$, respectively,
and let $\gamma''$ be the intersection of $X''_1$ and $X''_2$,
which is isomorphic to $\gamma$.
Notice that $\X''$ is smooth and the exceptional divisor $E''$ is a
smooth quadric (see Figure \ref{fig:1flam}).

\begin{figure}[ht]
\[
     \xymatrix@C=60pt@R=60pt{%
      *=0{} \ar@{-}[d] \ar@{-}[];[r]-<20pt,0pt>
        & *=0{E''} \ar@{-}[]-<0pt,20pt>;[d]_>>>{\gamma''} \ar@{-}[]+<20pt,0pt>;[r]
         \ar@{-}[]-<20pt,0pt>;[]-<0pt,20pt>   \ar@{-}[]+<20pt,0pt>;[]-<0pt,20pt>
         \ar@{-}[]-<20pt,0pt>;[]+<0pt,20pt>   \ar@{-}[]+<20pt,0pt>;[]+<0pt,20pt>
       & *=0{} \ar@{-}[r] \ar@{-}[d] & *=0{} \ar@{-}[d] \\
      *=0{} \ar@{-}[rrr] & *=0{} & *=0{} & *=0{}
      \save "1,1"."2,2"!C-<0ex,.5ex>*\txt{$X''_1$} \restore
      \save "1,2"."2,3"!C-<0ex,.5ex>*\txt{$X''_2$} \restore
      \save "1,3"."2,4"!C-<0ex,.5ex>*\txt{$Y'$} \restore
      \save "1,2"-<0pt,20pt>*\txt{$\bullet$} \restore
      \save "1,2"+<7pt,-25pt>*\txt{$p_2$} \restore
    }
\]
\caption{Blowing-up $\X'$ at $p_1$ when $p_1$ is ordinary
for both $\X'$ and $E'$}\label{fig:1flam}
\end{figure}
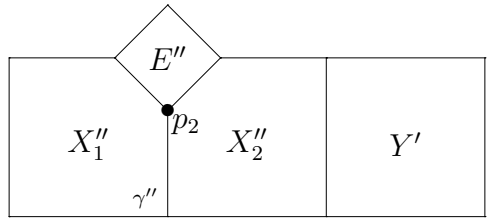

We remark that the central fibre of $\X''$ is now non-reduced,
since it contains $E''$ with multiplicity 2.
Thus we apply Corollary \ref{cor:step2} and we get
$$\Oc_{\gamma''} \cong \N_{\gamma''|X''_1} \otimes
\N_{\gamma''|X''_2} \otimes \Oc_{\gamma''} (F_{\gamma''}).$$
Since,
$$\deg(\N_{\gamma''|X''_1}) = \deg(\N_{\gamma|X_i}) -2, \; i=1,2,
\qquad \deg F_{\gamma''} = f_3(\gamma) + 2,$$ then,
\begin{equation}\label{eq:bii+2}
\deg(\N_{\gamma|X_1}) + \deg(\N_{\gamma|X_2})+  f_3(\gamma) - 1 =
d_{\gamma} + 1
> d_{\gamma}.
\end{equation}

If the point $p_1$ is not an ordinary double point, we again
blow-up $p_1$ as above. Now the exceptional divisor $E''$ of
$\X''$ is a singular quadric in $\Pp^3$, which
can only be either a quadric cone or it has to consist of
two distinct planes $E''_1$, $E''_2$.
Note that if $p_1$ lies on a double line of $E'$ (i.e.\ $p_1$ is in the intersection
of two irreducible components of $E'$),
then only the latter case occurs since $E''$ has to contain a curve $C_{E_4}$.

Let $p_2=E''\cap \gamma''$.
In the former case, if $p_2$ is not the vertex of the quadric cone,
then the total space $\X''$ is smooth at $p_2$
and we can apply Corollary \ref{cor:step2}
and we get \eqref{eq:bii+2} as before.

If $p_2$ is the vertex of the quadric cone,
then $p_2$ is a double point of $\X''$ and
we can go on blowing-up $\X''$ at $p_2$.
This blow-up procedure stops after finitely many, say $h$, steps
and one sees that Formula \eqref{eq:bii+2} has to be replaced by
\begin{equation}\label{eq:bii+h}
\deg(\N_{\gamma|X_1}) + \deg(\N_{\gamma|X_2})+  f_3(\gamma) - 1 =
d_{\gamma} + h
> d_{\gamma}.
\end{equation}

In the latter case, i.e.\ if $E''$ consists of two planes $E''_1$ and $E''_2$,
let $\lambda$ be the intersection line of $E''_1$ and $E''_2$. If $p_2$ does not belong to
$\lambda$ (see Figure \ref{fig:1flb}), then $p_2$ is a
smooth point of the total space $\X''$,
therefore we can apply Corollary
\ref{cor:step2} and we get again Formula \eqref{eq:bii+2}.

\begin{figure}[ht]
\[
     \xymatrix@C=70pt@R=60pt{%
      *=0{} \ar@{-}[d] \ar@{-}[];[r]-<30pt,0pt>
       & *=0{} \ar@{-}[]-<0pt,30pt>;[d]_>>>{\gamma''} \ar@{-}[]+<30pt,0pt>;[r]
         \ar@{-}[]-<30pt,0pt>;[]-<0pt,30pt>   \ar@{-}[]+<30pt,0pt>;[]-<0pt,30pt>
         \ar@{-}[]-<30pt,0pt>;[]+<0pt,30pt>   \ar@{-}[]+<30pt,0pt>;[]+<0pt,30pt>
         \ar@{-}[]+<30pt,0pt>;[]-<30pt,0pt>
       & *=0{} \ar@{-}[r] \ar@{-}[d] & *=0{} \ar@{-}[d] \\
      *=0{} \ar@{-}[rrr] & *=0{} & *=0{} & *=0{}
      \save "1,1"."2,2"!C-<0ex,.5ex>*\txt{$X''_1$} \restore
      \save "1,2"."2,3"!C-<0ex,.5ex>*\txt{$X''_2$} \restore
      \save "1,3"."2,4"!C-<0ex,.5ex>*\txt{$Y'$} \restore
      \save "1,2"-<0pt,30pt>*\txt{$\bullet$} \restore
      \save "1,2"+<9pt,-35pt>*\txt{$p_2$} \restore
      \save "1,2"-<0pt,15pt>*\txt{$E''_1$} \restore
      \save "1,2"+<0pt,15pt>*\txt{$E''_2$} \restore
      \save "1,2"+<-10pt,-3pt>*\txt{$\st\lambda$} \restore
    }
\]
\caption{$E''$ splits in two planes $E''_1$, $E''_2$ and
$p_2\not\in E''_1\cap E''_2$}\label{fig:1flb}
\end{figure}

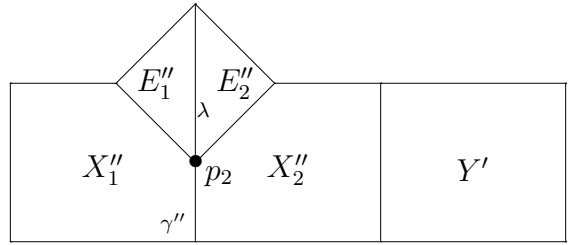
\begin{figure}[ht]
\[
     \xymatrix@C=70pt@R=60pt{%
      *=0{} \ar@{-}[d] \ar@{-}[];[r]-<30pt,0pt>
       & *=0{} \ar@{-}[]-<0pt,30pt>;[d]_>>>{\gamma''} \ar@{-}[]+<30pt,0pt>;[r]
         \ar@{-}[]-<30pt,0pt>;[]-<0pt,30pt>   \ar@{-}[]+<30pt,0pt>;[]-<0pt,30pt>
         \ar@{-}[]-<30pt,0pt>;[]+<0pt,30pt>   \ar@{-}[]+<30pt,0pt>;[]+<0pt,30pt>
         \ar@{-}[]+<0pt,30pt>;[]-<0pt,30pt>
       & *=0{} \ar@{-}[r] \ar@{-}[d] & *=0{} \ar@{-}[d] \\
      *=0{} \ar@{-}[rrr] & *=0{} & *=0{} & *=0{}
      \save "1,1"."2,2"!C-<0ex,.5ex>*\txt{$X''_1$} \restore
      \save "1,2"."2,3"!C-<0ex,.5ex>*\txt{$X''_2$} \restore
      \save "1,3"."2,4"!C-<0ex,.5ex>*\txt{$Y'$} \restore
      \save "1,2"-<0pt,30pt>*\txt{$\bullet$} \restore
      \save "1,2"+<9pt,-35pt>*\txt{$p_2$} \restore
      \save "1,2"-<15pt,0pt>*\txt{$E''_1$} \restore
      \save "1,2"+<15pt,0pt>*\txt{$E''_2$} \restore
      \save "1,2"+<3pt,-10pt>*\txt{$\st\lambda$} \restore
    }
\]
\caption{$E'''$ splits in two planes $E'''_1$, $E'''_2$ and
$p_3\in E'''_1\cap E'''_2$}\label{fig:1flc}
\end{figure}

If $p_2$ lies on $\lambda$, then $p_2$ is a double point for the
total space $\X''$ (see Figure \ref{fig:1flc}). We can thus
iterate the above procedure until the process terminates after finitely many, say $h$, steps
by getting rid of the singularities which are infinitely near to $p$
along $\gamma$. At the end, one gets again Formula \eqref{eq:bii+h}.
\end{proof}

\begin{remark}\label{rem:epsilongamma}
We observe that the proof of Theorem \ref{thm:BI} proves a
stronger result than what we stated in \eqref{eq:thmBI}.
Indeed, the idea of the proof is that we blow-up the total space
$\X$ at each Zappatic singularity $p$ in a sequence
of singular points $p$, $p_1$, $p_2$, \ldots, $p_{h_p}$, each infinitely near one
to the other along $\gamma$. Note that $p_i$, $i=1,\ldots,h_p$,
is a double point for the total space.

The above proof shows that the first inequality
in \eqref{eq:thmBI} is an equality if and only if each Zappatic
singularity of $\X$ has no infinitely near singular point.
Moreover \eqref{eq:bii+h} implies that
\[
\deg(\N_{\gamma|X_1}) + \deg(\N_{\gamma|X_2}) + f_3(\gamma) -
r_3(\gamma) - \sum_{n \geq 4} (\rho_n(\gamma) + f_n(\gamma)) =
d_\gamma + \sum_{p\in\gamma} h_p.
\]
In other words, as natural, every infinitely near double point
along $\gamma$ counts as a double point of the original total
space along $\gamma$.
\end{remark}

\section{On some results of Zappa}\label{S:last}

In \cite{Za1, Za1b, Za1c, Za2, Za2b, Za3, Za3b}, Zappa considered
degenerations of projective surfaces to a planar Zappatic surface
with only $R_3$-, $S_4$- and $E_3$-points. One of the results of
Zappa's analysis is that the invariants of a surface admitting a
good planar Zappatic degeneration with mild singularities are
severely restricted. In fact, translated in modern terms, his main
result in \cite{Za3} can be read as follows:

\begin{theorem}[Zappa]\label{thm:1Zappa}
Let $\X \to \D$ be a good, planar Zappatic degeneration, where the
central fibre $\X_0=X$ has at most $R_3$- and $E_3$-points. Then,
for $t\ne0$, one has
\begin{equation}\label{eq:1Zappa}
K^2:=K_{\X_t}^2\leq 8\chi+1-g,
\end{equation}
where $\chi=\chi(\Oc_{\X_t})$ and $g$ is the sectional genus of $\X_t$.
\end{theorem}

Theorem \ref{thm:1Zappa} has the following interesting
consequence:

\begin{corollary}[Zappa]\label{cor:1Zappa}
If $\X \to \D$ is a good, planar Zappatic degeneration of a scroll $\X_t$
of sectional genus $g\geq 2$ to $\X_0= X$, then $X$ has worse
singularities than $R_3$- and $E_3$-points.
\end{corollary}

\begin{proof}
For a scroll of genus $g$ one has $8\chi+1-g-K^2=1-g$.
\end{proof}

Actually Zappa conjectured that for most of the surfaces the
inequality $K^2 \le 8\chi+1$ should hold and even proposed a
plausibility argument for this. As well-known, the correct bound
for all the surfaces is $K^2 \le 9\chi$, proved by Miyaoka and Yau
(see \cite{Miyaoka, Yau}) several decades after Zappa.

We will see in a moment that Theorem \ref{thm:1Zappa} can be
proved as consequence of the computation of $K^2$ (see Theorem
\ref{thm:k2Gmain}) and the Multiple Point Formula (see Theorem
\ref{thm:BI}).

Actually, Theorems \ref{thm:k2Gmain} and \ref{thm:BI} can be used
to prove a stronger result than Theorem \ref{thm:1Zappa}; indeed:

\begin{theorem}\label{thm:1ZappaBis} (cf. \cite[Theorem 8.4]{CCFMk2})
Let $\X \to \D$ be a good, planar Zappatic degeneration, where the
central fibre $\X_0=X$ has at most $R_3$-, $E_3$-, $E_4$- and
$E_5$-points. Then
\begin{equation}\label{eq:1ZappaBis}
K^2\leq 8\chi+1-g.
\end{equation}
Moreover, the equality holds in \eqref{eq:1ZappaBis} if and only
if $\X_t$ is either the Veronese surface in $\Pp^5$ degenerating
to four planes with associated graph $S_4$ (i.e.\ with three
$R_3$-points, see Figure \ref{fig:zappequal}.a), or an elliptic
scroll of degree $n\geq 5$ in $\Pp^{n-1}$ degenerating to $n$
planes with associated graph a cycle $E_n$ (see Figure
\ref{fig:zappequal}.b).

Furthermore, if $\X_t$ is a surface of general type, then
\begin{equation}\label{eq:1ZappaTer}
K^2 < 8\chi -g.
\end{equation}
\end{theorem}

\begin{figure}[ht]
\[
\begin{array}{cc}
\qquad
\begin{xy}
0; <35pt,0pt>: 
{\xypolygon3{~>{}}}, "1"*=0{\bullet};"0"*{\bullet}**@{-} ,
"2"*=0{\bullet};"0"**@{-} , "3"*=0{\bullet};"0"**@{-} ,
\end{xy}
\qquad & \qquad
\begin{xy}
0; <35pt,0pt>: 
{\xypolygon9{\bullet}},
\end{xy}
\qquad  \\[35pt]
\text{(a)}  & \text{(b)}
\end{array}
\]
\caption{}\label{fig:zappequal}
\end{figure}

\begin{proof}
Notice that if $X$ has at most $R_3$-, $E_3$-, $E_4$- and
$E_5$-points, then Formulas \eqref{eq:K2Gbounds2} and
\eqref{eq:pK2Gbounds} give $K^2=9v-10e+6f_3+8f_4+10f_5 + r_3$.
Thus, by \eqref{eq:gplanar} and \eqref{eq:chiplan}, one gets
\begin{align*}
& 8\chi+1-g-K^2 = 8v-8e+8f_3+8f_4+8f_5+1-(e-v+1)-K^2 = e-r_3+2f_3-2f_5 = \\
& = \frac{1}{2}\,(2e-2r_3+3f_3-4f_4-5f_5)
+\frac{1}{2}\,f_3+2f_4+\frac{1}{2}\,f_5 \stackrel{(*)}{\ge}
\frac{1}{2}\,f_3+2f_4+\frac{1}{2}\,f_5 \ge 0
\end{align*}
where the inequality $(*)$ follows from \eqref{eq:corBI2}. This
proves Formula \eqref{eq:1ZappaBis} (and Theorem
\ref{thm:1Zappa}).

If $K^2=8\chi+1-g$, then $(*)$ is an equality, hence
$f_3=f_4=f_5=0$ and $e=r_3$. Therefore, by Formula \eqref{eq:r3},
we get
\begin{equation}\label{eq:contir3}
\sum_i w_i (w_i -1) = 2 r_3 = 2 e,
\end{equation}where $w_i$ denotes the valence of the vertex $v_i$ in the
graph $G_X$. By definition of valence, the right-hand-side of
\eqref{eq:contir3} equals $\sum_i w_i$. Therefore, we get
\begin{equation}\label{eq:contir32}
\sum_i w_i (w_i -2) = 0.
\end{equation}

If $w_i \geq 2$, for each $1 \leq i \leq v$, one easily shows that
only the cycle as in Figure \ref{fig:zappequal} (b) is possible.
This gives $$\chi=0, \; K^2 = 0 , \; g=1,$$which implies that
$\X_t$ is an elliptic scroll.

Easy combinatorial computations show that, if there is a vertex
with valence $w_i \ne 2$, then there is exactly one vertex with
valence 3 and three vertex of valence 1. Such a graph, with $v$
vertices, is associated to a planar Zappatic surface of degree $v$
in $\Pp^{v+1}$ with
$$\chi=0, \; p_g= 0 , \; g=0.$$
Thus, by hypothesis, $K^2=9$ and, by properties of projective
surfaces, the only possibility is that $v=4$, $G_X$ is as in
Figure \ref{fig:zappequal} (a) and $\X_t$ is the Veronese surface
in $\Pp^5$.

Suppose now that $\X_t$ is of general type. Then $\chi \ge 1$ and
$v = \deg(\X_t) < 2g-2$. Formulas \eqref{eq:gplanar} and
\eqref{eq:chiplan} imply that $\chi=f-g+1 \ge 1$, thus $f \ge g >
v/2 + 1$. Clearly $v \ge 4$, hence $f \ge 3$. Proceeding as at the
beginning of the proof, we have that:
\[
8\chi-g-K^2 \ge \frac{1}{2}\,f_3+2f_4+\frac{1}{2}\,f_5 - 1 \ge
\frac{1}{2}\,f - 1 > 0,
\]
or equivalently $K^2<8\chi-g$.
\end{proof}

\begin{remark}\label{rem:polito}
By following the same argument of the proof of Theorem
\ref{thm:1ZappaBis}, one can list all the graphs and the
corresponding smooth projective surfaces in the degeneration, for
which $K^2 = 8\chi -g$. For example, one can find $\X_t$ as a
rational normal scroll of degree $n$ in $\Pp^{n+1}$ degenerating
to $n$ planes with associated graph a chain $R_n$. On the other
hand, one can also have a del Pezzo surface of degree $7$ in
$\Pp^7$.
\end{remark}

Let us state some applications of Theorem \ref{thm:1ZappaBis}.

\begin{corollary}\label{cor:1Zappabis}
If $\X$ is a good, planar Zappatic degeneration of a scroll $\X_t$
of sectional genus $g\geq 2$ to $\X_0=X$, then $X$ has worse
singularities than $R_3$-, $E_3$-, $E_4$- and $E_5$-points.
\end{corollary}

\begin{corollary}\label{cor:torino}
If $\X$ is a good, planar Zappatic degeneration of a del Pezzo
surface $\X_t$ of degree $8$ in $\Pp^8$ to $\X_0=X$, then $X$ has
worse singularities than $R_3$-, $E_3$-, $E_4$- and $E_5$-points.
\end{corollary}

\begin{proof}
Just note that $K^2=8$ and $\chi=g=1$, thus $\X_t$ satisfies the
equality in \eqref{eq:1ZappaBis}.
\end{proof}

\begin{corollary}\label{cor:2Zappa}
If $\X$ is a good, planar Zappatic degeneration of a minimal
surface of general type $\X_t$ to $\X_0 = X$ with at most $R_3$-,
$E_3$-, $E_4$- and $E_5$-points, then
$$g\leq 6\chi+5.$$
\end{corollary}

\begin{proof}
It directly follows from \eqref{eq:1ZappaTer} and Noether's
inequality, i.e.\ $K^2\geq 2\chi-6$.
\end{proof}

\begin{corollary}\label{cor:3Zappa}
If $\X$ is a good planar Zappatic degeneration of an $m$-canonical
surface of general type $\X_t$ to $\X_0 = X$ with at most $R_3$-,
$E_3$-, $E_4$- and $E_5$-points, then
\begin{itemize}
\item[(i)] $m\leq 6$; \item[(ii)] if $m=5,6$, then $\chi=3$,
$K^2=1$; \item[(iii)] if $m=4$, then $\chi\leq 4$, $8\chi\geq
11K^2+2$; \item[(iv)] if $m=3$, then $\chi\leq 6$, $8\chi\geq
7K^2+2$; \item[(v)] if $m=2$, then $K^2\leq 2\chi-1$; \item[(vi)]
if $m=1$, then $K^2\leq 4\chi-1$.
\end{itemize}
\end{corollary}

\begin{proof}
Take $\X_t =S$ to be $m$-canonical. First of all, by Corollary
\ref{cor:2Zappa}, we immediately get (i). Then, by Formula
\eqref{eq:1ZappaTer}, we get
$$8 \chi -2 \geq \frac{(m^2+m+2)}{2} K^2.$$
Thus, if $m$ equals either $1$ or $2$, we find statements $(v)$
and $(vi)$.

Since $S$ is of general type, by Noether's inequality we get
$$8\chi-2 \geq (2 \chi - 6) \frac{(m^2+m+2)}{2}.$$
This gives, for $m\ge3$,
$$\chi \leq 3 +  \frac{22}{(m^2+m-6)}$$
which, together with the above inequality, gives the other cases
of the statement.
\end{proof}

It would be interesting to see whether the numerical cases listed
in the above corollary can actually occur.

We remark that Corollary \ref{cor:1Zappabis} implies in particular
that one cannot hope to degenerate all surfaces to unions of
planes with only global normal crossings, namely double lines and
$E_3$-points; indeed, one needs at least $E_n$-points, for
$n\ge6$, or $R_m$-, $S_m$-points, for $m\ge4$.

From this point of view, another important result of Zappa is the
following (cf. \cite[\S\;12]{Za1}):

\begin{theorem}[Zappa] \label{thm:2Zappa}
For every $g\geq 2$ and every $d\geq 3g+2$, there are families of scrolls of sectional
genus $g$, of degree $d$, with general moduli having a planar Zappatic
degeneration with at most $R_3$-, $S_4$- and $E_3$-points.
\end{theorem}

Zappa's arguments rely on a rather intricate analysis concerning degenerations of hyperplane sections of the scroll and, accordingly, of the branch curve of a general projection of the scroll to a plane. We have not been able to check all the details of this very clever argument. However, the idea which Zappa exploits, of degenerating the branch curve of a general projection to a plane, is a classical one which goes back to Enriques, Chisini, etc. and certainly deserves 
attention. In reading Zappa's paper \cite{Za1}, our attention has been attracted by other ingredients he uses to study 
the aforementioned degenerations, which look interesting on their own. Precisely, Zappa gives extendability conditions for a curve on a scroll which is not a cone. For an overview of these results with a modern approach, the reader is referred to \cite[\S\;6]{CCFMLincei}. In particular, we prove a slightly more general version of 
the following result of Zappa (cf. \cite[Theorem 6.8]{CCFMLincei}):

\begin{proposition}[Zappa] \label{prop:2Zappa}
Let $C\subset \P^2$ be a general element of the Severi variety
$V_{d,g}$ of irreducible curves of degree $d$ and geometric genus
$g \geq 2$, with $d\geq 2g+2$. Then $C$ is the plane section of a scroll
$S\subset \P^3$ which is not a cone.
\end{proposition}

Via completely different techniques, we prove in \cite{CCFMLincei} the
following result which generalizes Theorem \ref{thm:2Zappa} of Zappa:

\begin{theorem}\label{thm:introLincei} (cf. \cite[Theorem 1.2]{CCFMLincei}) Let $g \geq 0$ and either $d \geq 2$, if $g = 0$,  or
$d \geq  5$, if $g = 1$, or $d \geq 2g+4$, if $g \geq 2$. Then
there exists a unique irreducible component ${\mathcal H}_{d,g}$ of the
Hilbert scheme of scrolls of degree $d$ and sectional genus $g$ in
$\Pp^{d-2g+1}$, such that the general point of ${\mathcal H}_{d,g}$
represents a smooth scroll $S$ which is linearly normal and
moreover with $H^1(S,\Oc_S(1))=0$.

Furthermore,
\begin{itemize}
\item[(i)] ${\mathcal H}_{d,g}$ is generically reduced and $\dim ({\mathcal H}_{d,g})
= (d-2g+2)^2+7(g-1)$, \item[(ii)] ${\mathcal H}_{d,g}$ contains the Hilbert
point of a planar Zappatic surface having only either $d-2$
$R_3$-points, if $g= 0$, or $d-2g+2$ points of type $R_3$ and
$2g-2$ points of type $S_4$, if $g \geq 1$, as Zappatic
singularities, \item[(iii)] ${\mathcal H}_{d,g}$ dominates the moduli space
${\mathcal M}_g$ of smooth curves of genus $g$.
\end{itemize}
\end{theorem}

\noindent
For other results concerning the geometry of the general scroll parametrized 
by ${\mathcal H}_{d,g}$, cf. \cite{CCFMnonsp} and \cite{CCFMbrill}.

In a more general setting, it is a natural question to ask which Zappatic singularities are
needed in order to degenerate as many smooth,
projective surfaces as possible to good, planar Zappatic surface.
Note that there are some examples
(cf.\ \S \ref{S:7}) of smooth projective surfaces $S$
which certainly cannot be degenerated to Zappatic surfaces with
$E_n$-, $R_n$-, or $S_n$-points, unless $n$ is large enough.

However, given such a $S$,  the next result --- i.e.\ Proposition
\ref{prop:9chi-k2} --- suggests that there might be a birational
model of $S$ which can be Zappatically degenerated to a surface
with only $R_3$- and $E_n$-points, for $n \le 6$.

\begin{proposition}\label{prop:9chi-k2}
Let $\X \to \D$ be a good planar Zappatic degeneration and assume
that the central fibre $X$ has at most $R_3$- and $E_m$-points,
for $m\leq 6$. Then
$$K^2\leq 9\chi.$$
\end{proposition}

\begin{proof}
The bounds for $K^2$ in Theorem \ref{thm:k2Gmain} give $9 \chi -
K^2 = 9v - 9 e + \sum_{m=3}^6 9 f_m - K^2$. Therefore, we get:
\begin{equation}\label{eq:basforhelp}
2(9\chi - K^2) \ge 2 e + 6f_3 + 2f_4 - 2f_5 - 6 f_6 - 2 r_3
\end{equation}
If we plug \eqref{eq:corBI2} in \eqref{eq:basforhelp}, we get
$$2(9\chi - K^2) \geq (2e + 3f_3 - 4f_4 - 5 f_5 - 6 f_6 - 2r_3) +
(3 f_3 + 6 f_4 + 3 f_5 ),$$ where both summands on the
right-hand-side are non-negative.
\end{proof}

In other words, Proposition \ref{prop:9chi-k2} states that the
Miyaoka-Yau inequality holds for a smooth projective surface $S$
which can degenerate to a good planar Zappatic surface with at
most $R_3$- and $E_n$-points, $3\le n \le 6$.

Another interesting application of the Multiple Point Formula is
given by the following remark.

\begin{remark}\label{rem:class} Let $\X \to \D$ be a good, planar Zappatic degeneration.
Denote by $\delta$ the {\em class} of the general fibre $\X_t$ of
$\X$, $t \neq 0$. By definition, $\delta$ is the degree of the
dual variety of $\X_t$, $t \neq 0$. From Zeuthen-Segre
(cf.\ \cite{E} and \cite{Iv}) and Noether's Formula (cf.\ \cite{GH}, page
600), it follows that:
\begin{equation}\label{eq:salesian1}
\delta=\chi_{\rm top}+\deg(\X_t) +4(g-1)=(9\chi-K^2)+3f+e.
\end{equation}
Therefore, \eqref{eq:corBI2} implies that:
\[
\delta\geq 3f_3 + r_3 +\sum_{n \geq 4}(12-n)f_n+ \sum_{n \geq 4} (n-1) \rho_n - k.
\]
In particular, if X is assumed to have at most $R_3$- and
$E_3$-points, then \eqref{eq:salesian1} becomes
$$\delta = (2e + 3f_3 - 2 r_3) + (3 f_3 + r_3),$$
where the first summand in the right-hand side is non-negative by
the Multiple Point Formula; therefore, one gets
$$\delta \geq 3f_3 + r_3.$$
Zappa's original approach in \cite{Za1}, indeed, was to compute $\delta$ and
then to deduce Formula \eqref{eq:1Zappa} and Theorem
\ref{thm:1Zappa} from this.
\end{remark}

\section{Examples of degenerations of surfaces to Zappatic unions of planes}\label{S:7}

The aim of this section is to illustrate some interesting examples of
degenerations of smooth surfaces to good, planar Zappatic surfaces.
We also discuss some examples of non-smoothable Zappatic surfaces
and we pose open questions on the existence of degenerations to planar Zappatic surfaces
for some classes of surfaces.

\bigskip

\noindent
{\bf Product of curves (Zappa)} (cf.\ Example 3.1 in \cite{CCFMLincei}).
Let $C\subset \Pp^{n-1}$ and $C'\subset \Pp^{m-1}$ be general enough curves.
Consider the smooth surface
$$S= C \times C'\subset \Pp^{n-1}\times \Pp^{m-1} \subset \Pp^{nm-1}.$$
If $C$ and $C'$ can degenerate to stick curves, say to $C_0$ and $C'_0$
respectively,
then the surface $S$ degenerates to a union of quadrics $Y$
with only double lines as singularities in codimension one
and with Zappatic singularities.

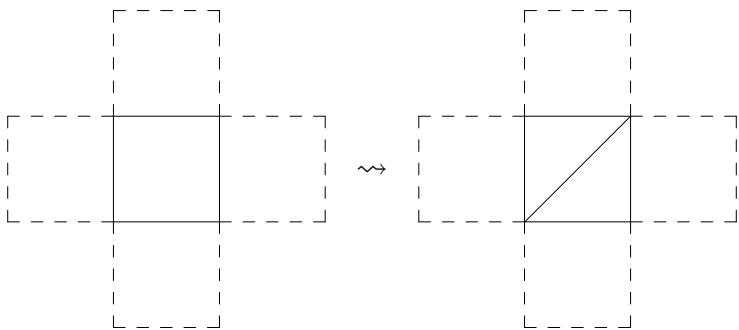
\begin{figure}[ht]
\[
\begin{xy}
0; <40pt,0pt>:
(1,1);(0,1)**@{--};(0,2)**@{--};(1,2)**@{--};(1,3)**@{--};(2,3)**@{--};
(2,2)**@{--};(3,2)**@{--};(3,1)**@{--};(2,1)**@{--};(2,0)**@{--};
(1,0)**@{--};(1,1)**@{--};
(2,1)**@{-};(2,2)**@{-};(1,2)**@{-};(1,1)**@{-} ,
\end{xy}
\quad
\raisebox{60pt}{$\rightsquigarrow$}
\quad
\begin{xy}
0; <40pt,0pt>:
(1,1);(0,1)**@{--};(0,2)**@{--};(1,2)**@{--};(1,3)**@{--};(2,3)**@{--};
(2,2)**@{--};(3,2)**@{--};(3,1)**@{--};(2,1)**@{--};(2,0)**@{--};
(1,0)**@{--};(1,1)**@{--};
(2,1)**@{-};(2,2)**@{-};(1,2)**@{-};(1,1)**@{-} ;(2,2)**@{-} ,
\end{xy}
\]
\caption{A quadric degenerating to the union of two planes}\label{fig:quadrics}
\end{figure}

If we could (independently) further
degenerate each quadric of $Y$ to the union of two planes,
then we would be able to get a degeneration $\X\to\Delta$
with general fibre $\X_t\cong S=C \times C'$ and central fibre a good, planar Zappatic
surface $\X_0\cong X$.
This certainly happens if each quadric of $Y$ meets the other quadrics of $Y$
along a union of lines of type $(a,b)$, with $a,b\leq 2$ (see Figure
\ref{fig:quadrics}).

For example $S =C\times C'$ can be degenerated to a good, planar Zappatic surface
if $C$ and $C'$ are either rational or elliptic normal curves and
we degenerate them to stick curves $C_{R_n}$ and $C_{E_n}$,
respectively.

Let us see some of these cases in detail.

\begin{itemize}

\item[(a)] \emph{Rational normal scrolls} (cf.\ Example 3.2 in \cite{CCFMLincei}).
Let $C\subset \Pp^{n}$ be a rational normal curve of degree $n$.
Since $C$ can degenerate to a stick curve $C_{R_n}$,
the surface $S= C \times \Pp^1 \subset \Pp^{2n+1}$ can degenerate
to a chain of $n$ quadrics as in Figure \ref{fig:nquadrics}, which is a good Zappatic
surface.

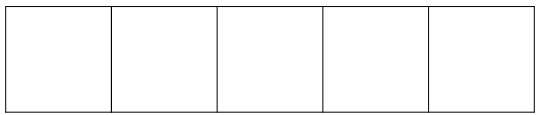
\begin{figure}[ht]
\[
\begin{xy}
0; <40pt,0pt>:
(0,0);(5,0)**@{-};(5,1)**@{-};(0,1)**@{-};(0,0)**@{-} ,
(1,0);(1,1)**@{-} ,
(2,0);(2,1)**@{-} ,
(3,0);(3,1)**@{-} ,
(4,0);(4,1)**@{-} ,
(1,0);(1,1)**@{-} ,
\end{xy}
\]
\caption{Union of $n$ quadrics with a chain as associated graph}
\label{fig:nquadrics}
\end{figure}

As we remarked above, the chain of quadrics can further degenerate
to a good, planar Zappatic surface $X$ which is the union of $2n$-planes
(see Figure \ref{fig:2nplanes}).
Note that the surface $X$ has only $R_3$-points as Zappatic
singularities and its associated graph $G_X$ is a chain $R_{2n}$.
In this way, one gets degeneration
of a rational normal scroll to a good, planar Zappatic surface
with only $R_3$-points as Zappatic singularities.

\begin{figure}[ht]
\[
\begin{xy}
0; <40pt,0pt>:
(0,0);(5,0)**@{-};(5,1)**@{-};(0,1)**@{-};(0,0)**@{-} ,
(1,0);(1,1)**@{-} ,
(2,0);(2,1)**@{-} ,
(3,0);(3,1)**@{-} ,
(4,0);(4,1)**@{-} ,
(1,0);(1,1)**@{-} ,
(0,0);(1,1)**@{-} ,
(1,0);(2,1)**@{-} ,
(2,0);(3,1)**@{-} ,
(3,0);(4,1)**@{-} ,
(4,0);(5,1)**@{-} ,
\end{xy}
\qquad
\raisebox{20pt}{$ %
\begin{xy}
0; <20pt,0pt>: 
(0,0)*=0{\bullet};(1,0)*=0{\bullet}**@{-};
(2,0)*=0{\bullet}**@{-};(3,0)*=0{\bullet}**@{-};
(4,0)*=0{\bullet}**@{-};(5,0)*=0{\bullet}**@{-};(6,0)*=0{\bullet}**@{-};
(7,0)*=0{\bullet}**@{-};(8,0)*=0{\bullet}**@{-};(9,0)*=0{\bullet}**@{-};
,
\end{xy}
$}
\]
\caption{Union of $2n$-planes with a chain as associated graph}
\label{fig:2nplanes}
\end{figure}
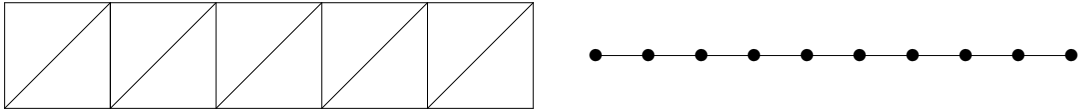

\item[(b)] \emph{Ruled surfaces} (cf.\ Example 3.3. in \cite{CCFMLincei}).
Let now $C \subset \Pp^{n-1}$ be a normal elliptic curve of degree $n$.
Since $C$ can degenerate to a stick curve $C_{E_n}$,
the surface $S=C \times \Pp^1 \subset \Pp^{2n-1}$ can degenerate
to a cycle of $n$ quadrics (see Figure \ref{fig:e2n}).

\begin{figure}[ht]
\[
\begin{xy}
0; <40pt,0pt>:
(0,1); {\xypolygon6"B"{~>{}}} ,
(0,0); {\xypolygon6"A"{~>{}}} ,
"A1";"A2"**@{-};"A3"**@{-};"A4"**@{-};"A5"**@{-};"A6"**@{-};"A1"**@{-} ,
"B1";"B6"**@{-};"B5"**@{-};"B4"**@{-} ,
"A1";"B1"**@{-} ,
"A2";"B2"**@{-} ,
"A3";"B3"**@{-} ,
"A4";"B4"**@{-} ,
"A5";"B5"**@{-} ,
"A6";"B6"**@{-} ,
"B1";"B2"**@{}?!{"A1";"A6"}="B12" ,
"B3";"B4"**@{}?!{"A4";"A5"}="B34" ,
"B1";"B12"**@{.};"B2"**@{-};"B3"**@{-};"B34"**@{-};"B4"**@{.} ,
\end{xy}
\]
\caption{A cycle of $n$ quadrics}\label{fig:e2n}
\end{figure}
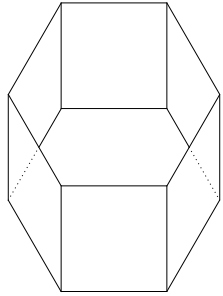

As before, such a cycle of quadrics can degenerate
to a good, planar Zappatic surface $X$ which is a union of $2n$ planes
with only $R_3$-points and whose associated graph corresponds to the cycle graph
$E_{2n+2}$ (see figure \ref{fig:e2n+2}).

\begin{figure}[ht]
\[
\begin{xy}
0; <40pt,0pt>:
(0,1); {\xypolygon6"B"{~>{}}} ,
(0,0); {\xypolygon6"A"{~>{}}} ,
"A1";"A2"**@{-};"A3"**@{-};"A4"**@{-};"A5"**@{-};"A6"**@{-};"A1"**@{-} ,
"B1";"B6"**@{-};"B5"**@{-};"B4"**@{-} ,
"A1";"B1"**@{-} ,
"A2";"B2"**@{-} ,
"A3";"B3"**@{-} ,
"A4";"B4"**@{-} ,
"A5";"B5"**@{-} ,
"A6";"B6"**@{-} ,
"B1";"B2"**@{}?!{"A1";"A6"}="B12" ,
"B3";"B4"**@{}?!{"A4";"A5"}="B34" ,
"B1";"B12"**@{.};"B2"**@{-};"B3"**@{-};"B34"**@{-};"B4"**@{.} ,
"B1";"A2"**@{}?!{"A1";"A6"}="AB1" ,
"B1";"AB1"**@{.};"A2"**@{-} ,
"B2";"A3"**@{-} ,
"B3";"A4"**@{-} ,
"B4";"A5"**@{-} ,
"B5";"A6"**@{-} ,
"B6";"A1"**@{-} ,
\end{xy}
\qquad
\qquad
\qquad
\raisebox{0pt}{$ %
\begin{xy}
0; <50pt,0pt>: 
(0,0.5) , {\xypolygon12{\bullet}}
\end{xy}
$}
\]
\caption{A union of $2n$ planes with a cycle as associated
graph}\label{fig:e2n+2}
\end{figure}
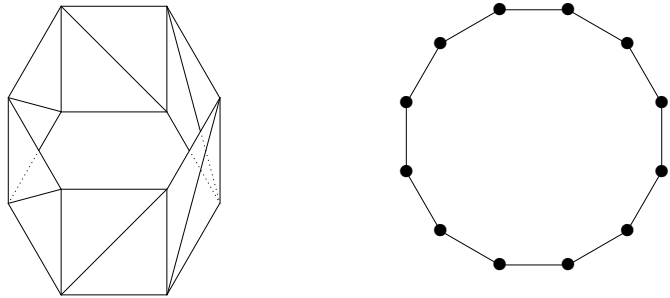

\item[(c)] \emph{Abelian surfaces} (cf.\ Example 3.4 in \cite{CCFMLincei}).
Let $C \subset \Pp^{n-1}$ and $C'\subset \Pp^{m-1}$ be
smooth, elliptic normal curves of degree respectively $n$ and $m$.
Then $C$ and $C'$ degenerate to the stick curves $C_{E_n}$ and $C_{E_m}$ respectively,
hence the abelian surface $S = C \times C'\subset \Pp^{nm-1}$ degenerates to
a Zappatic surface which is a union of $mn$ quadrics with only $E_4$-points as Zappatic
singularities, cf.\ e.g.\ the picture on the left in Figure \ref{fig:quadrics2},
where the top edges have to be identified with the bottom ones,
similarly the left edges have to be identified with the right ones.
Thus the top quadrics meet the bottom quadrics and
the quadrics on the left meet the quadrics on the right.

\begin{figure}[ht]
\[
\begin{xy}
0; <2em,0em>: (0,0);(5,0)**@{-} , (0,4);(5,4)**@{-} ,
(0,0);(0,4)**@{-} , (5,0);(5,4)**@{-} , (0,1);(5,1)**@{-} ,
(0,2);(5,2)**@{-} , (0,3);(5,3)**@{-} , (1,0);(1,4)**@{-} ,
(2,0);(2,4)**@{-} , (3,0);(3,4)**@{-} , (4,0);(4,4)**@{-} ,
(0,0);(1,0)**@{}, (0,4);(1,4)**@{}, (1,0);(2,0)**@{},
(1,4);(2,4)**@{}, (2,0);(3,0)**@{}, (2,4);(3,4)**@{},
(3,0);(4,0)**@{}, (3,4);(4,4)**@{}, (4,0);(5,0)**@{},
(4,4);(5,4)**@{}, (0,0);(0,1)**@{}, (5,0);(5,1)**@{},
(0,1);(0,2)**@{}, (5,1);(5,2)**@{}, (0,2);(0,3)**@{},
(5,2);(5,3)**@{}, (0,3);(0,4)**@{}, (5,3);(5,4)**@{},
\end{xy}
\qquad
\raisebox{4em}{$ \rightsquigarrow $}
\qquad
\begin{xy}
0; <2em,0em>:
(4,0);(5,1)**@{-} ,
(3,0);(5,2)**@{-} ,
(2,0);(5,3)**@{-} ,
(1,0);(5,4)**@{-} ,
(0,0);(4,4)**@{-} ,
(0,1);(3,4)**@{-} ,
(0,2);(2,4)**@{-} ,
(0,3);(1,4)**@{-} ,
(0,0);(5,0)**@{-} ,
(0,4);(5,4)**@{-} ,
(0,0);(0,4)**@{-} ,
(5,0);(5,4)**@{-} ,
(0,1);(5,1)**@{-} ,
(0,2);(5,2)**@{-} ,
(0,3);(5,3)**@{-} ,
(1,0);(1,4)**@{-} ,
(2,0);(2,4)**@{-} ,
(3,0);(3,4)**@{-} ,
(4,0);(4,4)**@{-} ,
(0,0);(1,0)**@{},
(0,4);(1,4)**@{},
(1,0);(2,0)**@{},
(1,4);(2,4)**@{},
(2,0);(3,0)**@{},
(2,4);(3,4)**@{},
(3,0);(4,0)**@{},
(3,4);(4,4)**@{},
(4,0);(5,0)**@{},
(4,4);(5,4)**@{},
(0,0);(0,1)**@{},
(5,0);(5,1)**@{},
(0,1);(0,2)**@{},
(5,1);(5,2)**@{},
(0,2);(0,3)**@{},
(5,2);(5,3)**@{},
(0,3);(0,4)**@{},
(5,3);(5,4)**@{},
\end{xy}
\]
\caption{$nm$ quadrics with $E_4$-points and $2nm$ planes with $E_6$-points}\label{fig:quadrics2}
\end{figure}
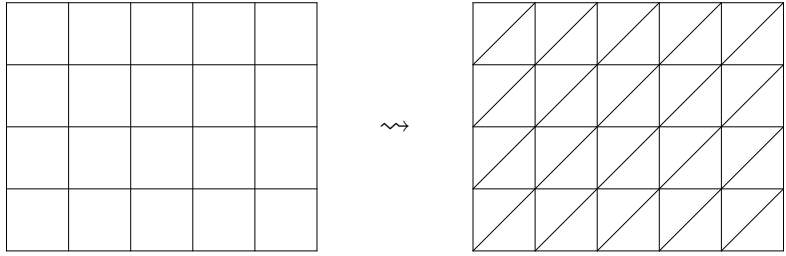

Again each quadric degenerates to the union of two planes. By
doing this as depicted in Figure \ref{fig:quadrics2}, one gets a
degeneration of a general abelian surface with a polarization of
type $(n, m)$ to a planar Zappatic surface of degree $2nm$ with
only $E_6$-points as Zappatic singularities.

\end{itemize}

Concerning the general case, suppose that either $C$ or $C'$ has genus greater
than 1. If $C$ and $C'$ degenerate to stick curves, then the surface $S=C\times C'$
degenerates as above to a union of quadrics.
Unfortunately it is not clear if it is possible to further
independently degenerate each quadric to two planes in the same way as above
and, in \cite{Za2}, Zappa left as an open problem to prove the degeneration to a union of planes.

Theorem \ref{thm:introLincei}, proved in the paper \cite{CCFMLincei},
show that $C \times \mathbb{P}^1$, suitably embedded as a non-special, linearly normal scroll,
really degenerates to a union of planes with only $R_3$- and $S_4$-points.
Indeed, let $C$ be any curve of genus $g$ and let $L$ be a very-ample non-special line bundle of degree
$d \geq g +3$. The global sections of $L$ determine an embedding of $C$ in $\mathbb{P}^{d-g}$. Consider
the Segre embedding of $C\times \mathbb{P}^1$. This gives a non-special, linearly normal, smooth
scroll $S$ of degree $2d$ in $\mathbb{P}^{2d-2g+1}$ and the corresponding point sits in
the irreducible component ${\mathcal H}_{2d,g}$ of the Hilbert scheme mentioned in Theorem \ref{thm:introLincei}.
By our results in \cite{CCFMLincei}, it follows that $S$ degenerates
to a planar Zappatic surface with $2(d-g+1)$ points of type $R_3$ and $2(g-1)$ points of type $S_4$
(cf.\ Construction 4.2 in \cite{CCFMLincei}).

\bigskip

\noindent
{\bf Degeneration to cones.} Recall the following result of Pinkham:

\begin{theorem}[Pinkham \cite{Pink1,Pink2}] \label{thm:pinkam}
Let $S \subset \Pp^n$ be a smooth, irreducible and projectively Cohen-Macaulay
surface.
Then $S$ degenerates to the cone over a hyperplane section of $S$.
\end{theorem}

Let $C$ be the hyperplane section of $S$.
Suppose that $C$ can be degenerated to a stick curve $C_0$.
In this case, $S$ can be degenerated to the cone $X$ over the stick curve $C_0$.
By definition, $X$ is a Zappatic surface only if $C$ has genus either $0$ or
$1$.
Therefore:

\begin{corollary} 
(i) Any surface $S$ of minimal degree (i.e.\ of degree $n$) in $\Pp^{n+1}$ can
be degenerated to the cone
over the stick curve $C_{T_n}$,
for any tree $T_n$ with $n$ vertices  (cf.\ Example \ref{ex:tngraphs}).

\noindent (ii) Any del Pezzo surface $S$ of degree $n$ in $\Pp^{n}$, $n\leq 9$,
can be degenerated to the cone over the stick curve $C_{Z_n}$,
for any connected graph  with $n\geq 3$ vertices and $h^1(Z_n, \CC)=1$
(cf.\ Example \ref{ex:zngraphs}).
\end{corollary}

For $n=4$, recall that the surfaces of minimal degree in $\Pp^5$
are either the Veronese surface (which has $K^2=9$)
or a rational normal scroll (which has $K^2=8$).
Therefore:

\begin{corollary}[Pinkham] 
The local deformation space of a $T_4$-singularity is reducible.
\end{corollary}

\bigskip

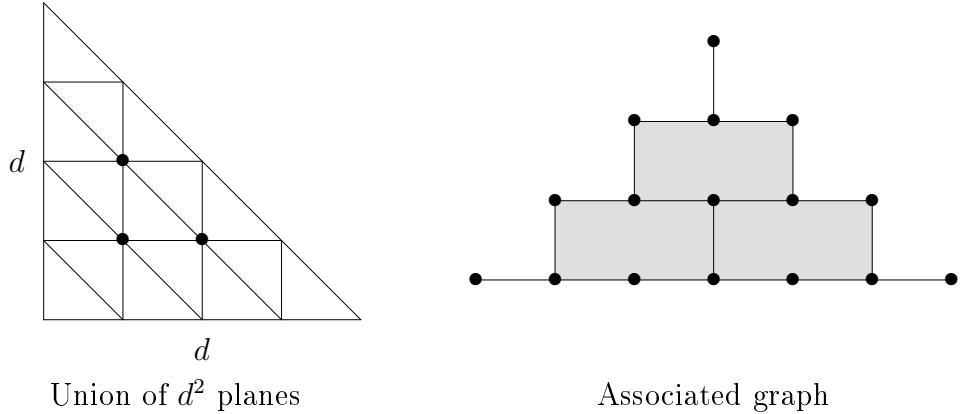
\begin{figure}[ht]
\[
\begin{array}{ccc}
\begin{xy}
0; <30pt,0pt>:
(0,0);(4,0)**@{-};(0,4)**@{-};(0,0)**@{-} ,
(0,1);(3,1)**@{-} ,
(0,2);(2,2)**@{-} ,
(0,3);(1,3)**@{-} ,
(1,0);(1,3)**@{-} ,
(2,0);(2,2)**@{-} ,
(3,0);(3,1)**@{-} ,
(1,0);(0,1)**@{-} ,
(2,0);(0,2)**@{-} ,
(3,0);(0,3)**@{-} ,
(1,1)*=0{\bullet} ,
(2,1)*=0{\bullet} ,
(1,2)*=0{\bullet} ,
(2,0)*++!U{ d} ,
(0,2)*++!R{ d} ,
\end{xy}
& \qquad &
\raisebox{15pt}{$
\begin{xy}
0; <60pt,0pt>:
(0.5,0.5).(1.5,1)*[!\colorxy(0.875 0.875 0.875)]\frm{*} ,
(1.5,0.5).(2.5,1)*[!\colorxy(0.875 0.875 0.875)]\frm{*} ,
(1,1).(2,1.5)*[!\colorxy(0.875 0.875 0.875)]\frm{*} ,
(0,0.5);(3,0.5)**@{-} ,
(0.5,1);(2.5,1)**@{-} ,
(1,1.5);(2,1.5)**@{-} ,
(0.5,0.5);(0.5,1)**@{-} ,
(1,1);(1,1.5)**@{-} ,
(1.5,1.5);(1.5,2)**@{-} ,
(1.5,0.5);(1.5,1)**@{-} ,
(2,1);(2,1.5)**@{-} ,
(2.5,0.5);(2.5,1)**@{-} ,
(0,0.5)*=0{\bullet} ,
(0.5,0.5)*=0{ \bullet} ,
(1,0.5)*=0{\bullet} ,
(1.5,0.5)*=0{ \bullet} ,
(2,0.5)*=0{\bullet} ,
(2.5,0.5)*=0{ \bullet} ,
(3,0.5)*=0{ \bullet} ,
(0.5,1)*=0{ \bullet} ,
(1,1)*=0{ \bullet} ,
(1.5,1)*=0{ \bullet} ,
(2,1)*=0{ \bullet} ,
(2.5,1)*=0{ \bullet} ,
(1,1.5)*=0{ \bullet} ,
(1.5,1.5)*=0{ \bullet} ,
(2,1.5)*=0{ \bullet} ,
(1.5,2)*=0{ \bullet} ,
\end{xy}
$}
\\[1mm]
\text{Union of $d^2$ planes} & & \text{Associated graph}
\end{array}
\]
\caption{Degeneration of the $d$-Veronese surface.}\label{fig:vero}
\end{figure}

\noindent
{\bf Veronese surfaces (Moishezon-Teicher).}
Consider $V_d \subset \Pp^{d(d+3)/2}$ be the $d$-Veronese surface,
namely the embedding of $\Pp^2$ via the linear system $|\Oc_{\Pp^2}(d)|$.
In \cite{MoiTei}, Moishezon and Teicher described a ``triangular'' degeneration
of $V_d$
such that the central fibre is a union of $d^2$ planes
with only $R_3$- and $E_6$-points as Zappatic singularities
(see Figure \ref{fig:vero}).

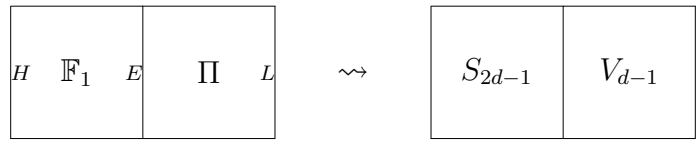
\begin{figure}[ht]
\[
\begin{xy}
0; <50pt,0pt>:
(0,0)="a";(1,0)="b"**@{-};(2,0)="c"**@{-};(2,1)="f"**@{-};
(1,1)="d"**@{-};(0,1)="e"**@{-};"a"**@{-} ,
"b";"d"**@{-}?*!R+{\st E} ,
"a"."d"!C*{\Ff_1} ,
"c"."d"!C*{\Pi} ,
"a";"e"**@{}?*!L+{\st H} ,
"c";"f"**@{}?*!R+{\st L} ,
\end{xy}
\qquad
\raisebox{25pt}{$\rightsquigarrow$}
\qquad
\begin{xy}
0; <50pt,0pt>:
(0,0)="a";(1,0)="b"**@{-};(2,0)="c"**@{-};(2,1)="f"**@{-};
(1,1)="d"**@{-};(0,1)="e"**@{-};"a"**@{-} ,
"b";"d"**@{-} ,
"a"."d"!C*{S_{2d-1}} ,
"c"."d"!C*{V_{d-1}} ,
\end{xy}
\]
\caption{Degenerating the $d$-Veronese surface $V_d$}\label{fig:proof}
\end{figure}

Let us explain how to get such a degeneration.
Consider the trivial family $\X = \Pp^2\times\Delta$, where $\D$ is a complex
disk.
Let $\mathcal{L} \cong \Oc_{\X}(d)$ be a line bundle on $\X$.
If we blow-up a point in the central fibre of $\X$, the new central
fibre becomes as the left-hand-side picture in Figure \ref{fig:proof},
where $E$ is the exceptional divisor, $\Pi = \Pp^2$ and $\Ff_1$ is the
Hirzebruch surface.

Let $\tilde{\X}$ be the blown-up family and let
$\tilde{\mathcal{L}}$ be
the line bundle on $\tilde{\X}$ given by
the pull-back of $\mathcal{L}$ twisted
by the divisor $- (d-1) \Pi$. Then $\tilde{\mathcal{L}}$ restricts to
$d H - (d-1) E$ on the surface $\Ff_1$ and to
$(d-1) L$ on the plane $\Pi$, respectively.

These line bundles embed, respectively,
$\Ff_1$ as a rational normal scroll $S_{2d-1}$
of degree $2d-1$ and $\Pi$ as the $(d-1)$-Veronese
surface $V_{d-1}$ meeting along a rational normal
curve of degree $d-1$ (see the right-hand-side picture
of Figure \ref{fig:proof}). One can independently degenerate $V_{d-1}$ ro a good, planar Zappatic surface, by
induction on $d$ (getting the left most bottom triangle in figure \ref{fig:vero})
and the rational normal scroll $S_{2d-1}$ as we saw before (getting the top strip in the triangle
in figure \ref{fig:vero}).

\bigskip

\noindent
{\bf K3 surfaces.} In the paper \cite{CMT}, the authors construct a specific projective
degeneration of a $K3$ surface of degree $2g-2$ in
$\Pp^g$ to a planar Zappatic surface which is a union of
$2g-2$ planes, which meet in such a way that the associated graph
to the configuration of planes is a triangulation of the $2$-sphere.

In the previous paper \cite{CLM}, planar Zappatic degenerations of $K3$ surfaces
were constructed in such a way that the general member of the degeneration was
embedded by a primitive line bundle.
In \cite{CMT} the general member of the degeneration is
embedded by a multiple of a primitive line bundle class (for details, the
reader is referred to the original articles).

Let $X$ denote a good, planar Zappatic surface which is a degeneration of a $K3$-surface
$S \subset \Pp^g$ of genus $g$ and let $G$ be the associated graph to $X$.
Then, $G$ is planar, because $p_g(S)=1$, and $3$-valent (see \cite{CLM}). By
using
Notation \ref{def:vefg}, we get
\begin{equation}\label{eq:k3}
v=2g-2, \qquad e=3g-3, \qquad f=g+1.
\end{equation}

Conversely, by starting from a planar graph $G$ with invariants as in
\eqref{eq:k3},
one can find a Zappatic numerical $K3$ surface $X$ whose associated graph is
$G$.
Such an $X$ is called a {\em graph surface}. Smoothable graph surfaces are
exhibited in \cite{CLM} and \cite{CMT}.

The specific degenerations constructed in \cite{CMT} depend on
two parameters and can be viewed as two rectangular
arrays of planes, joined along their boundary. For this reason, these are called
{\em pillow degenerations}.

Take two integers $a$ and $b$ greater than or equal to two and set
$g=2ab+1$. The number of planes in the pillow degeneration is then
$2g-2 = 4ab$. The projective space $\Pp^g$
has $g+1 = 2ab+2$ coordinate points,
and each of the $4ab$ planes is obtained as the span of three
of these.  These sets of three points are indicated in Figure \ref{fig:pillow},
which describes the bottom part of the ``pillow''
and the top part of the ``pillow'',
which are identified along the boundaries of the two
configurations.
The reader will see that the boundary is a cycle of $2a+2b$ lines.

\begin{figure}[ht]
\caption{Configuration of Planes, Top and Bottom}\label{fig:pillow}
\centerline{\begin{picture}(340,175)
\put(70,163){Top}
\put(30,150){\scriptsize Boundary points labeled from}
\put(30,140){\scriptsize $1$ through $2a+2b$, clockwise;}
\put(30,130){\scriptsize interior points labeled from}
\put(30,120){\scriptsize $2a+2b+1$ through $ab+a+b+1$}
\put(27,103){$\scriptstyle 1$}
\put(37,103){$\scriptstyle 2$}
\put(52,103){$\scriptstyle \cdots$}
\put(117,103){$\scriptstyle a$}
\put(132,103){$\scriptstyle a+1$}
\put(0,90){$\scriptstyle 2a+2b$}\put(132,90){$\scriptstyle a+2$}
\put(132,80){$\scriptstyle a+3$}
\put(10,60){$\scriptstyle \vdots$}\put(132,60){$\scriptstyle \vdots$}
\put(0,15){$\scriptstyle a+2b+1$}\put(132,15){$\scriptstyle a+b+1$}
\multiput(30,20)(10,0){11}{\line(0,1){80}}
\multiput(30,20)(0,10){9}{\line(1,0){100}}
\multiput(30,20)(0,10){8}{\multiput(0,0)(10,0){10}{\line(1,1){10}}}
%
\put(245,163){Bottom}
\put(210,150){\scriptsize Boundary points labeled from}
\put(210,140){\scriptsize $1$ through $2a+2b$, clockwise;}
\put(210,130){\scriptsize interior points labeled from}
\put(210,120){\scriptsize $ab+a+b+2$ through $2ab+2$}
\put(207,103){$\scriptstyle 1$}
\put(217,103){$\scriptstyle 2$}
\put(232,103){$\scriptstyle \cdots$}
\put(297,103){$\scriptstyle a$}
\put(312,103){$\scriptstyle a+1$}
\put(180,90){$\scriptstyle 2a+2b$}\put(312,90){$\scriptstyle a+2$}
\put(312,80){$\scriptstyle a+3$}
\put(190,60){$\scriptstyle \vdots$}\put(312,60){$\scriptstyle \vdots$}
\put(180,15){$\scriptstyle a+2b+1$}\put(312,15){$\scriptstyle a+b+1$}
\multiput(210,20)(10,0){11}{\line(0,1){80}}
\multiput(210,20)(0,10){9}{\line(1,0){100}}
\multiput(210,20)(0,10){8}{\multiput(0,0)(10,0){10}{\line(1,1){10}}}
\end{picture}}
\end{figure}

Note that no three of the planes meet along a line.
Also note that the set of bottom planes
lies in a projective space of dimension $ab+a+b$,
as does the set of top planes;
these two projective spaces meet
exactly along the span of the $2a+2b$ boundary points,
which has dimension $2a+2b-1$.
Finally note that the four corner points of the pillow
degeneration
(labelled $1$, $a+1$, $a+b+1$, and $a+2b+1$)
are each contained in three distinct planes,
whereas all the other points are each contained in six planes.
This property, that the number of lines and planes
incident on each of the points is bounded
is a feature of the pillow degeneration
that is not available in other previous degenerations (see \cite{CLM}).
We will call such a configuration of planes a
\emph{pillow of bidegree $(a,b)$}.

Observe that a pillow of bidegree $(a,b)$ is a planar Zappatic surface
of degree $2g-2$, having four $R_3$-points and
$2ab-2=g-3$ $E_6$-points as Zappatic singularities.

\begin{remark}
$E_n$-points, with $n\geq 6$, are
unavoidable for the degeneration of $K3$ surfaces with hyperplane
sections of genus $g$, if $g\geq 12$.
Indeed, by using Notation \ref{def:vefg}, Formula \eqref{eq:k3} and
the fact that $G$ is $3$-valent, we get
$$g+1 = f = \sum_n f_n \;\; {\rm and} \;\; 6(g-1) = 3v = \sum_n n f_n.$$
These give
$$\sum_n(6-n) f_n = 12.$$
If we assume that $f_n = 0$, for $n \geq 6$, the last equality
gives $2f_4 + f_5 = 12$ and so $10 f_4 + 5 f_5 = 60$.
This equality, together with
$4f_4 + 5f_5 = 6g-6$ gives $g+f_4 = 11$, i.e.\ $g \le 11$.
\end{remark}

\bigskip

\noindent
{\bf Complete intersections.} Consider a surface $S \subset \Pp^n$ which is
a \emph{general} complete intersection of type $(d_1, \ldots, d_{n-2})$.
Namely $S$ is defined as the zero-locus
$$f_1= \cdots =f_{n-2}=0,$$
where $f_i$ is a
general homogeneous polynomial of degree $d_i$, $ 1 \le i \le n-2$.

One can degenerate any hypersurface $f_i = 0$ to the union of $d_i$ hyperplanes.
This implies that $S$ degenerates
to a planar Zappatic surface $X$ with global normal crossings,
i.e.\ with only $E_3$-points as Zappatic singularities.

We remark that  degenerations of surfaces to good, planar Zappatic ones are possible also
when $S$ is projectively Cohen-Macaulay in $\Pp^4$,
by some results of Gaeta (see \cite{Gaeta}).

\bigskip

\begin{figure}[ht]
\[
\begin{xy}
0; <50pt,0pt>:
(0,0)="a"*=0{\bullet};(0,1)*=0{\bullet}**@{-};(1,1)="b"*=0{\bullet}**@{-};
(1,0)*=0{\bullet}**@{-};"a"**@{--} ,
"b";(2,1)*=0{\bullet}**@{-} ,
\end{xy}
\]
\caption{A non-smoothable planar Zappatic surface}\label{fig:nonsmooth}
\end{figure}
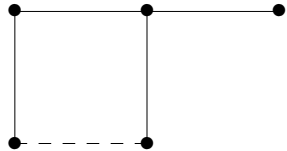

\noindent
{\bf Non-smoothable Zappatic surfaces.}
The results of the previous sections allow us to exhibit simple planar Zappatic
surfaces
which are not smoothable, i.e.\ which cannot be the central fibre of
a degeneration.
For example, the planar Zappatic surface $X$ with the graph of Figure
\ref{fig:nonsmooth}
as associated graph is not smoothable.
Indeed, if $X$ were the central fibre of a degeneration
$\X\to\Delta$, then Formulas \eqref{eq:K2Gbounds} and \eqref{eq:K2Gbounds2}  would
imply
\[
9\le K^2 \le 10,
\]
which is absurd because of the classification of
smooth projective surface of degree $5$ in $\Pp^6$ (see Theorem
\ref{thm:DelP1}).

It is an interesting problem to find more examples of smoothable
Zappatic surfaces with
only $R_3$- and $E_n$-points, $3\le n \le 6$.

E.g., does there exist a Zappatic degeneration with only $R_3$
and $E_n$, $3\le n \le 6$, for Enriques' surfaces?

\appendix

\section{Normal, Cohen-Macaulay and Gorenstein properties}\label{S:pCMpG}

The aim of this appendix is to briefly recall some well-known terminology
and results concerning normal, Cohen-Macaulay and Gorenstein properties
for projective varieties. These results are frequently used in Section \ref{S:MS}.
In the sequel, the term variety does not imply that the scheme under consideration is
supposed to be irreducible.

First, we focus on ``normality'' conditions.

\begin{definition}\label{def:1}
An algebraic scheme $X$ is said to be {\em normal} if, for each $p \in X$, the
local ring $\Oc_{X,p}$ is an integrally closed domain (cf.\ \cite{H}, pp. 23 and 91).
\end{definition}

Recall that if $A$ is a local ring, with maximal ideal $m$,
the {\em depth} of $A$ is the maximal
length of a regular sequence $x_1, \ldots, x_r$ with all
$x_i \in m$ (cf.\ \cite{H}, page 184).

Recall that, by the theorem of Krull-Serre (see, e.g., \cite{H}, Theorem 8.22A, page 185),
$X$ is normal at $p$ if, and only if, $X$ is non-singular in codimension one at $p$ and
the local ring $\Oc_{X,p}$ has the $S_2$-property, i.e.\ its depth is greater than
or equal to $2$.

Let $X \subset \PR$ be a projective variety and let
\begin{equation}\label{eq:ideal}
I(X):= H^0_*(\Ii_{X|\PR}) = \bigoplus_{n \in \ZZ} H^0(\Ii_{X|\PR} (n))
\end{equation}be the saturated ideal associated to $X$ in the homogeneous polynomial ring
$S: = \CC[x_0, \ldots, x_r]$.

\begin{definition}\label{def:2}
A projective variety $X \subset \PR$ is said to be
{\em projectively normal} (with respect
to the given embedding) if its homogeneous coordinate ring
\begin{equation}\label{eq:gamma}
\Gamma(X) :=  \CC[x_0, \ldots, x_r] /I(X)
\end{equation}is an integrally closed domain (cf.\ \cite{H}, pp. 23 and 126).
\end{definition}

Recall that projective normality is a property of the given embedding $X \subset \PR$ and
not only of $X$. Observe also that if $X$ is a projectively normal variety,
then it is also irreducible and normal (see \cite{H}, page 23).

If $X \subset \PR$ is a projective variety, we denote by
\begin{equation}\label{eq:cone}
{\mathcal C}(X) \subset \A^{r+1}
\end{equation}the affine cone over $X$ having vertex at the origin $\underline{o} \in \A^{r+1}$.
Observe that the homogeneous coordinate ring $\Gamma(X)$ in \eqref{eq:gamma}
coincides with the coordinate ring of the affine cone ${\mathcal C}(X)$.

We can characterize projective normality of a closed projective variety
$X \subset \PR$ in terms of its affine cone ${\mathcal C}(X)$.

\begin{proposition}\label{prop:2}
Let $X \subset \PR$ be a projective variety
and let ${\mathcal C}(X)$ be its
affine cone. Then,
the following conditions are equivalent:
\begin{itemize}
\item[(i)] $X$ is projectively normal;
\item[(ii)] ${\mathcal C}(X)$ is normal;
\item[(iii)] ${\mathcal C}(X)$ is normal at the vertex $\underline{o}$.
\end{itemize}
\end{proposition}
\begin{proof}
See \cite{Greco}, Prop.\ 7.10, p.\ 57, and
\cite{H}, p.\ 147.
\end{proof}

\begin{remark}\label{rem:0}
The above proposition is an instance of a general philosophy which states
that ``the properties
of the vertex of ${\mathcal C}(X)$ are equivalent to global properties
of ${\mathcal C}(X)$ as well as of $X$'' (cf.\ \cite{Greco}, page 54).
\end{remark}

\begin{proposition}\label{prop:2bis}
Let $X \subset \PR$ be a projective variety and let $H$ be a hyperplane section
of $X$. If $H$ is projectively normal, then $X$ is projectively normal.
\end{proposition}
\begin{proof} It is a trivial consequence of Krull-Serre's theorem
(see, e.g., \cite{Greco}, Theorem 4.27).
\end{proof}

\begin{definition}\label{def:3}
A projective variety $X \subset \PR$ is called \emph{arithmetically normal} if
the restriction map
\begin{equation}\label{eq:pn}
H^0(\Oc_{\PR}(j)) \to H^0(\Oc_X(j))
\end{equation} is surjective, for every $j \in \NN$.
\end{definition}

\begin{remark}\label{rem:A2}
The surjectivity of the map in
\eqref{eq:pn} is equivalent to
\begin{equation}\label{eq:pn2}
H^1(\Ii_{X|\PR}(j))=0, \qquad \text{for every $j\in\ZZ$.}
\end{equation}This follows from the sequence:
\begin{equation}\label{eq:*}
0 \to \Ii_{X|\PR}(j) \to \Oc_{\PR}(j) \to \Oc_X(j) \to 0
\end{equation}and from the cohomology of projective spaces.
\end{remark}

One has the following relationship among the above three notions:

\begin{proposition}\label{prop:1}
$X \subset \PR$ is projectively normal if, and only if, $X$ is normal and arithmetically
normal in $\PR$ (cf.\ \cite{H}, pg. 126, and \cite{Zar0}).
\end{proposition}

A fundamental property related to arithmetical normality is the following:

\begin{proposition}\label{lem:pn}
If a hyperplane section $H$ of a projective variety $X$ is
arithmetically normal, then $X$ is arithmetically normal.
\end{proposition}

\begin{proof}
Consider the following commutative diagram with two short exact
sequences of sheaves:
\begin{equation}
\label{eq:S-H}
\raisebox{4ex}{\xymatrix{%
 0 \ar[r] & \Oc_X(j-1) \ar[r] & \Oc_X(j) \ar[r] & \Oc_H(j) \ar[r] & 0 \\
 0 \ar[r] & \Oc_{\PR}(j-1) \ar[r]\ar[u] & \Oc_{\PR}(j) \ar[r]\ar[u] &
\Oc_{\Pp^{r-1}}(j) \ar[r]\ar[u] & 0
}}
\end{equation}
for all $j \in \ZZ$, where the vertical arrows are defined by the
usual short exact sequence \eqref{eq:*}. Diagram \eqref{eq:S-H}
induces in cohomology the following commutative diagram:
\[
\xymatrix{%
 H^0(\Oc_X(j)) \ar[r]^{\gamma_j} & H^0(\Oc_H(j)) \\
 H^0(\Oc_{\PR}(j)) \ar@{>>}[r]^-{\alpha_j}\ar[u]^{\delta_j} & H^0(\Oc_{\Pp^{r-
1}}(j)) \ar@{>>}[u]^{\beta_j}
}
\]
for every $j \in \ZZ$, where $\alpha_j$ is trivially surjective
and $\beta_j$ is surjective by hypothesis. The surjectivity of the
composite map $\gamma_j \circ \delta_j = \beta_j \circ \alpha_j$
forces the surjectivity of $\gamma_j$, but not yet that of
$\delta_j$. However, we have the injection:
\[
H^1(\Oc_X(j-1)) \hookrightarrow H^1(\Oc_X(j)),
\]
for every $j \in \ZZ$, which implies
\begin{equation}\label{eq:H1Ox}
H^1(\Oc_X(j))=0, \qquad\text{for each $j$,}
\end{equation}
because $H^1(\Oc_X(j))=0$ for $j \gg 0$ by Serre's Theorem. The
long exact sequences in cohomology induced by \eqref{eq:S-H} then
become:
\[
\xymatrix{%
 0 \ar[r] & H^0(\Oc_X(j-1)) \ar[r] & H^0(\Oc_X(j)) \ar[r] & H^0(\Oc_H(j)) \ar[r]
& 0 \\
 0 \ar[r] & H^0(\Oc_{\PR}(j-1)) \ar[r]\ar[u]^{\delta_{j-1}}
  & H^0(\Oc_{\PR}(j)) \ar[r]\ar[u]^{\delta_j} & H^0(\Oc_{\Pp^{r-1}}(j))
\ar[r]\ar@{>>}[u]^{\beta_j} & 0
}
\]
where we recall that $\beta_j$ is surjective for all $j$ by
hypothesis. The map $\delta_j$ is trivially surjective for $j \le
0$. Since $\delta_{j-1}=\delta_0$ is surjective for $j=1$, the map
$\delta_j=\delta_1$ is surjective too by elementary diagram chase.
Hence we conclude by induction on $j$.
\end{proof}

We consider now  ``Cohen-Macaulay'' conditions.

Recall that a local ring $(A, m)$
is said to be
{\em Cohen-Macaulay} ({\em CM} for short)
if $\depth(A) = \dim(A)$ (see \cite{H}, page 184).
This is equivalent to saying that a zero-dimensional local
ring is always CM and, if $\dim(A) >0$, then $A$ is CM if, and only if,
there is a non zero-divisor $x$ in $A$ such that $A/(x)$ is CM.
In this case, $A/(x)$ is CM for every non-zero divisor $x$ in $A$
(cf.\ \cite{Eisenbud}, \cite{Matsumura} page 107).

\begin{remark}\label{rem:freerescm}
Let $(R, M)$ be a regular local ring. Let $(A, m)$ be a local
$R$-algebra which is finitely generated as a $R$-module.
Set $c:= \codim_R(A)$. One has that $(A,m)$ is CM if, and only if,
there is a minimal, free resolution of $R$-modules
\begin{equation}\label{eq:dualres1}
{\mathcal F} \colon 0 \to F_c \to F_{c-1} \to \cdots \to F_1 \to R \to A \to 0,
\end{equation} (cf.\ Corollary 21.16  and
the Auslander-Buchsbaum Formula Theorem 19.9 in \cite{Eisenbud}).

If $(A,m)$ is CM, then one defines
\begin{equation}\label{eq:dual}
\omega_A:= Ext_R^c (A, R)
\end{equation}to be the {\em canonical module}
of $A$ (cf.\ Theorem 21.15 in \cite{Eisenbud}).
Then, the resolution \eqref{eq:dualres1} is such that
${\mathcal F}^*$ is a minimal free resolution of $\omega_A$ (cf.\ Corollary 21.16 \cite{Eisenbud}).
\end{remark}

A finitely generated $\CC$-algebra $B$ is said to be {\em Cohen-Macaulay} ({\em CM} for short)
if, and only if, for every prime ideal $p$ of $B$, $B_p$ is a CM local ring. This is
equivalent to saying that $B_p$ is a CM local ring, for every maximal
ideal $p$ in $B$ (cf.\ \cite{Eisenbud}, Prop. 18.8).

\begin{definition}\label{def:4} (cf.\ \cite{H}, page 185)
An algebraic, equidimensional scheme
$X$ is {\em Cohen-Macaulay at a point} $p \in X$ ({\em CM at} $p$, for short)
if $\Oc_{X,p}$ is a Cohen-Macaulay, local ring. X is {\em Cohen-Macaulay}
({\em CM}, for short) if
$X$ is Cohen-Macaulay at each $p \in X$.
\end{definition}

We have the following result:

\begin{theorem}\label{thm:1bis}
Let $X \subset \PR$ be an equidimensional projective variety and let $p $ be a point of $X$.
Then, the following
conditions are equivalent:
\begin{itemize}
\item[(i)] $X$ is CM at $p$;
\item[(ii)] each equidimensional hyperplane section of $X$ through $p$ is CM at $p$;
\end{itemize}
\end{theorem}
\begin{proof}
It directly follows from the definition of local, CM rings.
\end{proof}

\begin{definition}\label{def:5}
A projective, equidimensional variety $X \subset \PR$ is said to be
\emph{projectively Cohen-Macaulay} (\emph{pCM}, for short) if
the ring $\Gamma(X)$ as in \eqref{eq:gamma} is a Cohen-Macaulay ring.
\end{definition}

\begin{definition}\label{def:6}
A projective and equidimensional variety $X \subset \PR$ is said to be
\emph{arithmetically Cohen-Macaulay} (\emph{aCM}, for short), if $X$
is arithmetically normal and moreover
\begin{equation}\label{eq:defpCM}
H^i(\Oc_X(j))=0,\qquad\text{for every $j\in\ZZ$ and $1 \le i \le
n-1$,}
\end{equation}where $n = \dim(X)$.
\end{definition}

\begin{remark}\label{rem:73}
By standard exact sequences, $X$ is arithmetically Cohen-Macaulay
iff
\begin{equation}\label{eq:defpCM2}
H^i(\Ii_{X|\PR}(j))=0,\qquad\text{for every $j\in\ZZ$ and $1 \le i \le
n$.}
\end{equation}
\end{remark}

\begin{remark}\label{rem:curve}
For $n=1$, Formula \eqref{eq:defpCM} trivially holds. Thus a curve
is arithmetically Cohen-Macaulay if and only if it is arithmetically
normal.
\end{remark}

As in Proposition \ref{prop:2}, we have:

\begin{proposition}\label{prop:3}
Let $X \subset \PR$ be an equidimensional variety, ${\mathcal C}(X) \subset \A^{r+1}$ be
the affine cone over $X$ and let $\underline{o}$ be its vertex. Then, the following conditions
are equivalent:
\begin{itemize}
\item[(i)] $X$ is pCM;
\item[(ii)] ${\mathcal C}(X)$ is CM;
\item[(iii)] ${\mathcal C}(X)$ is CM at $\underline{o}$ (cf.\ Remark \ref{rem:0}).
\end{itemize}
\end{proposition}
\begin{proof}
Take $\Gamma(X)$ as in \eqref{eq:gamma}. Then, it coincides with the coordinate ring
of the affine cone ${\mathcal C}(X) \subset \A^{r+1}$. Thus, the claim follows from
Proposition 18.8 and Ex. 19.10 in \cite{Eisenbud}.
\end{proof}

\begin{proposition}\label{prop:74}
Let $X \subset \PR$ be an equidimensional variety. Then:
\begin{itemize}
\item[(i)] $X$ is pCM $\Leftrightarrow$ $X$ is aCM;
\item[(ii)] $X$ is aCM $\Rightarrow$ $X$ is CM.
\end{itemize}
\end{proposition}
\begin{proof} (i) From Proposition \ref{prop:3},
$X$ is pCM iff ${\mathcal C}(X)$ is CM at $\underline{o}$.
From \cite{H}, page 217, Ex. 3.4 (b), this implies that $H^i_{\underline{m}} (\Gamma(X)) = 0,$
for all $i < r-c = n$, where $  c = \codim_{\PR}(X)$, $n = \dim(X)$,
$\underline{m}$ is the maximal, homogeneous
ideal of $\Gamma(X)$ and, as usual, $H^i_{\underline{m}} (-)$ is the
local cohomology (see \cite{H}, Ex. 3.3 (a), page 217).

On the other hand, from \cite{Eisenbud}, Theorem A4.1,
$$H^i_{\underline{m}} (\Gamma(X)) \cong \bigoplus_{j \in \ZZ} H^i (X, \Oc_X(j)).$$

\noindent
(ii) Since $X$ is covered by affine open subsets which are hyperplane
sections of ${\mathcal C}(X)$, the assertion follows from Theorem
\ref{thm:1bis}.

\end{proof}

\begin{theorem}\label{thm:1ter}
Let $X \subset \PR$ be an equidimensional closed subscheme in $\PR$. Then
$X$ is aCM if, and only if, any hyperplane section $H$ of $X$ not containing
any component of $X$ is aCM.
\end{theorem}

\begin{proof}
It directly follows from Theorem \ref{thm:1bis} and Proposition \ref{prop:3}.
\end{proof}

\begin{proposition}\label{cor:pCM}
If a curve $H$, which is a
hyperplane section of a projective surface $S$, is
arithmetically normal then $S$ is aCM (equiv., pCM).
\end{proposition}

\begin{proof}
To prove that $S$ is aCM, one has to prove that
\eqref{eq:H1Ox} holds (see Remark \ref{rem:curve}).
This follows by the proof of Proposition \ref{lem:pn}.
\end{proof}

\begin{corollary}\label{eq:cor}
Let $X \subset \PR$ be a curve. Then:

\noindent
(i) $X$ is projectively normal $\Rightarrow$ $X$ is pCM.

\noindent
(ii) If, furthermore, $X$ is assumed to be
smooth, the implication in (i) is an equivalence.
\end{corollary}
\begin{proof}
By Proposition \ref{prop:74}, $X$ is pCM iff is aCM. On the other hand, since
$X$ is a curve, by Remark \ref{rem:curve} $X$ aCM is equivalent to $X$ arithmetically normal.
Therefore, the curve $X$ is in particular pCM if, and only if, it is arithmetically
normal. Only if $X$ is also smooth, then $X$ is projectively normal, as it follows
from Proposition \ref{prop:1}.
\end{proof}

\begin{proposition}\label{prop:freeacm}
Let $X \subset \PR$ be a projective, equidimensional variety
s.t. $\codim_{\PR}(X) = c$. Then, the following conditions are equivalent:
\begin{itemize}
\item[(i)] $X$ is pCM;
\item[(ii)] the {\em projective dimension} of $\Gamma(X)$ is equal to $c =\codim_{\PR}(X) $.
In other words, there is a minimal graded free resolution of $\Gamma(X)$,
\begin{equation}\label{eq:freeres}
 0 \to F_c \to F_{c-1} \to \cdots \to F_1 \to S \to \Gamma(X) \to 0
\end{equation}where $F_i$ is a free $S$-module, for $ 1 \leq i \leq c$.
\end{itemize}
\end{proposition}

\begin{proof}

From Proposition \ref{prop:3}, $X$ is pCM if, and only if,
$\Gamma(X)$ is CM. Let $S = \CC[x_0, \ldots , x_r]$ be the homogeneous graded polynomial ring
which is a finitely generated algebra over $\CC$. Let $M$ be the homogeneous maximal
ideal in $S$. Since $\Gamma(X)$ is a finitely generated graded $S$ module of finite projective
dimension, then by the {\em Auslander-Buchsbaum Formula} in the graded case, we have
\begin{equation}\label{eq:ausbus}
\pd_S (\Gamma(X)) = \depth_M(S) - \depth_{M\Gamma(X)}(\Gamma(X)) = r+1 - \depth_{M\Gamma(X)}(\Gamma(X)) ,
\end{equation}where $\pd_S (\Gamma(X))$ is the
projective dimension of $\Gamma(X)$ (cf.\ \cite{Eisenbud}, Ex. 19.8, page 485).

Thus, if $X$ is pCM, then $\Gamma(X)$ is CM and therefore
$\depth_{M\Gamma(X)}(\Gamma(X)) = n+1$ so, by \eqref{eq:ausbus},
$\pd_S(\Gamma(X)) = c = \codim_{\PR}(X).$
Conversely, if $\pd_S(\Gamma(X)) = c$, then $\depth_{M\Gamma(X)}(\Gamma(X)) = n+1$, hence
${\mathcal C}(X)$ is CM at $\underline{o}$. One concludes by Proposition \ref{prop:3}.
\end{proof}

\begin{remark}\label{remCMtype}
Observe that the ranks of the free modules $F_i$, $1 \leq i \leq c$, do not
depend on the minimal free resolution of $\Gamma(X)$.
In particular, the rank of $F_c$ in any minimal free resolution
of $\Gamma(X)$ is an invariant of $\Gamma(X)$
called the {\em Cohen-Macaulay type} of $X$.
\end{remark}

We now focus on ``Gorenstein'' conditions.
First, we recall some standard definitions.

\begin{definition}\label{def:gorring}
Let $(A,m)$ be a local, CM ring with residue field $K$.

If $\dim(A) = 0$, then $A$ is called a
{\em Gorenstein} ring if, and only if,
\begin{equation}\label{eq:gorring}
A \cong Hom_K(A, K).
\end{equation}

If $\dim(A) >0$, then $A$ is called a
{\em Gorenstein} ring if, and only if,
there is a non-zero divisor $x \in A$ s.t. $A/(x)$ is Gorenstein. In this case,
for every non-zero divisor $x \in A$, $A/(x)$ is Gorenstein.
\end{definition}

\begin{remark}\label{rem:gorring}

By using \eqref{eq:dual},
$(A,m)$ is Gorenstein if, and only if, $\omega_A \cong A$
(cf.\ \cite{Eisenbud}, Theorem 21.15). As in in Remark \ref{rem:freerescm},
this is equivalent to saying
that there is a minimal, free resolution of $R$-modules
\begin{equation}\label{eq:dualres}
{\mathcal F} \;\; : \; 0 \to F_c \to F_{c-1} \to \ldots F_1 \to R \to A \to 0,
\end{equation}which is symmetric in the sense that
${\mathcal F} \cong {\mathcal F}^*$. This, in turn, is equivalent to
saying that $F_c \cong R$ (see Corollary 21.16 from \cite{Eisenbud}).

Observe that, if $(A,m)$ is Gorenstein then, for each prime ideal $p$ in $A$,
$(A_p, pA_p)$ is Gorenstein. This follows by the above remark and by the flatness
of localization (see \cite{Eisenbud}, page 66).
\end{remark}

\begin{definition}\label{def:gorringgrad}
Let $K$ be a field.
Let $R$ be a graded, finitely generated $K$-algebra.

If $\dim(R) = 0$, then $R$ is called a {\em Gorenstein} graded
ring if, and only if, there is an integer $\delta$ such that
\begin{equation}\label{eq:gorringgrad}
R(\delta) \cong Hom_K(R,K).
\end{equation}

If $\dim(R) > 0$, then $R$ is called a {\em Gorenstein} graded
ring if, and only if, there is an homogeneous non-zero divisor $x \in R$
such that $R/(x)$ is Gorenstein.

\end{definition}

\begin{remark}\label{rem:gorringgrad}
Let $S = \CC[x_0, \ldots , x_{r+1}]$ be the homogeneous polynomial ring and let
$I$ be a homogeneous ideal, such that $A = S/I$ is a CM ring. Set
$c = \codim_S(A)$.

Then
\begin{equation}\label{eq:dualgrad}
\omega_A:= Ext_S^c (A, S(-r-1))
\end{equation}is called the graded
{\em dual module} of $A$. Then, $A$ is Gorenstein
if, and only if, there is an integer $\delta$
such that $\omega_A \cong A(\delta)$. This is equivalent to saying
that there is a minimal, free graded resolution
\begin{equation}\label{eq:dualgradres}
{\mathcal F} \colon 0 \to F_c \to F_{c-1} \to \cdots \to F_1 \to S \to A \to 0,
\end{equation}which is symmetric in the sense that
${\mathcal F} \cong {\mathcal F}^*$. This, in turn, is equivalent to
saying that there is an integer $\gamma$ such that
$F_c \cong S(\gamma)$ (see the proof of Corollary 21.16 and
\S 21.11 from \cite{Eisenbud}).

Observe also that, if $A$ as above is Gorenstein then,
for every prime ideal $P$ in $A$, $(A_P, PA_P)$ is a
local Gorenstein ring.

\end{remark}

In complete analogy with Definition
\ref{def:4}, we have:

\begin{definition}\label{def:6tris}
A projective
scheme $X$ is {\em Gorenstein} at a point $p \in X$ if
$\Oc_{X,p}$ is a Gorenstein, local ring. $X$ is {\em Gorenstein},
if it is Gorenstein at each point $p \in X$.
\end{definition}

\begin{theorem}\label{thm:1tris}
Let $X \subset \PR$ be an equidimensional projective variety and let $p $ be a point of $X$.
Then, the following
conditions are equivalent:
\begin{itemize}
\item[(i)] $X$ is Gorenstein at $p$;
\item[(ii)] each equidimensional hyperplane section of $X$ through $p$ is Gorenstein.
\end{itemize}
\end{theorem}
\begin{proof}
It directly follows from the definition of local, Gorenstein
rings.
\end{proof}

\begin{definition}\label{def:6bis}
A projective and equidimensional
variety $X\subset \PR$ is called \emph{projectively Gorenstein}
(\emph{pG}, for short) if its homogeneous coordinate ring $\Gamma(X)$
is Gorenstein.
\end{definition}

In complete analogy with Proposition \ref{prop:3}, we have:

\begin{proposition}\label{prop:4}
Let $X \subset \PR$ be an equidimensional variety,
let ${\mathcal C}(X) \subset \A^{r+1}$ be
the affine cone over $X$ and let $\underline{o}$ be its vertex. Then, the following conditions
are equivalent:
\begin{itemize}
\item[(i)] $X$ is pG;
\item[(ii)] ${\mathcal C}(X)$ is Gorenstein;
\item[(iii)] ${\mathcal C}(X)$ is Gorenstein at $\underline{o}$;
\item[(iv)] $X$ is pCM and the dualizing sheaf
\begin{equation}
\label{eq:PG} \omega_X \cong \Oc_X(a), \qquad\text{for some $a\in\ZZ$.}
\end{equation}
\end{itemize}
\end{proposition}

\begin{proof}
We prove the following implications.
\begin{itemize}
\item (i) $\Rightarrow $ (ii): it directly follows from what recalled in Remark
\ref{rem:gorringgrad};
\item (ii) $ \Rightarrow $ (iii): trivial;
\item (iii) $ \Rightarrow $ (i): $X$ is pCM from Proposition \ref{prop:3}. Now, let
$${\mathcal F} \;\; : \; 0 \to F_c \to F_{c-1} \to \ldots F_1 \to S \to \Gamma(X)
\to 0$$be a graded, minimal free resolution of $\Gamma(X)$. By localizing
${\mathcal F}$ at the homogeneous maximal ideal
$M$ one still obtains a minimal free resolution. The assertion follows by Remarks
\ref{rem:gorring} and \ref{rem:gorringgrad}.
\item (iv) $ \Leftrightarrow $ (i): it directly follows from the definition
of dual module (see \eqref{eq:dualgrad}) and the definition of dualizing sheaf.
\end{itemize}
\end{proof}

\begin{remark}\label{rem:90}
Observe that, if $X$ is pG then ${\mathcal C}(X)$ is
Gorenstein, hence $X$ is Gorenstein, since it is an equidimensional
hyperplane section of its affine cone.
\end{remark}

Clearly, by Adjunction Formula and by
Theorem \ref{thm:1ter} and Proposition \ref{prop:4},
if $X$ is pG, then each hyperplane section $H$ of $X$  not containing
any component of $X$ is pG. Conversely:

\begin{proposition}\label{lem:PG}
Let $X$ be an equidimensional, projective variety.
If a equidimensional hyperplane section $H$ of a projective variety $X$ is
pG, then $X$ is pG too.
\end{proposition}

\begin{proof}
It directly follows from Theorem \ref{thm:1tris} and Proposition \ref{prop:4}.
\end{proof}


\begin{thebibliography}{ADSE}

\bibitem{AW} Artin, M., Winters, G., Degenerate fibres and stable reduction of curves,
{\em Topology}, {\bf 10} (1971), 373-383.

\bibitem{Bar} Bardelli F., Lectures on stable curves, in {\em Lectures on Riemann
surfaces - Proc. on the college on Riemann surfaces, Trieste - 1987}, 648--704.
Cornalba, Gomez-Mont, Verjovsky, World Scientific, Singapore, 1988.

\bibitem{BPV}
Barth W., Peters C., and Van de Ven A., {\em Compact Complex
Surfaces}, Ergebnisse der Mathematik, 3. Folge, Band 4, Springer,
Berlin, 1984.


\bibitem{BCK} Batyrev, V. V., Ciocan-Fontanine, I., Kim, B. \&
van Straten, D., Mirror symmetry and toric
degenerations of partial flag manifolds. {\em  Acta Math.},
{\bf 184} (2000), no. 1, 1--39.


\bibitem{BS}
Beltrametti M. C., Sommese A. J., {\em The adjunction theory of
complex projective varieties}, de Gruyter expositions in mathematics, 16, de
Gruyter, Berlin--New York, 1995.


\bibitem{Bert} Bertini, E., {\em Introduzione alla geometria proiettiva
degli iperspazi}, Messina, 1923.

\bibitem{CCFMto} Calabri, A., Ciliberto C., Flamini, F., Miranda, R.,
On the geometric genus of reducible surfaces and degenerations of surfaces to unions of planes,
in {\em The Fano Conference}, Univ.\ Torino, Turin, 2004,  277--312.  


\bibitem{CCFMLincei} Calabri, A., Ciliberto C., Flamini, F., Miranda, R.,
Degenerations of scrolls to union of planes, {\em Rend. Lincei Mat. Appl.},
{\bf 17} (2006), no. 2, 95--123. 

\noindent
{\bf Available on the web}: 
http://www.ems-ph.org/journals/all$_-$issues.php?issn=1120-6330, Vol. 17, Issue 2. 

\bibitem{CCFMk2} Calabri, A., Ciliberto C., Flamini, F., Miranda, R.,
On the $K^2$ of degenerations of surfaces and the Multiple Point Formula,
{\em  Annals of Mathematics}, {\bf 165} (2007), no. 2, 335--395.

\noindent
{\bf Available on the web}: 
http://projecteuclid.org/DPubS?service=UI$\&$version=1.0$\&$verb=Display $\&$handle=euclid.annm/1185985451

\bibitem{CCFMpg} Calabri, A., Ciliberto C., Flamini, F., Miranda, R.,
On the genus of reducible surfaces and degenerations of surfaces,
{\em  Annales de l'Institut Fourier (Grenoble)}, {\bf 57} (2007), no. 2, 491--516.



\noindent
{\bf Available on the web}: 
http://aif.cedram.org/,  Vol. 57 (2), 2007, 491--516.


\bibitem{CCFMnonsp} Calabri, A., Ciliberto C., Flamini, F., Miranda, R.,
Non-special scrolls with general moduli, 
to appear on {\em Rend. Circolo Mat. Palermo}, (2008), 24 pages.

\noindent
{\bf Available on the web}: 
http://arxiv.org/PS$_-$cache/arxiv/pdf/0712/0712.2105v1.pdf

\bibitem{CCFMbrill} Calabri, A., Ciliberto C., Flamini, F., Miranda, R.,
Brill-Noether theory and non-special scrolls, 
to appear on {\em Geometriae Dedicata}, (2008), 16 pages.

\noindent
{\bf Available on the web}: 
http://arxiv.org/PS$_-$cache/arxiv/pdf/0712/0712.2106v1.pdf




\bibitem{CDM} Ciliberto, C., Dumitrescu, O., Miranda, R., Degenerations of the Veronese and applications, to be published on {\em Proceedings of the Conference "Linear systems and subschemes"}, Gent, 2007.

\bibitem{CLM} Ciliberto, C., Lopez, A. F., Miranda, R., Projective degenerations
of $K3$ surfaces, Gaussian maps and Fano threefolds, {\em Inv.
Math.}, {\bf 114} (1993), 641--667.

\bibitem{CMT} Ciliberto, C., Miranda, R., Teicher, M., Pillow degenerations of
$K3$ surfaces, in {\em Application of
Algebraic Geometry to Computation, Physics and Coding Theory}, Nato Science Series II/36,
Ciliberto et al.\ (eds.), , Kluwer Academic Publishers, 2002.

\bibitem{CS} Cohen, D. C., Suciu, A. I., The braid monodromy of plane algebraic curves and
hyperplane arrangements, {\em Comment. Math. Helv.}, {\bf 72} (1997), 285--315.


\bibitem{DelP} Del Pezzo, P. , Sulle superficie dell'$n^{mo}$ ordine immerse
nello spazio a $n$ dimensioni, {\em Rend. Circ. Mat. Palermo}, {\bf 1} (1887),
241--247.

\bibitem{Eisenbud} Eisenbud, D., {\em Commutative Algebra, with a view toward
Algebraic Geometry},
Springer, New York, 1995.

\bibitem{E} Enriques, F., {\em Le superficie algebriche},
Zannichelli, Bologna, 1949.

\bibitem{Frie} Friedman, R., Global smoothings of varieties with normal
crossings, {\em Ann. Math.}, {\bf 118} (1983), 75--114.

\bibitem{FM} Friedman, R., Morrison, D.R., (eds.,)
\emph{The birational geometry of degenerations}, Progress in
Mathematics 29, Birkhauser, Boston, 1982.


\bibitem{Fulton}
Fulton, W., \emph{Introduction to toric varieties}, Annals of
mathematics studies {\bf 131}, Princeton University Press, 1993.







\bibitem{Gaeta} Gaeta, F.,
Nuove ricerche sulle curve sghembe algebriche
di residuale finito e sui gruppi di punti del piano,
An.\ Mat.\ Pura Appl.\ (4) 31, (1950). 1--64.

\bibitem{GonciuL} Gonciulea, N., Lakshmibai, V., Schubert varieties,
toric varieties, and ladder determinantal varieties.
{\em Ann. Inst. Fourier (Grenoble)}, {\bf 47} (1997),
no. 4, 1013--1064.

\bibitem{GMac} Goresky, M., MacPherson, R., {\em Stratified Morse Theory},
Ergebnisse der Mathematik und ihrer Grenzgebiete (3), {\bf 14},
Springer Verlag, Berlin, 1988.




\bibitem{Greco} Greco, S., {\em Normal varieties.} Notes written with the
collaboration of A. Di Sante. Inst. Math. {\bf 4}, Academic Press, London-New
York, 1978.

\bibitem{GH} Griffiths, P., Harris, J., {\em Principles of Algebraic Geometry},
Wiley Classics Library, New York, 1978.

\bibitem{Harris} Harris, J., A bound on the geometric genus of projective
varieties, {\em Ann. Scu. Norm. Pisa}, {\bf 8} (1981), 35--68.

\bibitem{Hart2} Hartshorne, R.,
On the De Rham cohomology of algebraic varieties,
{\em Inst. Hautes Études Sci. Publ. Math.}, {\bf 45} (1975), 5--99.

\bibitem{H}
Hartshorne, R., {\em Algebraic Geometry}, (GTM, No. 52), Springer-Verlag, New
York-Heidelberg, 1977.

\bibitem{Hart}
Hartshorne, R., {\em Families of curves in $\Pp^3$ and Zeuthen's
problem}, Mem.\ Amer.\ Math.\ Soc.\ {\bf 130} (1997), no.\  617.



\bibitem{Iv} Iversen, B.,
Critical points of an algebraic function,
{\em Invent. Math.}, {\bf 12} (1971), 210--224.


\bibitem{Kempf}
Kempf, G., Knudsen, F.F., Mumford, D., and Saint-Donat, B.,
Toroidal embeddings. I., Lecture Notes in Mathematics 339,
Springer-Verlag, Berlin-New York, 1973.

\bibitem{Kollar} Koll\'ar, J., Toward moduli of singular varieties,
{\em Compositio Mathematica}, {\bf 56} (1985), 369--398.

\bibitem{KolSB} Koll\'ar, J., Shepherd-Barron, N.I.,
Threefolds and deformations of surface singularities, {\em Invent.
math.}, {\bf 91} (1988), 299--338.



\bibitem{MacRae} MacRae, R.E.,
On an application of Fitting invariants,
{\em J. Algebra}, {\bf 2} (1965), 153--169.

\bibitem{Matsumura} Matsumura, H., {\em Commutative Algebra}, Reading Massachussets, B.
Cummings Publishing Program, 1980.

\bibitem{Ran} Melliez, F., Ranestad, K.,
Degenerations of $(1,7)$-polarized abelian surfaces,
{\em Math. Scand.}, {\bf 97} (2005), 161--187.

\bibitem{Migliore}
Migliore, J.C., {\em An introduction to deficiency modules and
liaison theory for subschemes of projective space},  Lecture Notes Series, {\bf 24}.
Seoul National University,
Research Institute of Mathematics,
Global Analysis Research Center, Seoul, 1994


\bibitem{Miyaoka} Miyaoka, Y.,
On the Chern numbers of surfaces of general type,
{\em Invent. Math.}, {\bf 42} (1977), 255--237.

\bibitem{Moi} Moishezon, B., Stable branch curves and braid monodromies,
in {\em Algebraic geometry (Chicago, Ill., 1980)}, Lecture Notes in Math.,
{\bf 862} (1981), Springer , 107--192.

\bibitem{MoiTei} Moishezon, B., Teicher, M., Braid group techniques in
complex geometry III: Projective degeneration of $V_3$, in {\em
Classification of Algebraic Varieties, Contemporary Mathematics}
{\bf 162}, AMS (1994), 313--332.

\bibitem{Morr} Morrison, D.R., The Clemens-Schmid exact sequence and
applications, in {\em Topics in Trascendental Algebraic Geometry},
{\em Ann. of Math. Studies}, {\bf 106} (1984), 101--119.

\bibitem{Mum} Mumford D., A remark on the paper of M. Schlessinger, {\em
Rice Univ. Studies}, {\bf 59/1} (1973), 113--118.

\bibitem{Nami} Namikawa, Y., Toroidal degeneration of abelian
varieties. II., {\em Math. Ann.}, {\bf 245} (1979), no. 2, 117--150.


\bibitem{Pers} Persson, U., {\em On degeneration of algebraic surfaces},
Memoirs of the American Mathematical Society, {\bf 189}, AMS,
Providence, 1977.

\bibitem{Pink1} Pinkham, H.C., {\em Deformations of algebraic varieties
with $\mathbb{G}_m$ action}, Ph.D. thesis, Harvard University,
Cambridge, Massachussets, 1974.

\bibitem{Pink2} Pinkham, H.C., Deformation of cones with negative grading,
{\em J. Algebra}, {\bf 30} (1974), 92--102.

\bibitem{Sally} Sally, J.D., Tangent cones at Gorenstein
singularities, {\em Compositio Math.}, {\bf 40} (1980), no.2,
167-175.




\bibitem{Sev} Severi, F., {\em Vorlesungen \"uber algebraische Geometrie},
vol. 1, Teubner, Leipzig, 1921.

\bibitem{S} Sernesi, E., {\em Topics on families of projective schemes},
Queen's Papers in Pure and Applied
Mathematics, n. {\bf 73},
Queen's University Press, Canada, 1986.

\bibitem{Teicher} Teicher, M., Hirzebruch surfaces: degenerations, related
braid monodromy, Galois covers. in
{\em Algebraic geometry: Hirzebruch 70 (Warsaw, 1998)},  {\em Contemp. Math.},
{\bf 241} (1999), 305--325.

\bibitem{X} Xamb\'o, S., On projective varieties of minimal degree,
{\em Collectanea Math.}, {\bf 32} (1981), 149--163.

\bibitem{Yau} Yau, S.-T.,
Calabi's conjecture and some new results in algebraic geometry,
Proc.\ Natl.\ Acad.\ Sci.\ USA 74 (1977), 1789--1799.

\bibitem{Za1} Zappa, G., Caratterizzazione delle curve di diramazione delle
rigate e spezzamento di queste in sistemi di piani, {\em Rend.
Sem. Mat. Univ. Padova}, {\bf 13} (1942), 41--56.

\bibitem{Za1b} Zappa, G., Su alcuni contributi
alla conoscenza della struttura topologica delle superficie
algebriche, dati dal metodo dello spezzamento in sistemi di piani,
{\em Acta Pont.\ Accad.\ Sci.}, {\bf 7} (1943), 4--8.

\bibitem{Za1c} Zappa, G., Applicazione della teoria delle matrici di
Veblen e di Poincar\'e allo studio delle superficie spezzate in
sistemi di piani, {\em Acta Pont.\ Accad.\ Sci.}, {\bf 7} (1943),
21--25.

\bibitem{Za2} Zappa, G., Sulla degenerazione delle superficie
algebriche in sistemi di piani distinti, con applicazioni allo
studio delle rigate, {\em Atti R.\ Accad.\ d'Italia, Mem.\ Cl.\
Sci.\ FF., MM.\ e NN.}, {\bf 13} (2) (1943), 989--1021.

\bibitem{Za2b} Zappa, G., Invarianti numerici
d'una superficie algebrica e deduzione della formula di
Picard-Alexander col metodo dello spezzamento in piani, {\em
Rend.\ di Mat.\ Roma}, {\bf 5} (5) (1946), 121--130.

\bibitem{Za3} Zappa, G., Sopra una probabile disuguaglianza tra i caratteri
invariantivi di una superficie algebrica, {\em Rend. Mat. e
Appl.}, {\bf 14} (1955), 1--10.

\bibitem{Za3b} Zappa, G., Alla ricerca di nuovi significati topologici dei
generi geometrico ed aritmetico di una superficie algebrica, {\em
Ann.\ Mat.\ Pura Appl.}, {\bf 30} (4) (1949), 123--146.

\bibitem{Zar0} Zariski, O., Normal varieties and birational correspondences, {\em
Bull. \ Amer\ Math. \ Soc.}, {\bf 48} (1942), 402--413.

\bibitem{Zar} Zariski, O., {\em Algebraic surfaces, Second supplemented
editions, with Appendices S.S. Abhyankar, J. Lipman and D.
Mumford}, Springer-Verlag, Berlin, 1971.

\end{thebibliography}
\end{document}